\documentclass{article}
\usepackage{fancyhdr} 
\usepackage{color} 
\usepackage{hyperref} 
\usepackage{graphicx} 

\usepackage{pstricks}
\usepackage{amssymb}

\usepackage{amsthm,amsmath}

\usepackage{amsmath}
\usepackage{amsthm}
\usepackage{amsbsy}
\usepackage{bbm}
\usepackage{bm}
\usepackage{amssymb}
\usepackage{amsfonts}
\usepackage{cancel}
\usepackage{centernot}
\usepackage[colorinlistoftodos]{todonotes}
\usepackage[dvips]{lscape}
\usepackage{enumerate}
\usepackage{graphicx}
\usepackage{tikz}
\usepackage{hyperref}
\usepackage{verbatim}
\usepackage[lmargin=3cm,rmargin=3cm]{geometry}
\usepackage{changepage}

\definecolor{aleacolor}{rgb}{0.16,0.59,0.78}

\hypersetup{
breaklinks,
colorlinks=true,
linkcolor=aleacolor,
urlcolor=aleacolor,
citecolor=aleacolor}

\theoremstyle{plain}
\newtheorem{theorem}{Theorem}[section]

\theoremstyle{definition}
\newtheorem{definition}[theorem]{Definition}
\theoremstyle{remark}

\makeatletter \@addtoreset{equation}{section} \makeatother

\let\d\delta
\let\eps\varepsilon
\let\n\noindent
\newcommand{\ov}{   \overline{\Pi}   }

\let\b\begin
\let\e\end
\let\f\frac
\let\bb\mathbb
\let\cal\mathcal
\let\bbm\mathbbm
\let\l\left
\let\r\right
\let\p\paragraph
\newcommand{\pp}{\mathbb{P}}
\newcommand{\oo}{\mathcal{O}}
\newcommand{\D}{\Delta_{1}}
\newcommand{\ks}{\cite{ks13}}

\begin{document}

\title{Transience and Recurrence of Markov Processes with Constrained Local Time}

\author{Adam Barker}
\date{}





\maketitle

\begin{abstract}
 We study Markov processes conditioned so that their local time must grow slower than a prescribed function. Building upon recent work on Brownian motion with constrained local time in \cite{bb11,ks13}, we study transience and recurrence  for a broad class of Markov processes. 

    In order to understand the   local time, we determine the distribution of a non-decreasing L\'evy process (the inverse local time) conditioned to remain above a given level which varies in time. We study a time-dependent     region, in contrast to previous works in which a process is conditioned to remain in a fixed region  (e.g.\   \cite{dw15,g09}), so we must  study boundary crossing probabilities for a family of curves, and thus obtain uniform asymptotics for such a family.
     
 Main results include necessary and sufficient conditions for transience or recurrence of the conditioned Markov process.   We will   explicitly determine the distribution of the  inverse local time for the conditioned process, and in the transient case, we explicitly determine the law of the conditioned Markov process. In the recurrent case, we characterise the ``entropic repulsion envelope'' via necessary and sufficient conditions. 
\end{abstract}

\section{Introduction}\label{secLabeled}

We study the asymptotic behaviour of a Markov process whose local time is constrained to grow slower than $f$, an     increasing function.     
 %
 %
The (right-continuous) inverse of the local time process is a subordinator (a non-decreasing L\'evy process)
, so our study of the behaviour of the local time process  is effectively equivalent to studying a  subordinator conditioned to grow faster than the inverse function $f^{-1}$.

This work is hence related to a number of works on stochastic processes conditioned to remain in a certain fixed region, such as cones in \cite{dw15,g09}, 
 and Weyl chambers in \cite{dw10,ks10}.   We highlight the fact that our subordinator is conditioned to remain in a region which varies in time, whereas the aforementioned works consider fixed regions, as appears to be the case for all works prior~to~\cite{ks13}.

We emphasise that in constraining the local time of a Markov process, the extent to which our constraint affects the process varies over time, depending on the past behaviour of the process.     This   is a ``weak'' constraint, in constrast to ``strong'' constraints such as conditioning a process to avoid a point, where the constraint does not change (see  e.g.\ \cite{b92,bbc03,c96}).    When our conditioned Markov process is recurrent, our constraint    varies over all time, whereas when our conditioned Markov process is transient, the constraint varies for only a finite window of time. So our   results, especially   in the recurrent case,  offer a significant contrast to many  prior works with a ``strong'' constraint.

Many works, e.g.\   \cite{aks12,bl16,dw14,kl13,l13,m16,pw01}, 
   consider a time-dependent region, and the passage time out of this region is studied.   In this paper, we study the  boundary crossing probability  for  a   family of curves  and study the time at which this crossing occurs (the function $f^{-1}$ forms our boundary of interest). We study these asymptotics,  uniformly, among a family of curves   in Lemma \mbox{\ref{lemma2}} (in contrast to prior non-uniform asymptotic results),   and   consider the deeper problem of determining the law of a subordinator conditioned to remain in this time-dependent region.  
%
Our studies are also similar, in spirit, to  various other works on Brownian motion \cite{ks17,rvy06,p75}, 
   L\'evy processes \cite{b93,bbc03,c96,p13,yyy09},   
    and more general diffusions \cite{c18,rvy07,sv09} with restricted path behaviour.

  Specifically, this work is motivated by previous works on Brownian motion with constrained local time, such as \cite{bb11}, in which a 1-dimensional Brownian motion is conditioned so the local time at zero, $(L_t)_{t\geq0}$, satisfies $L_t\leq f(t)$ for all $t\geq0$, for a given function $f$, and a sufficient condition for transience of the conditioned Brownian motion is found. 
  
  In \cite{ks13}, it is shown that the condition  is necessary and sufficient for transience of the conditioned Brownian motion, and the law of the conditioned inverse local time process is explicitly determined in both the transient and recurrent cases. In the transient case, an explicit formulation for the conditioned Brownian motion is found, and in the recurrent case the ``entropic repulsion envelope'' is found.   
  
  This paper builds upon \cite{ks13} in particular, providing analogous results for a much broader class of processes than Brownian motion, with some mild regularity conditions.  It was conjectured in \cite[Remark 9]{ks13} that such analogous results hold when the inverse local time process has L\'evy measure with regularly varying tail, which we confirm in this paper. We   extend beyond     this conjecture   by including a much  more general setting, see Assumption~\mbox{\ref{case2}}.

Constraining local time from above  imposes weak repellence on a Markov process.   Many such processes   are studied in works related to polymer physics, see e.g.\ \cite{bhp18,chp12,hk01,hk01a,w80}.   
%
%
 %
%
Particularly important is the  transition
between a localised phase, where the polymer remains close to a point, and a
delocalised phase, where it moves away.

The goal is often to understand
when the transition occurs as underlying model parameters vary, as in e.g.\ \cite{b02,bhp18,chp12,gk97,hkw11}. 
 This motivates our study of   transience and recurrence of  Markov   processeses, transience and recurrence corresponding to  delocalised and localised phases, respectively. 
 
 Now we  provide a brief exposition of the main result, before introducing some key definitions.


\p{Main Result} Starting with a recurrent Markov process, we constrain its local time $(L_t)_{t\geq0}$ so that $L_t\leq f(t)$ for all $t$. The following necessary and sufficient condition tells us if the constraint is strong enough to change the process to become transient: it is transient   if 
\b{equation}
\label{criterion}
 \int_1^\infty     f(x)  \Pi(dx)     <\infty, 
\e{equation}
and the process remains recurrent   otherwise, where $\Pi(dx)$ denotes the L\'evy measure of the inverse local time subordinator. 
Our criterion $\l(\ref{criterion}\r)$ can also be understood in terms of the rate of growth of the inverse local time $(X_s)_{s\geq0}$ as $s\to\infty$, since it is known \cite[Theorem III.13]{b98} that
\[
\int_1^\infty     f(x)  \Pi(dx)     <\infty   \iff      \lim_{s\to\infty} \frac{X_s}{f^{-1}(s)} =0, \textup{ almost surely. }
\]
So the boundary choice of $f$, at which the conditioned process changes from recurrent to transient, coincides with the boundary at which $X_s$  grows to infinity faster or slower than $f^{-1}(s)$.

  The remainder of the paper is structured as follows: Section \mbox{\ref{keydefns}} provides key definitions; Section \hbox{\ref{mainresults}} outlines the statements of  the main results and the conditions under which they hold, including the necessary and sufficient conditions for transience/recurrence, the   distribution of the conditioned process, and the characterisation of the entropic repulsion envelope; Section  \hbox{\ref{mainresultsproofs}} contains the proofs of the main results; Sections    \hbox{\ref{lemma2proof}} and \hbox{\ref{l6proof}} contain the proofs of 2 key lemmas required for the main results; Section  \hbox{\ref{lemmasproofs}} contains the proofs of the remaining auxiliary lemmas. 

\section{Key Definitions}
\label{keydefns} We shall provide some definitions,  following   conventions   of \cite[Chapter IV]{b98}. 
 \begin{definition}
 \label{markovprocessdefn}
 A Markov process $(M_t)_{t\geq0}$, is a $\bb{R}^d$-valued stochastic process   such that for each (almost surely) finite stopping time $T$, under the conditional law $ \bb{P}\l(   \cdot |   M_{T} =x   \r)   =   \bb{P}_{x}\l( \cdot  \r)$, the shifted process  $(M_{s+T})_{s\geq0}$ is independent of $\cal{F}_{T}$
 and has  the same as the law, $\bb{P}_x$, as the process $M$ started from $x$. 
  Moreover, we impose that $M$ has right-continuous sample paths, $M_0=0$, and that the origin is \textit{regular} and \textit{instantaneous}. 
   Regular means that for each (almost surely) finite stopping time $T$, if $M_T=0$, then $\inf \l\{   t>T :  M_t=0      \r\} =T$ almost surely.
 Instantaneous means that for each (almost surely) finite stopping time $T$,  if $M_T=0$, then $\inf \l\{   t>T :  M_t \neq 0      \r\} =T$ almost surely.
 \end{definition}
 \begin{definition} For a Markov process, and for an arbitrary choice of $c\in(0,\infty)$, let $l_1(x)$ denote the length of the first excursion interval (away from zero) of length $l>x>0$, and define 
  \[
  P(a):= \begin{cases} 1/\bb{P}(l_1(a)>c), \quad \hspace{0.5pt} 0<a\leq c, \\ \bb{P}(l_1(c)>a), \qquad  a> c.
  \end{cases}
  \]
  Let $g_n(a)$ be the start time of the $n$th excursion of length  $l>a>0$, and write $N_a(t):= \sup\{ n\in \bb{N} : g_n(a)<t  \}$. Then the local time process at zero, $(L_t)_{t\geq0}$,   is defined by $L_t:= \lim_{a\to0} N_a(t)/P(a)$.
 \end{definition}

  %
  %
    A subordinator is defined to be a non-decreasing real-valued stochastic process with stationary independent increments, started from 0.  
 The right continuous inverse local time, defined by $X_t:= \inf\{ s>0 :  L_s > t   \}$,   is a subordinator.     The jumps of $(X_t)_{t\geq0}$ correspond to excursions of $(M_t)_{t\geq0}$ away from zero.

  The Laplace exponent $\phi$ of a subordinator $X$ is defined by   $ e^{-\phi(\lambda)} =   \mathbb{E}[  e^{- \lambda X_1}  ] $,    $\lambda\geq0$. By the   L\'evy-Khintchine formula \cite[p72]{b98},   $\phi$ can   be written 
 \begin{equation*}
 \label{lk}   \phi(\lambda) = \ \text{d} \lambda +  \int_0^\infty  (1- e^{-\lambda x} ) \Pi(dx) , 
 \end{equation*}
 where d is the linear drift, and $\Pi$ is the L\'evy measure, which determines the size and rate of the jumps   of $X$, and satisfies $\int_0^\infty (1\wedge x)\Pi(dx)<\infty$.  
We refer to  \cite{b98} 
 for  background on subordinators.

Next, we define some important classes of functions with which we shall work. 

   \b{definition}[Regular Variation and Related Properties]  \  \\ 
 \begin{enumerate}[(i)]  
   \item  \label{rv}  
         A   function $h:\bb{R}\to \bb{R}$ is regularly varying at $\infty$ with index $\alpha\in\bb{R}$ if    for all $\lambda>0$,  we have $\lim_{t\to\infty} h(\lambda t)/ h(t)   = \lambda^{\alpha} $. We refer to  \mbox{\cite{bgt89}} for background on regular variation.      
  \\
  
\item \label{sv}  A   function $L:\bb{R}\to \bb{R}$ is slowly varying at $\infty$ if $\lim_{t\to\infty} L(\lambda t)/ L(t)   =1$ for each $\lambda>0$.  A  function $h$, regularly varying at $\infty$ of index $\alpha$, can always be written as $h(x) = x^\alpha L(x)$, where $L$ is slowly varying at $\infty$, see \mbox{\cite{bgt89}}.
\\

  \item \label{matus} The lower index, $\beta(h)$, of a function $h:\bb{R}\to \bb{R}$ is  the supremum of $\beta\in\mathbb{R}$ for which there exists $     C  > 0$  so that     for all  $    \Lambda  > 1$,    $  h(\lambda x)/  h(x)         \geq      (1+o(1))  C  \lambda^\beta$,  uniformly in   $ \lambda\in [1,\Lambda] $, as $x\to\infty$,~see~\cite[p68]{bgt89}. 
\\
 
\item \label{crv}    A function $h$ is   CRV at $\infty$ if    $   \lim_{\lambda\to1} \lim_{t\to \infty }   h(\lambda t)/   h(t)       =1$.  The class of CRV functions  lies between ``extended regularly varying'' functions and   $\cal{O}$-regularly varying functions. See \cite{d98} for   details.  
\\ 
 
 \item   \label{orv}      A function $h:\bb{R}\to \bb{R}$ is    $\cal{O}$-regularly varying at $\infty$ if   for each $\lambda>0$,    both  $   \limsup_{ t\to \infty }        h(\lambda t) /   h(t)     <\infty$ and   $   \liminf_{ t\to \infty }   \   h(\lambda t) /   h(t)        >0$.   See \cite{d98} for further details.          
\end{enumerate}
\e{definition}

\section{Statements of Main Results} \label{mainresults}    We aim to constrain the local time so that $L_t \leq f(t)$ for all $t\geq0$, where $f:[0,\infty)\to(0,\infty)$ is  increasing, $f(0)\in(0,1)$,     and $\lim_{t\to\infty}f(t)=\infty$. This work concerns  the behaviour  of our process as $t\to\infty$,   which is unaffected  by the  condition on $f(0)$. Before stating our main results, we   define the regularity conditions under which these~results~hold. 
%
%
%
%
 \subsection{Regularity Conditions}  We shall  impose regularity  conditions on  the function  $f$, its   inverse function $g:=f^{-1}$ (extended so that for $x\in[0,f(0))$, $g(x)=0$), and the tail    $\ov(x):=\Pi(x,\infty)$ in two main cases of interest. Our conditions are imposed on the inverse local time subordinator rather than directly on the Markov process. 
  Now let us define the  main cases of interest for our results:
\b{assumption}[Case (i)]  \label{case1}     We impose on our subordinator that the drift is zero, and the  tail    $\ov(x)=\Pi(x,\infty)$ is regularly varying at $\infty$ with index $-\alpha\in(-1,0)$, so  $\ov(x)=x^{-\alpha} L(x)$ for $L$ slowly varying at $\infty$. We further impose   there exist $B,N>0$  such that the function  $x\mapsto  x^N L(x) $ is non-decreasing on   $(B,\infty)$.

\n We impose that   $f(0)\in(0,1)$,    $\lim_{t\to\infty} f(t) = \infty$, $f$ is  differentiable,  $tf'(t)  \ov(t) $ decreases to 0 as $t\to\infty$ (so $f$ is increasing), the inverse $g:=f^{-1}$ satisfies $\lim_{t\to\infty} g(t+\eps)/g(t)=1$ for all $\eps>0$,  and there exists some value $\beta>  (1+2\alpha)/(2 \alpha + \alpha^2)>1 $ such that 
  \b{equation} \label{condt}
  \lim_{t\to\infty}    t \ov\l(    \f{g(t)}{  \log(t)^\beta  } \r)    =0.
  \e{equation}
   \e{assumption}
  
\b{remark}
Case (\hbox{\hyperref[case1]{i}}) includes stable subordinators, and subordinators whose L\'evy measure has similarly well-behaved tail asymptotics. Thus the set of Markov processes corresponding to  case    (\hbox{\hyperref[case1]{i}})  includes   Bessel processes, stable L\'evy processes of index $\alpha\in(1,2)$, and other Markov processes with similarly well-behaved asymptotics.   Case (\hbox{\hyperref[case2]{ii}}) corresponds to a much broader class~of~processes. 
\e{remark}



   \b{assumption}[Case (ia)]  \label{case1a} Under the assumptions of case (\hbox{\hyperref[case1]{i}}),     define   ``case  (ia)'' by     imposing   $f$, $f'$ are \mbox{\hyperref[orv]{$\cal{O}$-regularly varying}} at $\infty$,         the densities    
    $ f_t(x) dx:=  \pp(X_t\in dx)  $ and  $u(x)dx := \Pi(dx)$ exist, $u$ has bounded increase and bounded decrease (see \cite[p71]{bgt89} for precise definitions),  and          there exist  constants $a,x_0\in(0,\infty)$,    such that for all $t\in(0,\infty)$ and  $x\geq  g(t) +  x_0$, where~$g=f^{-1}$,        
        \b{equation}   \label{l6cond}         f_t(x)    \leq      a  t  u(x).
    \e{equation}  
  %
  \e{assumption}

    \b{remark} If $\ov$ is \mbox{\hyperref[rv]{regularly varying}} at $\infty$ and the density $f_t$   exists, then     $\l(\ref{l6cond}\r)$ 
    holds  for each fixed $t$ and $x>x(t)$, where $x(t)$ may depend on $t$            (see e.g.\     \cite[Theorem 1]{y94}). Here we further impose a bound on    $x(t)$, so that $\l(\ref{l6cond}\r)$ holds  uniformly    among sufficiently many  $x$ and $t$ for us to prove Theorem \mbox{\ref{entrrep}}.  For a stable subordinator of index $\alpha\in (0,1)$,  the density $f_t$ exists (see \cite[p227]{b98}) and $\l(\ref{l6cond}\r)$  
    holds  (see Corollary \mbox{\ref{entrrepstable}}), so  case (\hbox{\hyperref[case1a]{ia}}) includes stable subordinators.
    \e{remark}

\b{assumption}[Case (ii)]  \label{case2}  We impose on our subordinator that the drift is zero, and the   tail function    $\ov(x)=\Pi(x,\infty)$ is   CRV   at $\infty$, with lower   index  $\beta(\ov)>-1$. 

\n  We impose  that  $f(0)  \in (0,1)$,    for  $f$   increasing,  and that there exists $\eps> 0$ such that for $g := f^{-1}$,  
   \b{equation}    \label{cond2}
    \lim_{t\to\infty}    t ^{1+\eps} \ov\l(   g(t)  \r)    =0.
    \e{equation}
  \e{assumption}

 
       \b{remark} In Assumption  \mbox{\ref{case2}}, we impose    $\beta(\ov)>-1$, which is equivalent to imposing that the function $ \int_0^x   \ov(y) dy$ has positive increase as $x\to\infty$. This has many equivalent formulations   \cite[Ex. III.7]{b98}, \cite[Section 2.1]{bgt89}, and appears naturally in a range of contexts \cite[p2]{b18}, \cite[p87]{b98}.  
 \e{remark}

 \n  Now let us introduce some notation required to formulate our results. 
      Recall that   $f:[0,\infty)\to(0,\infty)$ is  increasing,   $f(0)\in(0,1)$, and $g:=f^{-1}$ is the inverse of $f$, where  we take $g(x)=0$ for $x\in[0,f(0))$. The event $ \cal{O}_{u}$ corresponds to bounding the inverse local time until time $u$ (or equivalently, bounding the local time until time $g(u)$). We will study the asymptotics of $\pp(\cal{O}_u)$ as $u\to\infty$, and those of the integral $\Phi(s)$ of this probability.
 \b{align}
  \label{o}
  \cal{O}_{u} \    &:=    \{   X_s \geq   g(s),   \forall \  0\leq s\leq u      \}  ,    \phantom{\Big(}
  \\ 
\label{Phi}
  \Phi(s) \   &:= \int_0^s \pp(  \cal{O}_u )du     .              \phantom{\Big(}
   \intertext{We also study the event $ \cal{O}_{u}$ for the process  $X^{(0,a)}$ with truncated L\'evy measure $\Pi(dx)\mathbbm{1}_{\{x\in(0,a)\}}$,  }
 \label{truncoo} 
    \cal{O}_{u,X^{(0,a)} }   &:= \l\{     X_s^{(0,a)}     \geq g(s), \forall \   0 \leq s \leq u   \r\}.
   %
   \intertext{The time of our subordinator's first jump of size larger than $x>0$, or in the interval $(a,b)$ for $b>a>0$, are respectively denoted by}
 \label{D}
  \Delta_1^x \ \ \   \  &:=\inf\l\{ t\geq 0 : X_t - X_{t-} >x \r\},       \phantom{\Big(}
  \\
  \label{D()}
  \Delta_1^{(a,b)}  \  &:=\inf\l\{ t\geq 0 : X_t - X_{t-} \in(a,b) \r\}.
     \intertext{In Proposition \mbox{\ref{transientcorollary}} and Proposition \mbox{\ref{lemma5}}, we determine that   $I(f)<\infty$ is a necessary and sufficient condition for transience of the conditioned process,~where}
     \label{}
   I(f)    &:=    \int_1^\infty   f(x)    \Pi(dx). \phantom{\Big(}
   \intertext{\b{remark} The necessary and sufficient condition    $I(f)  <    \infty$ arises naturally in a number of contexts, including rate of growth of subordinators \cite[Theorem III.13]{b98} and spectrally negative L\'evy processes \cite[Theorem 3]{p08}.
   \e{remark}
   }
    \intertext{For $h<t$, by the stationary independent increments property, $ \{X_t>g(t) | X_h =y\} $ is effectively equivalent to $\{X_{t-h} > g(t) -y\} = \{X_{t-h} > g((t-h)+h) -y\}$, in terms of probability.   The new   boundary  for $X$ to stay above is given by $g_{y}^{h}(\cdot)$, with $ \cal{O}_u^{g_{y}^{h}} $, $   \Phi_y^h(s)$    corresponding to $ \cal{O}_u$, $\Phi(s)$, where} 
  \label{gyh}
   g_{y}^{h}(t)  &:= g(t+h)-y   ,     \phantom{\Big(}
     \\ 
       \label{ogyh}
  \cal{O}_{u}^{g_{y}^{h}}   &  :=   \{   X_s \geq   g_{y}^{h}(s),   \forall \  0\leq s\leq u      \}      ,            \phantom{\Big(}
  \\
  \label{phiyh}
  \Phi_y^h(s)    &:=\int_0^s   \pp(\cal{O}_{u}^{g_{y}^{h}})du     .            \phantom{\Big(}
     \end{align}
     
     \n    The functions $\rho(\cdot)$, $\rho_y^h(\cdot)$  are   error terms  in the   upcoming ODEs $\l(\ref{ODE1}\r)$~and~$\l(\ref{ODE2}\r)$.    \b{align}     
       \label{rho}
  \rho(t) &:=\f{  \pp\l( \cal{O}_t\r) }{ \Phi(t)}    - \ov(g(t))     ,        \phantom{\Big(}
\\
 \label{rhoyh}
  \rho_y^h(t) &:=   \f{  \pp( \cal{O}_{t}^{g_{y}^{h}}) }{ \Phi_y^h(t)}    - \ov(g_{y}^{h}(t)).       \phantom{\Big(}
\end{align}

  \n     The law of our conditioned process will be found by taking limits. Recalling the notation (\ref{o}) and (\ref{ogyh}), for the measure $\mathbb{Q}(\cdot):=\lim_{t\to\infty} \mathbb{P}(\cdot | \cal{O}_t) $, for all $\cal{B}_h \subseteq \cal{O}_h, \cal{B}_h \in \cal{F}_h$, where $(\cal{F}_u)_{u\geq0}$ is the natural filtration of $X$,
 \begin{align} 
 \mathbb{Q}(X_h\in dy ; \cal{B}_h )  : &=  \lim_{t\to\infty}  \mathbb{P}\l(    X_h \in dy ; \cal{B}_h      |   \cal{O}_t     \r) \nonumber
\\
 &=   \lim_{t\to\infty}    \f{  \mathbb{P}\l(   \cal{O}_t    |  X_h \in dy ; \cal{B}_h       \r)     \mathbb{P}\l( X_h \in dy ; \cal{B}_h   \r)   }{  \mathbb{P}(\cal{O}_t)    }    \nonumber
\\
 &=     \mathbb{P}\l( X_h \in dy ; \cal{B}_h    \r)     \lim_{t\to\infty}    \f{  \mathbb{P}(   \cal{O}_{t-h}^{g_{y}^{h}}        )      }{  \mathbb{P}(\cal{O}_t)    }   .      \label{doobh}
\end{align} 
\n  We shall see in Theorems \mbox{\ref{transientthm}} and  \mbox{\ref{recurrentthm}}  that the limit in $\l(\ref{doobh}\r)$ exists and is finite. In order to understand the   behaviour of  $X$ under $\mathbb{Q}$, we   study  the probabilities $ \mathbb{P}(     \cal{O}_{t}^{g_{y}^{h}}       )$ and $ \mathbb{P}(   \cal{O}_t       )$ as $t\to\infty$.   Corollary  \hbox{\ref{lemma1}} relates the asymptotics of $\pp(\cal{O}_t)$ to  $\Phi(t)$,  and   $ \mathbb{P}(     \cal{O}_{t}^{g_{y}^{h}}       )$ to $\Phi_y^h(t)$.
 We   obtain the~ODEs 
\b{align}
\label{ODE1}
      \pp(\cal{O}_t )   &=   \f{d}{dt} \Phi(t)   =    \l(   \ov(g(t)) + \rho(t) \r)  \Phi(t) ,
      \\
      \label{ODE2}
      \pp(\cal{O}_{t}^{g_{y}^{h}} )   &=  \f{d}{dt}   \Phi_{y}^h (t)   =    \l(   \ov(g_{y}^{h}(t)) + \rho_y^h(t) \r)  \Phi_y^h(t) .
\intertext{These ODEs are easily solved, yielding (for $B>0$ as in Assumption~\mbox{\ref{case1}}) the expressions}
\label{soln1}
\Phi(t) &=     \Phi(1) \exp\l(       \int_1^t  \ov(g(s))ds +        \int_1^t  \rho(s)ds        \r) ,
\\
\label{soln2}
\Phi_y^h(t) &=     \Phi_y^h(t_0(y)) \exp\l(       \int_{t_0(y)}^t  \ov(g_{y}^{h}(s))ds +        \int_{t_0(y)}^t  \rho_y^h(s)ds        \r) ,
\\ 
     \label{t0y}
    t_0(y) :\hspace{-2.5pt}&=  f(Ay)\vee  f(1 + 2/A),   \qquad A>3\vee (B-1).      \phantom{\Big(}   
    \end{align}
The error terms  $\rho(\cdot)$ and  $\rho_y^h(\cdot)$ are later shown to be integrable (the latter uniformly in $y$ and $h$), which is key for determining the distribution of our conditioned process. The required bound for  $\rho(\cdot)$ is given  in Remark \mbox{\ref{yh0}}, and we provide a uniform bound for $\rho_y^h(\cdot)$   in Lemma \hbox{\ref{lemma2}},  proven in Section    \hbox{\ref{lemma2proof}}.

\subsection{Results in the \textit{I}(\textit{f})  $<\infty$ Case} \label{transientsection}   
We shall see  in Theorem \mbox{\ref{Mistransient}} that when $I(f)<\infty$, any possible weak limit of our conditioned process $M$ is transient. Theorem \mbox{\ref{transientthm}} finds the distribution of the process $X$ in this case.

 \b{theorem} \label{transientthm}   With assumptions on $f$ and $\ov$  in case (\hbox{\hyperref[case1]{i}}) or case (\hbox{\hyperref[case2]{ii}}), if $I(f)  <\infty$,    
  \n then the measure $\mathbb{Q}(\cdot):=\lim_{t\to\infty}\bb{P}(\cdot|\cal{O}_t)$ exists for the space $\cal{D}[0,\infty)$, to which $X$ belongs,  of  c\`adl\`ag  paths on $[0,\infty)$,  in the sense that   for all    $h>0$, $y>g(h)$,  for all $\cal{B}_h \subseteq \cal{O}_h, \cal{B}_h \in \cal{F}_h$, where $(\cal{F}_u)_{u\geq0}$ is the natural filtration of $X$,
  \[
  \mathbb{Q}(X_h \in dy ; \cal{B}_h )= \f{  \Phi_y^h(\infty) }{  \Phi(\infty)  }     \pp\l(     X_h \in dy ; \cal{B}_h       \r) ,
  \]
   where  $\Phi(\infty)<\infty$, $\Phi_y^h(\infty) <\infty$.    Define, independently of $M$ or $X$, the random variable $\frak{C}$  by
 \b{equation} \label{frakC}
 \pp\l( \mathfrak{C}  \in ds\r) :=       \f{   \pp(\oo_s)   }{ \Phi(\infty)   } ds    ,  \quad s\geq0,
 \e{equation}
which exists since $\Phi(\infty)<\infty$. Then   for all $h\geq 0$, $\mathbb{Q}(X_h <\infty)= \pp( \frak{C}>h)$.
  \e{theorem}

    \b{remark} \label{transientremark}    Since $\mathbb{Q}(X_h <\infty)= \pp( \frak{C}>h)$ for all $h>0$,    $X$ under $\bb{Q}$ is finite until a random time, which we denote by $T_\infty$, and which has  the same distribution under $\bb{Q}$ as $\frak{C}$ under $\bb{P}$. In particular, $\bb{Q}(T_\infty \in ds) = \pp(\frak{C} \in ds)$ for all $s\geq0$.  In Theorem \mbox{\ref{Mistransient}}, we will show that under any possible weak limit measure, the process $M$ never returns to 0 after time $X_{T_{\infty}-}= \lim_{s\uparrow T_{\infty}} X_s$.
  \e{remark}
   
 We shall now     determine    the behaviour of the conditioned Markov process $M$ until its first excursion longer than  $g(t)$, and the time at which this excursion occurs.  Describing the behaviour until the first excursion  longer than  $g(t)$, as $t\to\infty$, in fact gives full knowledge of how $M$ behaves under $\mathbb{Q}$ until the time of its final, infinite excursion. We   verify in Proposition  \hbox{\ref{transientlemma}} that $M$ is transient under $\bb{Q}$ by showing that the process $M$ never returns to 0 after the start of this final excursion.   Remark \mbox{\ref{panti}} considers the behaviour after this time.
  
   \b{proposition}   \label{transientcorollary}   In cases (\hbox{\hyperref[case1]{i}}) and (\hbox{\hyperref[case2]{ii}}),       $  ((X_u)_{\D^{g(t)}> u\geq0} , \D^{g(t)}), t\geq0 $  under $\bb{P}(\cdot|\cal{O}_t)$ converges  as $t\to\infty$, 
   in the sense that there exists a unique limit measure $\bb{Q}'(\cdot)$, on the space $\cal{D}[0,\infty)\times(0,\infty)$, such that  for all $y > g(x)$, $ t > b > a > x >0$, with $a,b,x,y$   fixed, and for all events $\cal{B}_X \in \cal{F}_{x}$, where $(\cal{F}_u)_{u\geq0}$ is the natural filtration of $X$, such that $\cal{B}_X\subseteq \cal{O}_x$,  for $\frak{C}$ as in $\l(\ref{frakC}\r)$, 
  \b{equation}  
 \begin{split} 
    \lim_{t\to\infty}    \pp(       X_x \in dy   ;  \cal{B}_X ;     \D^{g(t)} \in (a,b) | \cal{O}_{t}   ) 
   \!     &=  \!  
    \int_a^b  \!   \bb{P}(X_x \! \in  \!   dy ; \cal{B}_X | \cal{O}_s)\pp(\frak{C} \in ds)    \\   
      \label{statementofproblemlemma}
    &=: \bb{Q}'\l(  X_x \in dy   ;  \cal{B}_X ;    T_{\infty}' \in (a,b)   \r),
    \end{split}
  \e{equation}
 for the explosion time $T_\infty'$ defined under $\bb{Q}'$ as $T_\infty$ is defined under $\bb{Q}$ in Remark \mbox{\ref{transientremark}}. The projection of $\bb{Q}'$  onto $(0,\infty)$ agrees with $\bb{Q}$ in the sense that  $\bb{Q}'( T_\infty' \!\in\! ds )\!=\!\bb{Q}( T_\infty \!\in\! ds )\!=\!\pp(\frak{C}\!\in\! ds)$.  
  \e{proposition}
  

\n  We shall now     determine    the behaviour of the conditioned process $M$ until a time corresponding to the point at which $X$ becomes infinite. 
    Theorem \mbox{\ref{Mistransient}} and Remark \mbox{\ref{panti}}    consider the behaviour after this time. Proposition  \mbox{\ref{transientlemma}}   
 requires some understanding of excursion theory of Markov processes.  For  background on excursion theory, we direct the reader to \cite[Chapter IV]{b98}.  
  \b{proposition} \label{transientlemma}   In cases (\hbox{\hyperref[case1]{i}}) and (\hbox{\hyperref[case2]{ii}}),  if $I(f)<\infty$ then there exists a measure $\mathbb{Q}''(\cdot)$ on  the  product space of the space containing the excursion process with $\cal{D}[0,\infty)\times(0,\infty)$,  such that   for all fixed $b>a>h>0$, and for $\cal{B}\subseteq \cal{O}_h$, $\cal{B} \in \cal{F}_h$,  where $\cal{F}$ denotes the natural filtration of $X$, with $F_1$    a bounded continuous functional on the excursion process $(\eps_s)_{s\geq0}$ of $M$,  
 defining the operator $\pi_h( (Z_u)_{u\geq0} ) :=   (Z_u)_{h\geq u \geq 0}$, and  letting $F_1$ satisfy $F_1((\eps_s)_{s\geq0}) = F_1 ( \pi_h((\eps_s)_{s\geq0} )))$, so $F_1$ depends only on the excursion process of $M$ up to time $h$,  we have
     \b{align} \label{E}   
 & \color{black}   \lim_{t\to\infty} \      \mathbb{E}  \l[   \ F_1(\pi_h((\eps_s)_{s\geq0} ))   \      \bbm{1}_{ \{ \pi_h(X) \in \cal{B}  \}  }       \bbm{1}_{ \{ 
\D^{g(t)} \in (a,b)  \}  }          \    \big|  \  \cal{O}_t  \ \r]   
     \\
 &\color{black} =         \int_{\nu \in \cal{B}}     \int_{u\in(a,b)}      \mathbb{E}  \l[       F_1(\pi_h((\eps_s)_{s\geq0} ) )       \big|      \pi_h(X) = \nu               \r]          \mathbb{Q}'\l(     \pi_h(X) \in d\nu        ;    T_\infty' \in du  \r)                  \nonumber
\\
  &=:       \mathbb{E}_{\mathbb{Q}''}  \l[       F_1(\pi_h((\eps_s)_{s\geq0} ) )        \    \bbm{1}_{\{ \pi_h(X) \in \cal{B} \} }  \     \bbm{1}_{ \{     T_\infty''  \in (a,b) \} } \r],                     \nonumber  
  \end{align}  
 where $T_\infty'' $ is the    explosion time under the measure $\mathbb{Q}''(\cdot)$, and the projection of $\bb{Q}''( \cdot )$ onto $\cal{D}[0,\infty)\times(0,\infty)$ agrees with $\bb{Q}'( \cdot )$. 
In particular, with $\Delta:=\Delta_1^{g(t)}$,    $((M_t)_{X_{\Delta-} > t\geq0}, (X_s)_{\Delta > s \geq0}, \Delta)$ under  $\pp(\cdot|\cal{O}_t)$ converges weakly as $t\to\infty$ to $((M_t)_{X_{T_{\infty}''-} > t\geq0}, (X_s)_{\infty > s \geq0}, T_\infty'')$ under $\mathbb{Q}''( \cdot )$. The behaviour of $M$ under $\bb{Q}''( \cdot )$ before time $X_{T_{\infty}''-}$ has the same distribution as the following construction, expressed in terms of the original measure $\bb{P}$ as follows: sampling the random time  $\frak{C}=s$ under $\bb{P}( \cdot )$, we run $X$ conditioned on $\cal{O}_s$ until time $s$,  take $X_u=\infty$ for all $u\geq s$, then construct $ M$ via its excursions using $(X_u)_{\infty>u\geq0}$ to determine the timing and length of each excursion,  where we sample each excursion of $ M$  until time $X_{s-}$ using the excursion measure conditional on the given excursion length. 
  \e{proposition}
  
 \b{theorem} 
  \label{Mistransient} 
  In cases (\hbox{\hyperref[case1]{i}}) and (\hbox{\hyperref[case2]{ii}}),  if $I(f)<\infty$ then $M$ is transient under any possible weak limit of the measure $\pp(\cdot|\cal{O}_t)$ as $t\to\infty$. 
   \e{theorem}

\b{remark}  \label{panti}   
While the last excursion of the Markov process $M$ is not   dealt with explicitly here, the behaviour of $M$ from time $X_{\frak{C}-}$ onwards  should be the same as that of $M$ conditioned to avoid zero. Proving this requires existence of the limit as $g(t)\to\infty$ of the excursion measure conditioned on the length (lifetime) of an excursion being longer than $g(t)$, which is beyond the scope of this work. This is verified in the simple, single case where $M$ is a 1-dimensional Brownian motion in \cite[p8]{ks13}.  When $M$ is a L\'evy process,  the behaviour of the process conditioned to avoid zero is well understood,  see   \cite[Theorem 8]{p13}.  
 There is some technical difficulty in applying results from \cite{p13} to our final excursion.  The measures $\mathbb{Q},\mathbb{Q}',\mathbb{Q}''$ are constructed by conditioning until a deterministic time $t\to\infty$, but in \cite{p13}, the measure is constructed by conditioning until an independent exponential random time with   parameter $q\to0$. Equivalence of such deterministic and random limits is a separate~matter, beyond the scope of this work. 
 \e{remark}   


  %


\subsection{Results in the  \textit{I}(\textit{f})    $=\infty$ Case}  
    \label{recurrentsection} 
     We now restrict our attention to case (\hbox{\hyperref[case1]{i}}).  We will see that when $I(f)=\infty$, our conditioned Markov process is recurrent. Theorem \mbox{\ref{recurrentthm}} finds the distribution of the conditioned inverse local time subordinator in this case.

 \b{theorem} \label{recurrentthm}  In case (\hbox{\hyperref[case1]{i}}), if $I(f)   =\infty$, then  the   law  $\mathbb{Q}(\cdot)= \lim_{t\to\infty} \pp(   \cdot  | \cal{O}_t  )$ exists for the process $X$ in the sense that for all $h>0$ and $y\geq g(h)$,    for all    $h>0$, $y>g(h)$,  for all $\cal{B}_h \subseteq \cal{O}_h, \cal{B}_h \in \cal{F}_h$, where $(\cal{F}_u)_{u\geq0}$ is the natural filtration of $X$,
\b{align} \nonumber  
 \mathbb{Q}\l( X_h \in dy ; \cal{B}_h    \r)     &=    \lim_{t\to\infty}    \f{  \mathbb{P}(   \cal{O}_{t-h}^{g_{y}^{h}}        )      }{  \mathbb{P}(\cal{O}_t)     }       \pp\l( X_h \in dy ; \cal{B}_h    \r)
  \\ \label{qqq}
      &=:      q_h(y)   \pp\l( X_h \in dy ; \cal{B}_h    \r),  
  \end{align}
   where  $q_h(y)$ is finite, non-decreasing in $y$,  and satisfies 
  \[
q_h(y)       =       \f{    \Phi_y^h(t_0(y))}{    \Phi(1)}
      \exp\l(    \int_{  t_0(y)}^{\infty}   \l(   \ov(g_{y}^{h}(s)) -  \ov(g(s))  + \rho_y^h(s)   \r)  ds     -  \int_{1}^{\infty}            \rho(s)    ds      
 -        \int_1^{ t_0(y)   }       \ov(g(s)) ds     \r)    ,
 \] 
 \b{equation*}    
    t_0(y) :=  f(Ay)\vee  f(1 + 2/A),   \qquad A>3\vee (B-1),      \phantom{\Big(}   
    \e{equation*}
     for $B>0$ as in Assumption~\mbox{\ref{case1}}, and where each integral is finite. In the case that $t_0(y)<1$, the final integral should be interpreted as $- \int_1^{ t_0(y)   }       \ov(g(s)) ds =  \int_{ t_0(y)   }^1       \ov(g(s)) ds $. 
   \e{theorem}

   \n    We now verify that when $I(f)=\infty$, $M$ is recurrent  under the new measure $\mathbb{Q}''(\cdot)$, as $X$ never hits infinity at a finite time, $\mathbb{Q}$-almost surely, and then $M$ under $\bb{Q}''(\cdot)$ is constructed from its excursion process and $X$. 

 \b{proposition} \label{lemma5}  
  In   case (\hbox{\hyperref[case1]{i}}), if $I(f) =\infty$, then  for each $h>0$,
  \(
  \mathbb{Q}\l( X_h \in (g(h),\infty)    \r)        =1.
  \)
  \e{proposition}

 \b{proposition} \label{recurrentmess}  
In   case (\hbox{\hyperref[case1]{i}}), if $I(f) =\infty$,   then there exists a measure $\mathbb{Q}''(\cdot)$ on the product space of the space containing the excursion process with the space $\cal{D}[0,\infty)$ of c\`adl\`ag paths on $[0,\infty)$, such that for all fixed $h>0$, and for $\cal{B}\subseteq \cal{O}_h$, $\cal{B} \in \cal{F}_h$,  where $\cal{F}$ denotes the natural filtration of $X$, let $F_1$   be a bounded continuous functional on the excursion process $(\eps_s)_{s\geq0}$ of $M$,  
 defining the operator $\pi_h( (Z_u)_{u\geq0} ) :=   (Z_u)_{h\geq u \geq 0}$, with $F_1$ such that $F_1((\eps_s)_{s\geq0}) = F_1 ( \pi_h((\eps_s)_{s\geq0} )))$,  
     \b{align} \label{Erec}   
 &     \lim_{t\to\infty} \      \mathbb{E}  \l[   \ F_1(\pi_h((\eps_s)_{s\geq0} ))   \      \bbm{1}_{ \{ \pi_h(X) \in \cal{B}  \}  }               \    \big|  \  \cal{O}_t  \ \r]   
     \\
 &    =         \int_{\nu \in \cal{B}}         \mathbb{E}  \l[       F_1(\pi_h((\eps_s)_{s\geq0} ) )       \big|      \pi_h(X) = \nu               \r]          \mathbb{Q}\l(     \pi_h(X) \in d\nu         \r)                  \nonumber
\\
  &  =:       \mathbb{E}_{\mathbb{Q}''}  \l[     F_1(\pi_h((\eps_s)_{s\geq0} ) )         \    \bbm{1}_{\{ \pi_h(X) \in \cal{B} \} }   \r],                     \nonumber  
  \end{align}  
 where  the projection of $\bb{Q}''( \cdot )$ onto $\cal{D}[0,\infty)$ agrees with $\bb{Q}( \cdot )=\lim_{t\to\infty}\bb{P}(\cdot|\cal{O}_t)$. 
In particular, as $t\to\infty$, $((M_t)_{\infty> t\geq0}, (X_s)_{\infty > s \geq0})$ under  $\pp(\cdot|\cal{O}_t)$ converges weakly  to $((M_t)_{\infty > t\geq0}, (X_s)_{\infty > s \geq0})$ under $\mathbb{Q}''(\cdot)$. We  construct $ M$ via its excursions using $(X_u)_{\infty>u\geq0}$ to determine the timing and length of each excursion,  where we sample the excursions of $ M$  using the excursion measure conditional on each excursion length. Moreover,  
 $M$ visits 0 at arbitrarily large times, so
  $M$ is recurrent  under  $\mathbb{Q}''(\cdot)$.
   \e{proposition}   
 

 \n Now we shall determine the entropic repulsion envelope through Theorem \mbox{\ref{entrrep}}.  

     \b{definition}   \label{entdefn} A non-decreasing function $w$, with $\lim_{h\to\infty}w(h)  =  \infty$, is in the entropic repulsion envelope $R_g$ (for the function $g=f^{-1}$) if 
  \b{equation}   
  \label{entdef}
  \lim_{h\to\infty}  \mathbb{Q}''\l(  X_h \geq   w(h) g(h)    \r)   =1.
  \e{equation}
  \e{definition}
\b{theorem} \label{entrrep} In case (\hbox{\hyperref[case1a]{ia}}), a necessary and sufficient condition for non-decreasing  $w$,  for which $\lim_{h\to\infty}w(h)=\infty$, to be in $R_g$ (for the function  $g=f^{-1}$)  is
  \[
   w\in R_g \iff  \lim_{h\to\infty}   \int_h^{  f( w(h) g(h)   )  }     \ov(g(s)) ds =0.
  \]
  \e{theorem}
 \n As a result, one can    verify that the entropic repulsion envelope is always non-empty in case (\hbox{\hyperref[case1a]{ia}}),  by finding suitable  $w$.  We illustrate the generality of Theorem \mbox{\ref{entrrep}} via  the following corollary,   expanding upon \cite[Theorem 4]{ks13}.

  \b{corollary} \label{entrrepstable} In case (\hbox{\hyperref[case1]{i}}), with $f,f'$    \mbox{\hyperref[orv]{$\cal{O}$-regularly varying}}   at $\infty$, for a stable subordinator of index $\alpha\in(0,1)$,
  a necessary and sufficient condition for non-decreasing  $w$,   with $\lim_{h\to\infty}w(h)=\infty$, to be in  $R_g$, for   $g=f^{-1}$  is
  \[
   w\in R_g \iff  \lim_{h\to\infty}   \int_h^{  f( w(h) g(h) )  }   g(s)^{-\alpha} ds =0.
  \]
  \e{corollary}



\section{Proofs of Main Results} \label{mainresultsproofs} 
This section contains the proofs of the results stated in Section \mbox{\ref{mainresults}}. First we state Lemmas 
 \mbox{\ref{lemma2}} and \mbox{\ref{lowerboundforrho}}, which are proven in Sections 
  \mbox{\ref{lemma2proof}} and \mbox{\ref{lemmasproofs}}, respectively.

   \b{definition} In this paper we   use the following asymptotic notation:
 
 $f(x) \sim g(x) $ as $x\to \infty$ if $\lim_{x\to \infty} f(x)/g(x) =1$.

   $f(x) \lesssim g(x) $ if  there exists $C\in(0,\infty)$ such that for all large enough $x$, $f(x) \leq C g(x)$.
 
 \n   Moreover, we write   $f(x) \gtrsim g(x) $   if $g(x)\lesssim f(x)$,  and
 
 $f(x) \asymp g(x) $ if both $f(x) \gtrsim g(x) $ and $f(x) \lesssim g(x) $.

   \e{definition}
  
  \b{lemma} \label{lemma2}  In cases (\hbox{\hyperref[case1]{i}}) and (\hbox{\hyperref[case2]{ii}}),  there exists a function $u(t)$ with $\lim_{t\to\infty} u(t)    =   0$, and there exists $\eps>0$  such that for all    $A>3$, uniformly in $h>0$, $y>g(h)$, and  $ t> t_0(y)$ as defined in $\l(\ref{t0y}\r)$,  
  \b{equation} \label{1}
  \rho_y^h(t)     \lesssim  \f{1}{   t \log(t)^{1+\eps} }    \l( 1 + \f{1}{f(y)-h} \r) ,
  \e{equation}
  \b{equation}\label{2}
    \rho_y^h(t)   \leq         u(t)\ov(g(t))   \l( 1 + \f{1}{f(y)-h} \r) .
  \e{equation}
  The inequalities $\l(\ref{1}\r)$ and $\l(\ref{2}\r)$ also hold when $y=h=0$, for $t>t_0(0)>0$,  with $\rho(t)$ in place of $\rho_y^h(t)$.
   \e{lemma} 
  \b{lemma} \label{lowerboundforrho} In cases (\hbox{\hyperref[case1]{i}}) and (\hbox{\hyperref[case2]{ii}}), for the function $\rho$ as  defined in $\l(\ref{rho}\r)$, $\rho(t)=o(\ov(g(t))$ as $t\to\infty$, and moreover, $\int_1^\infty \rho(s)ds >-\infty$.   
 \e{lemma}
 Corollary \mbox{\ref{lemma1}} follows immediately from Lemmas \mbox{\ref{lemma2}} and  \mbox{\ref{lowerboundforrho}}.
 \b{corollary} \label{lemma1}
 In cases (\hbox{\hyperref[case1]{i}}) and (\hbox{\hyperref[case2]{ii}}),  as $t\to\infty$,
   \b{equation*} 
   \pp\l( \oo_{t}   \r) =  \l(  \ov(g(t))+\rho(t)\r) \Phi(t)     = ( 1+ o(1)) \  \ov(g(t)) \Phi(t).
   \e{equation*}
 \e{corollary}
  \b{remark} \label{yh0}  
  Taking $y=h=0$ in $\l(\ref{1}\r)$, it follows  that $\int_{1}^\infty   \rho(s)ds<\infty$. Then as  $\int_1^\infty \rho(s)ds >-\infty$ by Lemma \mbox{\ref{lowerboundforrho}}, it follows immediately from  $\l(\ref{soln1}\r)$ that  as $t\to\infty$, 
  \b{equation}
  \label{phiasymp}
  \Phi(t)   =   \Phi(1) \exp     \l(  \int_{1}^t (\ov(g(s))  +    \rho(s))ds       \r)    \asymp     \exp   \l(   \int_{1}^t \ov(g(s))ds     \r)  .
  \e{equation}
  \e{remark} 



 \subsection{Proofs in the \textit{I(f)} $<\infty$ Case }

  \subsubsection{Proof of Theorem  \hbox{\ref{transientthm}}}
  \b{proof}[Proof of Theorem  \hbox{\ref{transientthm}}] First, let us verify that   $\Phi(\infty)<\infty$. Recalling that $g  :=  f^{-1}$, we have  
\b{equation}   \label{If}
 I(f)    
 =           \int_1^\infty       f(x) \Pi(dx)  
 =       \hspace{-0.1cm}    \int_1^\infty      \int_0^{f(x)}     dy \Pi(dx) 
 =       \int_{0}^{\infty}       \int_{1\vee g(y)}^\infty       \Pi(dx) dy 
 =       \int_{0}^\infty  \ov(1\vee g(y))dy    .
\e{equation}
Now, recall from $\l(\ref{soln1}\r)$ that 
\[
\Phi(\infty) =    \Phi(1)    \exp \l(    \int_{1}^\infty    \ov(g(s))   ds     +       \int_{1}^\infty    \rho(s)   ds        \r).
\]
By Corollary  \hbox{\ref{lemma1}},  as $s\to\infty$,   $\rho(s) =o(   \ov(g(s)) )$.   Then by $\l(\ref{If}\r)$, since $I(f)<\infty$,
\[
   \int_{1}^\infty    \ov(g(s))   ds     +       \int_{1}^\infty    \rho(s)   ds     
 \overset{\ref{lemma1}}\lesssim       \int_{1}^\infty    \ov(g(s))   ds    \overset{\l(\ref{If} \r)}<  \infty ,
 \]
  so $\Phi(\infty)<\infty$. Now,  $\int_1^\infty \rho(s)ds >-\infty$ by Lemma \mbox{\ref{lowerboundforrho}}, and hence 
\begin{equation}  \label{IfPhi}
   I(f) < \infty    \iff     \Phi(\infty)<\infty.
\end{equation}
To show $\Phi_y^h(\infty)<\infty$,  with  $t_0(y)$ as defined in $\l(\ref{t0y}\r)$,  recall that by $\l(\ref{soln2}\r)$,
  \[
  \Phi_y^h(\infty) =    \Phi_y^h(t_0(y))    \exp \l(    \int_{t_0(y)}^\infty    \ov(g_{y}^{h}(s))   ds     +       \int_{t_0(y)}^\infty   \rho_y^h(s)   ds        \r).
  \]
Now, observe that for each fixed $y,h>0$, $g(s)\sim g_{y}^{h}(s)$ as $s \to\infty$, by  $\l(\ref{gyh}\r)$ and the properties of $g$ introduced in Assumption \mbox{\ref{case1}}, so  $\ov(g(s)) \sim \ov(g_{y}^{h}(s))$ as $s\to\infty$, since $\ov$ is  \mbox{\hyperref[crv]{CRV}} at $\infty$. Now,  applying  $\l( \ref{2} \r)$ and $\l(\ref{If}\r)$, since $y$ and $h$ are fixed and $y>g(h)$ implies $f(y)-h>0$, noting $s>t_0(y)$ ensures $g_y^h(s)>0$, we get
  \[
      \int_{t_0(y)}^\infty   \ov(g_{y}^{h}(s))   ds      +       \int_{t_0(y)}^\infty    \rho_y^h(s)   ds     \overset{   \l( \ref{2} \r)   }{\lesssim}        \l(   1 + \f{1}{f(y)-h} \r) \int_{t_0(y)}^\infty   \ov(g(s))   ds     \overset{\l(\ref{If}\r)}<    \infty,
  \]
\n    so $\Phi_y^h(\infty)<\infty$. By $\l(\ref{doobh}\r)$ and Corollary \mbox{\ref{lemma1}},  since $\ov(g(t)) \sim \ov(g_{y}^{h}(t-h))$ as $t\to\infty$,
 \b{align}
 \nonumber
\mathbb{Q}(X_h\in dy ; \cal{B}_h)    &\overset{  \l(\ref{doobh}\r)}=  \pp\l(    X_h \in dy      ;  \cal{B}_h    \r)     \lim_{t\to\infty}     \f{ \pp\big(   \oo_{t-h}^{g_{y}^{h}}          \big) }{   \pp\l(  \oo_{t}  \r)    }   
\\
  &\overset{  \hspace{1.9pt}   \ref{lemma1}     \hspace{1.9pt}    }=      \pp\l(    X_h \in dy      ; \cal{B}_h    \r)     \lim_{t\to\infty}     \f{\ov(g_{y}^{h}(t-h)) \Phi_y^h(t-h) }{ \ov(g(t)) \Phi(t)   } 
  \label{PQ}
 =      \pp\l(    X_h \in dy      ; \cal{B}_h    \r)      \f{  \Phi_y^h(\infty) }{  \Phi(\infty)   } .  
  \end{align}
 
  \n Now we show  
  \(
     \mathbb{Q}(X_h     <          \infty)      \hspace{-2pt}       =   \hspace{-2pt}      \pp\l( \frak{C}                        >                       h   \r)   
  \). 
   Applying $\l(\ref{PQ}\r)$  with $\cal{B}_h=\cal{O}_h$,
 \b{align*}
\mathbb{Q}(X_h <\infty) &=     \int_{g(h)}^\infty  \mathbb{Q}(X_h \in dy)  =   \int_{g(h)}^\infty     \f{\Phi_y^h(\infty)    }{  \Phi(\infty)   }      \pp\l(    X_h \in dy  ; \cal{O}_h        \r)
 \\
  &=     \f{1}{\Phi(\infty)}  \int_{g(h)}^\infty    \int_0^\infty \pp\l(   \oo_{v}^{g_{y}^{h}}    \r) dv       \pp\l(    X_h \in dy         ; \cal{O}_h \r)    \\ 
  &= \f{1}{\Phi(\infty)}      \int_0^\infty    \int_{g(h)}^\infty    \pp\l( \oo_v^{g_{y}^{h}}    \r) \pp\l(   X_h\in dy   ; \cal{O}_h \r)  dv    .
 \intertext{Now,  $  \pp\big( \oo_v^{g_{y}^{h}}    \big) \pp\l(   X_h\in dy   ; \cal{O}_h \r) = \pp\l( \oo_{v+h}  ;  X_h\in dy   \r) $ by $\l(\ref{ogyh}\r)$. Then      by  the definition $\l(\ref{o}\r)$ of $\cal{O}_{v+h}$,  }   
 \mathbb{Q}(X_h <\infty)     
 &\overset{\phantom{(\ref{o})}}=      \f{1}{\Phi(\infty)}      \int_0^\infty    \int_{g(h)}^\infty  \hspace{-5pt}   \pp\l( \oo_{v+h}  ;  X_h\in dy   \r)  dv  
&&\hspace{-65pt}\overset{\phantom{\ref{o}}}=     \f{1}{\Phi(\infty)}      \int_0^\infty   \hspace{-5pt}     \pp\l( \oo_{v+h}  ;  X_h > g(h)   \r)  dv\\ 
      &\overset{(\ref{o})}=   \hspace{30pt}  \f{1}{\Phi(\infty)}      \int_0^\infty        \pp\l( \oo_{v+h}      \r)  dv \hspace{-30pt}
       &&\hspace{-65pt}\overset{\phantom{\ref{o}}}=     \f{1}{\Phi(\infty)}      \int_h^\infty    \hspace{-5pt}    \pp\l( \oo_{u}     \r)  du =: \pp\l( \frak{C}   >h   \r)  .
  \end{align*}
\e{proof}

  \subsubsection{Proof of Proposition  \hbox{\ref{transientcorollary}}}

  \b{proof}[Proof of Proposition  \hbox{\ref{transientcorollary}}]

  \n For $y > g(x)$, $ t > b > a > x >0$, with $a,b,x,y$   fixed, and an event $\cal{B}_X \in \cal{F}_{x}$, where $(\cal{F}_u)_{u\geq0}$ is the natural filtration of $X$, such that $\cal{B}_X\subseteq \cal{O}_x$, consider 
  \b{equation}    \label{Jt}
    \lim_{t\to\infty}    \pp\l(       X_x \in dy   ;  \cal{B}_X ;     \D^{g(t)} \in (a,b) | \cal{O}_{t}   \r).
  \e{equation}
   If $\D^{g(t)} \in ds$, then $X_s \! >\!  g(t)$, so   $\cal{O}_t$ is fully  attained by time $s$, and    $\cal{O}_t$  can be replaced by $\cal{O}_s$, so
  \b{align}
  \nonumber
\l(\ref{Jt}\r) &=   \lim_{t\to\infty}    \f{1}{\pp(\oo_t )}   \int_a^b         \pp\l(     X_x   \in dy  ;  \cal{B}_{X } ;     \D^{g(t)} \in ds ; \cal{O}_{t}       \r) 
\\
\nonumber
  &=   \lim_{t\to\infty}    \f{1}{\pp(\oo_t )}   \int_a^b         \pp\l(    X_x   \in dy    ;  \cal{B}_X  ;     \D^{g(t)} \in ds ; \cal{O}_{s}       \r) .
  \intertext{Recall the definition $\l(\ref{truncoo}\r)$. Given $\Delta_1^{g(t)} >a>x$, we can replace $ \{ X_x \in dy \},   \cal{B}_X$, $\cal{O}_s$ by     corresponding events  $ \{ X_x^{(0,g(t))} \in dy \} , \cal{B}_{X^{(0,g(t))}}$, $\cal{O}_{s,X^{(0,g(t))}}$  for the   process $X^{(0,g(t))}$ with L\'evy measure  restricted to $(0,g(t))$, i.e.\    all jumps  larger than $g(t)$ are removed. These events are each independent~of~$\Delta_1^{g(t)}$, and since $\Delta_1^{g(t)}$ is exponentially distributed with parameter $\ov(g(t))$,}
  \nonumber
 \l(\ref{Jt}\r)  &=    \lim_{t\to\infty}     \f{1}{\pp(\oo_t )}     \int_a^b         \pp\l(    X_x^{(0,g(t))} \in dy   ;  \cal{B}_{X^{(0,g(t))}} ;     \D^{g(t)} \in ds ; \cal{O}_{s,X^{(0,g(t))} }       \r) 
\\
\nonumber
   &=   \lim_{t\to\infty}      \f{1}{\pp(\oo_t )}     \int_a^b         \pp\l(    X_x^{(0,g(t))} \in dy   ; \cal{B}_{X^{(0,g(t))}}  ; \cal{O}_{s,X^{(0,g(t))} }       \r)   \pp\l(   \D^{g(t)} \in ds\r)
\\
     &=    \lim_{t\to\infty}    \f{     \ov(g(t))    }{\pp(\oo_t )}     \int_a^b         \pp\l(    X_x^{(0,g(t))} \in dy   ;  \cal{B}_{X^{(0,g(t))}} ; \cal{O}_{s,X^{(0,g(t))} }       \r)   e^{- \ov(g(t)) s} ds.
 \intertext{Now, since $\lim_{t\to\infty} e^{- \ov(g(t)) s}   =1$, uniformly among $s\in(a,b)$, }
 \nonumber
     \l(\ref{Jt}\r)      &=   \lim_{t\to\infty}       \f{     \ov(g(t))    }{\pp(\oo_t )}    \int_a^b         \pp\l(    X_x^{(0,g(t))} \in dy   ;  \cal{B}_{X^{(0,g(t))}}  ; \cal{O}_{s,X^{(0,g(t))} }       \r)  ds.
   \intertext{Applying   Corollary  \hbox{\ref{lemma1}}, and recalling  from $\l(\ref{IfPhi}\r)$ that   $\Phi(\infty)<\infty$  when $I(f)<\infty$, }
 \l(\ref{Jt}\r)   &=   \f{ 1}{   \Phi(\infty)     }    \lim_{t\to\infty}       \int_a^b         \pp\l(    X_x^{(0,g(t))} \in dy   ;  \cal{B}_{X^{(0,g(t))}}  ; \cal{O}_{s,X^{(0,g(t))} }       \r)  ds
\\
\nonumber
  &=      \f{1}{   \Phi(\infty)     }  \lim_{t\to\infty}      \int_a^b         \pp\l(    X_x^{(0,g(t))} \in dy   ;  \cal{B}_{X^{(0,g(t))}}   \big|  \cal{O}_{s,X^{(0,g(t))} }       \r)   \pp\l( \cal{O}_{s,X^{(0,g(t))} }  \r)  ds.
     \intertext{Now, $\lim_{t\to\infty}   \pp( \cal{O}_{s,X^{(0,g(t))} }  )  /    \pp\l( \cal{O}_{s }  \r)     =1$, uniformly among $s \in (a,b)$, so}
 \nonumber
 \l(\ref{Jt}\r)   &=      \lim_{t\to\infty}          \int_a^b         \pp\l(    X_x^{(0,g(t))} \in dy   ;  \cal{B}_{X^{(0,g(t))}}   \big|  \cal{O}_{s,X^{(0,g(t))} }       \r)   \f{ \pp\l( \cal{O}_{s  }  \r)}{   \Phi(\infty)     }     ds,
  \intertext{and similarly $       \pp(    X_x^{(0,g(t))} \in dy   ;  \cal{B}_{X^{(0,g(t))}}   \big|  \cal{O}_{s,X^{(0,g(t))} }       )  \sim          \pp\l(    X_x \in dy   ;  \cal{B}_X   \big|  \cal{O}_{s }       \r) $ as $t\to\infty$, uniformly among $s\in(a,b)$. Then  by the definition of $\frak{C}$ in $\l(\ref{frakC}\r)$,         }
  \begin{split}
 \l(\ref{Jt}\r) &=           \int_a^b         \pp\l(    X_x  \in dy   ;  \cal{B}_X  \ \big| \ \cal{O}_{s }       \r)    \f{  \pp\l( \cal{O}_{s }  \r)   }{  \Phi(\infty)       } ds
 =         \int_a^b         \pp\l(    X_x  \in dy   ;  \cal{B}_X  \ \big| \ \cal{O}_{s }       \r)   \bb{P}\l(   \frak{C} \in ds\r)
\\
\label{pqpq}
 &=:       \bb{Q}'\l(    X_x  \in dy   ;  \cal{B}_X   ;  T'_\infty  \in (a,b)\r). 
 \end{split}
    \end{align}
  %
  %
It is clear that $\l(\ref{pqpq}\r)$ uniquely determines the limit measure $\bb{Q}'(\cdot)$ on $\cal{D}[0,\infty)\times(0,\infty)$. To verify that $T'_\infty$ under $\bb{Q}'(\cdot)$ has  the desired properties,  by $\l(\ref{pqpq}\r)$ with $\cal{B}_X = \cal{O}_x$, since $x<a< s$,
   \[
\bb{Q}'( T'_\infty \in (a,b))\!=\! \int_{g(x)}^\infty  \bb{Q}'\l(    X_x  \in dy   ;  \cal{O}_x   ;  T'_\infty  \in (a,b)\r)  
\!=\!
 \int_{g(x)}^\infty    \int_a^b         \pp\l(    X_x  \in dy   ;   \cal{O}_x  \ \big| \ \cal{O}_{s }       \r)   \bb{P}\l(   \frak{C} \in ds\r)
   \]
   \[
   =     \int_a^b         \pp\l(    X_x  > g( x)    ;  \cal{O}_x  \ \big| \ \cal{O}_{s }       \r)   \bb{P}\l(   \frak{C} \in ds\r) =  \int_a^b           \bb{P}\l(   \frak{C} \in ds\r)  = \pp(  \frak{C} \in(a,b)) = \bb{Q}(T_\infty\in(a,b)). 
   \]
   Similarly, by $\l(\ref{pqpq}\r)$ with $\cal{B}_X = \cal{O}_x$, taking limits as $a\to x$ and $b\to\infty$, since $x<s$, we also have
\[
 \bb{Q}'( X_x < \infty) = \int_{g(x)}^\infty  \int_x^\infty  \pp(X_x\in dy ; \cal{O}_x | \cal{O}_s) \pp(\frak{C} \in ds)
\]
\[
=     \int_x^\infty  \pp\l(    X_x  > g( x)    ;  \cal{O}_x  \ \big| \ \cal{O}_{s }       \r)   \bb{P}\l(   \frak{C} \in ds\r) =     \int_x^\infty    \bb{P}\l(   \frak{C} \in ds\r)= \pp(\frak{C}>x) = \bb{Q}'(T'_\infty >x),
\]
 so that $T'_\infty $ is indeed the explosion time for the process $X$ under $\bb{Q}'(\cdot)$.
%
 %
 %
 %
 %
 %
 %
   \e{proof}

\subsubsection{Proof of Proposition  \hbox{\ref{transientlemma}}}
  \b{proof}[Proof of Proposition  \hbox{\ref{transientlemma}}]

Recall $\D^{g(t)}$ is the time of $X$'s first jump bigger than $g(t)$,  $\pi_h(X)$ is the sample  path of $X$ up to time $h$, $F_1$ is a  functional on the excursion process, 
and $\cal{B}\subseteq \cal{O}_h$, $\cal{B} \in \cal{F}_h$, where $(\cal{F}_u)_{u\geq0}$ is $X$'s natural filtration.  For  fixed $b>a>h>0 $,   disintegrating on the values of  $\D^{g(t)}$ and $\pi_h(X)$,   
  \b{equation} \label{E}
  \mathbb{E}\l[       \bbm{1}_{ \{ \pi_h(X) \in \cal{B}  \}  }    \    \bbm{1}_{ \{ 
\D^{g(t)} \in (a,b)  \}  }       F_1(\pi_h((\eps_s)_{s\geq0}))     \    \big|  \  \cal{O}_t  \ \r]   
  \e{equation}
  \[
  =    \underset{\nu\in\cal{B}}\int    \    \underset{ u \in(a,b)}\int      \mathbb{E} \big[       F_1(\pi_h((\eps_s)_{s\geq0}))      \big|      \cal{O}_t   ;  \pi_h(X) = \nu        ;    \D^{g(t)} =u     \big]      \pp\big(   \pi_h(X) \in d\nu        ;    \D^{g(t)} \in du    | \cal{O}_t   \big).
  \]
Given a fixed path  $\pi_h(X)= \nu$, $\pi_h((\eps_s)_{s\geq0})$ depends only on $\nu$, so  $\pi_h((\eps_s)_{s\geq0})$ is conditionally independent   of $\D^{g(t)}$ and $\cal{O}_t$.  Here, $h<a<u$, so the   excursion  process $(\eps_s)_{s\geq0}$  contains only excursions of length at most $g(t)$, so we may   replace $\pi_h((\eps_s)_{s\geq0}) $ by $ \pi_h( (\eps_s^{g(t)})_{s\geq0} ) $,  where    $(\eps_s^{g(t)})_{s\geq0}   $    is the excursion process sampled using the conditional excursion measure on the space    of excursions of length  at most $ g(t)$, so
  \b{align*}
  \l(\ref{E}\r)
  &=     \int_{\cal{B}}     \int_a^b       \mathbb{E} \l[      F_1(\pi_h((\eps^{g(t)}_s)_{s\geq0} ))       \big|    \pi_h(X) = \nu          \r]    \pp\l(   \pi_h(X) \in d\nu        ;    \D^{g(t)} \in du    | \cal{O}_t   \r),
  \end{align*}
  Now,  $\lim_{t\to\infty} \mathbb{E}  [  \hspace{0.01cm}     F_1(\pi_h( (\eps^{g(t)}_s)_{s\geq0} ))   \      \big|   \   \pi_h(X) = \nu  \     ] =  \mathbb{E}  [       F_1(\pi_h((\eps_s)_{s\geq0} ))       \big|      \pi_h(X) = \nu               ]$, and by Proposition \mbox{\ref{transientcorollary}}, $  \lim_{t\to\infty}      \pp(      \pi_h(X) \in d\nu        ;    \D^{g(t)} \in du    |    \cal{O}_t   ) 
      = \mathbb{Q}'\l(     \pi_h(X) \in d\nu        ;    T_\infty'  \in du  \r)    $, so 
  \b{equation}
  \b{split}
  \label{pointingatyou}
  \lim_{t\to\infty}    \l(\ref{E}\r) &=         \int_{\cal{B}}     \int_a^b       \mathbb{E}  \l[       F_1(\pi_h((\eps_s)_{s\geq0} ))       \big|      \pi_h(X) = \nu               \r]     \mathbb{Q}'\l(     \pi_h(X) \in d\nu        ;    T_\infty'  \in du  \r) 
\\
  &=:       \mathbb{E}_{\bb{Q}''}  \l[       F_1(\pi_h((\eps_s)_{s\geq0} ) )         \    \bbm{1}_{\{ \pi_h(X) \in \cal{B} \} }  \     \bbm{1}_{ \{     T_\infty''  \in (a,b) \} } \r],
  \end{split}
  \e{equation}
    where we are able to exchange the order of limits and integration since  $F_1$ is bounded. 
    Taking $F_1\equiv 1$, it follows immediately that $\bb{Q}''(\cdot)$ and $\bb{Q}'(\cdot)$ agree on $\cal{D}[0,\infty)\times (0,\infty)$.
    The weak convergence  of $((M_t)_{X_{\Delta-} > t\geq0}, (X_s)_{\Delta > s \geq0}, \Delta)$ under  $\pp(\cdot|\cal{O}_t)$  to $((M_t)_{X_{T_{\infty}''-} > t\geq0}, (X_s)_{\infty > s \geq0}, T_\infty'')$ under $\mathbb{Q}''( \cdot )$ as $t\to\infty$  then follows immediately from the fact (see e.g.\  \cite[Ex.\ IV.6.3]{b98} or \cite[p4113]{ks13}) that for all $x>0$, $(M_t)_{X_{x-}>t\geq0}$  is uniquely determined by $(\eps_s)_{x>s\geq0}$ and $(X_s)_{x>s\geq0}$, and both of $(\eps_s)_{x>s\geq0}$ and $(X_s)_{x>s\geq0}$ have weak limits as determined in $\l(\ref{pointingatyou}\r)$. That is, we construct $ M$ pathwise via its excursions using $(X_u)_{\infty>u\geq0}$ to determine the timing and length of each excursion,  where we sample the excursions of $ M$  until time $X_{s-}$ using the excursion measure conditional on each excursion length.  Similarly, the explicit description of the behaviour of $M$ until time $X_{T''_\infty-}$ under $\bb{Q}''(\cdot)$ follows immediately from the definition of $\bb{Q}'(\cdot)$ in $\l(\ref{statementofproblemlemma}\r)$, using the fact that $\bb{Q}''(\cdot)$ and $\bb{Q}'(\cdot)$ agree on $\cal{D}[0,\infty)\times (0,\infty)$.
  \e{proof}

\subsubsection{Proof of Theorem  \hbox{\ref{Mistransient}}}
  \b{proof}[Proof of Theorem  \hbox{\ref{Mistransient}}]
As $X$ determines the lengths and timings of excursions of $M$ (see \cite[Ex.\ IV.6.3]{b98}), it follows that for all $K>0$ and $ b>a>0$, for all $t$ large enough that $g(t)>K$,
 \b{equation}
 \label{Kbitt} \{ 
\D^{g(t)} \in (a,b)  \} =  \{ 
\D^{g(t)} \in (a,b)  \} \cap \{ M_v \neq 0, \textup{ for all }  v\in (X_{ \Delta_1^{g(t)}-},X_{ \Delta_1^{g(t)}-} +K) \} .
\e{equation}
Let us assume that a weak limit measure $\hat{\bb{Q}}(\cdot)=\lim_{t\to\infty}\pp(\cdot|\cal{O}_t)$ exists on the space containing $(M_t)_{t\geq0}$.
Such a measure must agree with $\bb{Q}'(\cdot)$ on $\cal{D}[0,\infty)\times(0,\infty)$, 
as we proved in Proposition  \mbox{\ref{transientcorollary}} that any such limit measure is uniquely determined on $\cal{D}[0,\infty)\times (0,\infty)$ by $\l(\ref{statementofproblemlemma}\r)$. It follows that for all $K>0$ and $ b>a>0$,
\b{equation}
\label{limMM}
\lim_{t\to\infty}\pp(  M_v \neq 0, \textup{ for all }  v\in (X_{ \Delta_1^{g(t)}-},X_{ \Delta_1^{g(t)}-} +K) ; \D^{g(t)} \in (a,b)  | \cal{O}_t )
\e{equation}
\[
=
\hat{\bb{Q}}(  M_v \neq 0, \textup{ for all }  v\in (X_{ \hat{T}_\infty-},X_{ \hat{T}_\infty-} +K) ; \hat{T}_\infty \in (a,b)    ),
\]
where $\hat{T}_\infty$ is the explosion time for $X$ under $\hat{\bb{Q}}(\cdot)$. But also by $\l(\ref{Kbitt}\r)$ and uniqueness of the limit measure on $\cal{D}[0,\infty)\times(0,\infty)$, we have for all $K>0$ and $ b>a>0$,
\[
\hat{\bb{Q}}(  M_v \neq 0, \textup{ for all }  v\in (X_{ \hat{T}_\infty-},X_{ \hat{T}_\infty-} +K) ; \hat{T}_\infty \in (a,b)    )=\l(\ref{limMM}\r) 
\]
\[ =\lim_{t\to\infty}\pp(    \D^{g(t)} \in (a,b)  | \cal{O}_t )  = \hat{\bb{Q}}(  \hat{T}_\infty \in (a,b)    ),
\]
from which it follows immediately that $M$ is transient under $\hat{\bb{Q}}(\cdot)$, as required.
\e{proof}

  \subsection{Proofs in the \textit{I}(\textit{f}) $=\infty$ Case }

  \n  The next three proofs require Lemma  \hbox{\ref{lemma4}},  proven in Section  \hbox{\ref{lemmasproofs}}.
     
 \b{lemma} \label{lemma4} 
In case  (\hbox{\hyperref[case1]{i}}),     for $t_0(y)$ as in~$\l(\ref{t0y}\r)$, uniformly in    $h  >  0, y  >  g(h)$, and $t\in (t_0(y),\infty]$,
 \b{equation} \label{l4}
\int_{t_0(y)}^t     \l(   \ov(g(s+h)-y) - \ov(g(s))     \r) ds     \lesssim   y f'(y) \ov(y)     .
 \e{equation}
 \e{lemma}

    \subsubsection{Proof of Theorem \mbox{\ref{recurrentthm}}}
  
  \b{proof}[Proof of Theorem  \hbox{\ref{recurrentthm}}] 
     For  fixed $h>0$, $y>g(h)$,  we  will   prove   that $    q_h(y):=      \lim_{t\to\infty}      \pp(       \cal{O}_{t-h}^{g_{y}^{h}}       ) / \pp(\cal{O}_t) < \infty$. For   each $h  >  0$, $y  >  g(h)$, note that $g(t)\sim g_{y}^{h}(t-h)$ by the properties of $g(t)$ given in Assumption \mbox{\ref{case1}}. Hence   $\ov(g(t)) \sim \ov(g_{y}^{h}(t-h))$ as $t\to\infty$, since  $\ov$ is \mbox{\hyperref[crv]{CRV}} at $\infty$. Thus, applying  Corollary \hbox{\ref{lemma1}},
  \b{equation} \label{firstlim}
  \lim_{t\to\infty}  \f{\pp\big(       \cal{O}_{t-h}^{g_{y}^{h}}       \big) }{ \pp(\cal{O}_t)}      =   \lim_{t\to\infty}      \f{ \ov( g_{y}^{h}(t-h)) \Phi_y^h(t)   }{   \ov( g(t)) \Phi(t) }    =   \lim_{t\to\infty}          \f{  \Phi_y^h(t)   }{    \Phi(t) }.
  \e{equation}
Then by $\l(\ref{soln1}\r)$ and $\l(\ref{soln2}\r)$, for   $t_0(y)$ as defined in (\ref{t0y}),  
  \b{align*}
\l( \ref{firstlim}\r)  \hspace{-0.05cm} &=    \hspace{-0.05cm}  \f{        \Phi_y^h(   t_0(y)  )    }{   \Phi(1)   }    \lim_{t\to\infty}     \exp\l(    \hspace{-0.05cm}  \int_{   t_0(y) }^{t}   \hspace{-0.2cm}     \l(  \ov(g_{y}^{h}(s))     \hspace{-0.05cm}  +  \hspace{-0.05cm}  \rho_y^h(s)   \r)   ds      
 -  \hspace{-0.15cm}    \int_1^{ t   }  \hspace{-0.1cm}  \l(  \ov(g(s)) +   \rho(s)   \r) ds   \hspace{-0.05cm}    \r)   \hspace{-0.1cm} .
\intertext{By $\l(\ref{1}\r)$ in Lemma  \hbox{\ref{lemma2}},   the integral $  \int_{   t_0(y) }^{\infty}      \rho_y^h(s) ds    $ is uniformly bounded for all $h>0$, $y>g(h)$.
By Lemma \mbox{\ref{lowerboundforrho}}, it follows that     $- \int_1^\infty \rho(s) ds <\infty$,~so}
  \l(\ref{firstlim} \r) &\lesssim    \f{        \Phi_y^h(   t_0(y)  )    }{   \Phi(1)   }  \lim_{t\to\infty}   \exp\l(   \int_{   t_0(y) }^{t}      \ov(g_{y}^{h}(s))     ds      
 -     \int_1^{ t   }  \ov(g(s))   ds     \r)
\\
  &\lesssim       \f{        \Phi_y^h(   t_0(y)  )    }{   \Phi(1)   }    \lim_{t\to\infty}     \exp\l(   \int_{   t_0(y) }^{t}   \l(   \ov(g_{y}^{h}(s))     -        \ov(g(s))   \r)  ds     \r).
\intertext{Applying Lemma  \hbox{\ref{lemma4}}, and recalling that $ y f'(y)  \ov(y)$  decreases to zero as $y\to\infty$,}
   \l(\ref{firstlim} \r) &\lesssim   \f{        \Phi_y^h(   t_0(y)  )    }{   \Phi(1)   }  \exp\l(  y f'(y) \ov(y)  \r)  <\infty.
  \end{align*}
  Now, $q_h(y) :=    \lim_{t\to\infty}      \pp(       \cal{O}_{t-h}^{g_{y}^{h}}       ) / \pp(\cal{O}_t) $ is non-decreasing in $y$ since for all $y<y'$, $g_{y}^{h}(t)= g(t+h) -y > g(t+h) -y'  = g_{y'}^{h}(t)$, and so   $ \pp(       \cal{O}_{t-h}^{g_{y}^{h}}       ) \leq     \pp(       \cal{O}_{t-h}^{g_{y'}^{h}}       )$. Finally, we conclude by  $\l(\ref{doobh}\r)$ that
  \(
  \mathbb{Q}\l( X_h \in dy    ;  \cal{B}_h  \r)         =       \pp\l(X_h \in dy   ;  \cal{B}_h \r) q_h(y),
  \) as required.

  \e{proof}

  \b{proof}[Proof of Proposition  \hbox{\ref{lemma5}}] For   $h>0$,       $    \lim_{t\to\infty}  \pp\l(    X_h \in (g(h),\infty) | \cal{O}_t    \r)=1$.
    We  will prove by dominated convergence  \cite[Theorem 1.21]{k06} that     limits and integration can be exchanged   from $\l(\ref{intt1}\r)$ to $\l(\ref{intt2}\r)$,   so by $\l(\ref{doobh}\r)$,~for~all~$h>0$,
  \b{align}          \nonumber
 1=  \lim_{t\to\infty}   \pp\l(    X_h \in (g(h),\infty) | \cal{O}_t    \r)  &= 
       \lim_{t\to\infty}     \int_{g(h) }^\infty     \pp\l(    X_h \in dy | \cal{O}_t    \r)  
      \\
       &=         \label{intt1}
          \lim_{t\to\infty}     \int_{g(h) }^\infty     \f{  \pp(    \cal{O}_{t-h}^{g_{y}^{h}} )   }{   \pp(\cal{O}_t)  }    \pp\l(    X_h \in dy ; \cal{O}_h    \r) 
  \\
    &=            \label{intt2}
          \int_{g(h) }^\infty     \lim_{t\to\infty}    \f{  \pp(    \cal{O}_{t-h}^{g_{y}^{h}} )   }{   \pp(\cal{O}_t)  }    \pp\l(    X_h \in dy ; \cal{O}_h    \r)   
   \\
       &=      \label{exchange}
           \mathbb{Q}\l(     X_h \in (g(h),\infty)    \r)  ,
  \end{align}
  
  \n as required. For $A  >  3\vee   (B-1)$, we will  bound the integral over $(g(h),\infty)$ via:
  \b{equation}   \label{I1I2I3}
     \l[  \f{ g(t-h)}{ A}  ,  \infty  \r)  \cup      (g(h) ,  g(h+1)] \cup \l( g(h+1) , \f{ g(t-h)}{ A}  \r)
  =:           I_1 \cup I_2   \cup  I_3.
  \e{equation}

    \p{Proof for \mbox{\hyperref[I1I2I3]{$I_1$}}}  Since $y\in  I_1 $ if and only if $t \leq  f(Ay ) + h$,   by  Corollary  \hbox{\ref{lemma1}},
  \b{align*}
 \int_{g(h)}^\infty     \bbm{1}_{\{ y\in I_1  \}}          \f{  \pp(    \cal{O}_{t-h}^{g_{y}^{h}} )   }{   \pp(\cal{O}_t)  }    \pp\l(    X_h \in dy ; \cal{O}_h    \r)     \lesssim
  \int_{g(h)}^\infty                \f{ \bbm{1}_{\{ y\in I_1  \}}  }{  \ov(g(t))   \Phi(t) }  \pp\l(    X_h \in dy     \r) 
=  \f{       \pp\l(    X_h \geq      \f{g(t-h)}{A}        \r)     }{  \ov(g(t))   \Phi(t) }        .
  \end{align*}
  Now,  $I(f)=\infty$, so   $\lim_{t\to\infty} \Phi(t) = \infty$ by $\l(\ref{IfPhi}\r)$, and it suffices  to show that
  \b{equation} \label{markov}
  \limsup_{t\to\infty}  \f{    \pp\l(    X_h \geq      \f{g(t-h)}{A}     \r)    }{   \ov(g(t))  }     
  \e{equation}
 is finite for each fixed $h>0$, as   the integal in $\l(\ref{intt1}\r)$ over the region $I_1$ tends to 0 as $t\to\infty$, so the dominated convergence theorem applies, trivially, on~$I_1$.

  Recall the notation in $\l(\ref{D}\r)$. Observe  that   $\D^{g(t)}$ has exponential distribution of~rate~$\ov(g(t))$,~so
  \b{align}    \label{partitioning}     \nonumber
    \pp      \l(     X_h \geq       \f{g(t-         h)}{A}          \r)    &=   
      \pp       \l(               X_h \geq        \f{g(t-         h)}{A}      ;   \D^{g(t)} >        h         \r)    
      +     \pp          \l(     X_h \geq        \f{g(t-         h)}{A}     ;   \D^{g(t)}       \leq h      \r) 
 \\\nonumber
  &\leq     \ \    \ \pp\l(    X_h^{(0,g(t))}  \geq    \f{g(t-h)}{A}        \r)    \ \ \ \    +       \pp\l(     \D^{g(t)} \leq h  \r) 
 \\\nonumber
  &=    \ \  \  \pp\l(    X_h^{(0,g(t))}  \geq       \f{g(t-h)}{A}       \r)   \ \ \ \  +      1 -  e^{ -  h    \ov(g(t))}
  \\
  &\leq  \ \   \   \pp\l(    X_h^{(0,g(t))}\geq        \f{g(t-h)}{A}      \r)     \ \ \ \   +      h    \ov(g(t)) ,
  \end{align}
  where $X^{(0,g(t))}$ has   the same L\'evy measure  as $X$, but restricted to $(0,g(t))$, so $X^{(0,g(t))}$ has no jumps larger than $g(t)$.
By $\l(\ref{partitioning}\r)$ and Markov's inequality, 
  \b{align*}
    \l(  \ref{markov}  \r)    -  h   
   &\leq    
       \limsup_{t\to\infty}  \f{    \pp\big(    X_h^{(0,g(t))}   \geq        \f{g(t-h)}{A}     \big)    }{   \ov(g(t))  }            
      \lesssim          \limsup_{t\to\infty}     \f{  A   \mathbb{E}[   X_h^{(0,g(t))} ]}{     \ov(g(t)) g(t  -  h)      } 
 %
      =          \limsup_{t\to\infty}    \f{ A h \int_0^{g(t)}   x \Pi(dx)      }{   \ov(g(t))    g(t  -  h)    }   .   
      \intertext{Now,   $ \int_0^{g(t)}   x \Pi(dx)     =   \int_{x=0}^{g(t)} \int_{y=0}^x dy    \Pi(dx)     =      \int_{y=0}^{g(t)}   \int_{x=y}^{g(t)} \Pi(dx) dy  \leq    \int_0^{g(t)}  \ov(y) dy     $. Then because $\ov$ is \mbox{\hyperref[rv]{regularly varying}} at $\infty$, by   Karamata's theorem \cite[Prop 1.5.8]{bgt89}, we deduce, as required for  dominated convergence  on $I_1$,~that~for~each $h>0$,  }
           \l(  \ref{markov}  \r)        &\lesssim      h     +  \limsup_{t\to\infty}   \f{  A h \int_0^1 \ov(y)dy +  Ah \int_1^{g(t)} \ov(y)dy  }{   \ov(g(t)) g(t-h)    }      \lesssim   h     +  \limsup_{t\to\infty}   \f{     g(t) \ov(g(t))      }{     \ov(g(t)) g(t-h)   }  
 <\infty.
 \end{align*}

    \p{Proof for \mbox{\hyperref[I1I2I3]{$I_2$}}}  
By  Theorem \mbox{\ref{recurrentthm}},  $q_h(y) $ is non-decreasing in $y$, so by $\l(\ref{qqq}\r)$,  $ \pp(       \cal{O}_{t-h}^{g_{y}^{h}}       )$ is non-decreasing in $y$.  Now,    $ \lim_{t\to\infty}     \pp(    \cal{O}_{t-h}^{g_{g(h+1)}^{h}} )  /  \pp(\cal{O}_t)          = q_h(g(h+1))<\infty$ for each fixed $h$ by     Theorem \mbox{\ref{recurrentthm}}, and we conclude that
  \b{align*}
  \lim_{t\to\infty}   \int_{g(h)}^\infty    \bbm{1}_{\{ y\in I_2  \}}         \f{  \pp(    \cal{O}_{t-h}^{g_{y}^{h}} )   }{   \pp(\cal{O}_t)  }    \pp\l(    X_h \in dy ; \cal{O}_h    \r)          &\leq    \lim_{t\to\infty}       \f{  \pp(    \cal{O}_{t-h}^{g_{g(h+1)}^h} )   }{   \pp(\cal{O}_t)  }        \int_{g(h)}^\infty    \bbm{1}_{\{ y\in I_2  \}}        \pp\l(    X_h \in dy ; \cal{O}_h    \r) 
\\
  &=   q_h(g(h+1))         \pp\l(    X_h \in I_2 ; \cal{O}_h    \r)  ,
  \end{align*}
which is finite for each   $h   >      0$, so   dominated~convergence applies~on~$I_2$.

    \p{Proof for \mbox{\hyperref[I1I2I3]{$I_3$}}}  By $\l(\ref{ODE2}\r)$ and Corollary \mbox{\ref{lemma1}}, for all large enough $t$,
 \b{align} \nonumber
 \f{  \pp(    \cal{O}_{t-h}^{g_{y}^{h}} )   }{   \pp(\cal{O}_t)  }  
  &\leq 2      \f{    \l[   \ov(g_{y}^{h}(t-h))   +     \rho_y^h(t-h)    \r] \Phi_y^h(t-h)             }{     \ov(g(t))  \Phi(t)}    .
  \intertext{For $y\in I_3$, $y>g(h+1)$, so $f(y)-h> f(g(h+1))-h =1$,  and $ 1+ 1/(f(y)-h)   <2 $.       By $\l(\ref{2}\r)$,  as $\lim_{t\to\infty} u(t)=0$, for all large enough $t$ and for all $y \in I_3$,    }
\nonumber
  \f{  \pp(    \cal{O}_{t-h}^{g_{y}^{h}} )   }{   \pp(\cal{O}_t)  }    &\leq 2      \f{ \l( 1 +  u(t)  (1 + \f{1}{f(y)-h}  )  \r) \ov(g_{y}^{h}(t-h))     \Phi_y^h(t-h)             }{     \ov(g(t))  \Phi(t)}       
  \\
  \nonumber
  \f{  \pp(    \cal{O}_{t-h}^{g_{y}^{h}} )   }{   \pp(\cal{O}_t)  }    &\leq 6      \f{     \ov(g_{y}^{h}(t-h))     \Phi_y^h(t-h)             }{     \ov(g(t))  \Phi(t)}       .
  \intertext{Now, $y<g(t-h)/A< g(t)/3$   for $y\in I_3$, so by $\l(\ref{gyh}\r)$, since $\ov$ is 
  regularly varying at $\infty$,        $\ov(g_{y}^{h}(t-h))= \ov(g(t)- y)  \lesssim \ov(g(t))$, uniformly among $y\in I_3$.   So   for each fixed $h>0$, for all large enough $t$, uniformly in $y \in I_3$,    }
    \label{uno}
    \f{  \pp(    \cal{O}_{t-h}^{g_{y}^{h}} )   }{   \pp(\cal{O}_t)  }      &\lesssim        \f{    \Phi_y^h(t-h)    }{   \Phi(t)}.  
 \intertext{Now, by $\l(\ref{soln2}\r)$, for $t_0(y)$ as defined in $\l(\ref{t0y}\r)$, we have}
\nonumber
\Phi_y^h(t-h) &= \Phi_y^h(t_0(y)) \exp\l(      \int_{t_0(y)}^{t-h}  \ov(g_{y}^{h}(s))ds  +     \int_{t_0(y)}^{t-h}   \rho_y^h(s)ds       \r).
\intertext{Applying    $\l(\ref{1}\r)$ and  recalling  that $ 1+ 1/(f(y)-h)   <2 $ for $y\in \mbox{\hyperref[I1I2I3]{$I_3$}}$,   the integral  $\int_{t_0(y)}^{\infty}   \rho_y^h(s)ds $ is uniformly bounded among $y$, so uniformly among $y\in I_3$,     }     
\nonumber
 \Phi_y^h(t-h) &\lesssim 
 \Phi_y^h(t_0(y)) \exp\l(     \int_{t_0(y)}^{t-h}  \ov(g_{y}^{h}(s))ds  \r)      .
 \intertext{By Lemma \mbox{\ref{lowerboundforrho}},       $ \liminf_{t\to\infty} \int_1^t   \rho(s) ds > - \infty$,  so~by~$\l(\ref{ODE1}\r)$,}   
\nonumber
 \Phi(t) &=   \Phi(1) \exp\l(    \int_1^t \ov(g(s))ds  +     \int_1^t \rho(s)ds  \r)     \gtrsim   \exp\l(    \int_1^t \ov(g(s))ds    \r).
 \intertext{Since $y> g(h+1)$ in $I_3$, recalling $\l(\ref{t0y}\r)$,   $t_0(y)   \geq  f(Ay) >  f(y) \geq h+1 >1$,~so
     }
\nonumber
 \f{ \Phi_y^h(t-h) }{    \Phi(t)}    &\lesssim  \Phi_y^h(t_0(y)) \exp  \hspace{-0.05cm}   \l( \hspace{-0.05cm}    \int_{t_0(y)}^{t-h} \hspace{-0.05cm}    \ov(g_{y}^{h}(s))ds     - \hspace{-0.1cm}    \int_1^t   \hspace{-0.05cm}  \ov(g(s))ds      \r)
\\ 
\nonumber
&\leq       \Phi_y^h(t_0(y) )  \exp\l(      \int_{t_0(y)}^{t}   \l(\ov(g_{y}^{h}(s))   - \ov(g(s))    \r) ds          \r).
\intertext{Now, by Lemma  \hbox{\ref{lemma4}}, since $  y  f'(y) \ov(y) $ decreases to zero as $y\to\infty$,  we have uniformly among $y > g(h)$,
\[
    \int_{t_0(y)}^{t}                    \l( \ov(g_{y}^{h}(s))                  -                    \ov(g(s))    \r)       ds                 \leq       \int_{t_0(y)}^{\infty}                    \l( \ov(g_{y}^{h}(s))                  -                    \ov(g(s))    \r)       ds                    \lesssim                                 \sup_{y>g(h)}         y  f'(y) \ov(y)                 <\infty  .
\]
 So    for each fixed $h>0$,  uniformly among $y > g(h)$,  by   $\l(\ref{phiyh}\r)$,     }
\\
    \lim_{t\to\infty}   \f{ \Phi_y^h(t-h) }{    \Phi(t)} 
  &\lesssim    \Phi_y^h(t_0(y))   e^{-  \int_1^{t_0(y)} \ov(g(s))ds}    \leq       t_0(y) e^{-  \int_1^{t_0(y)} \ov(g(s))ds}  .             \label{dos}
\end{align}
Now, it follows from $\l(\ref{uno}\r)$ and  $\l(\ref{dos}\r)$ that for each fixed $h>0$, uniformly for all $t$ large enough,
\b{equation}    \label{eqn}
           \int_{g(h)}^\infty    \bbm{1}_{\{ y \in I_3 \}}     \f{  \pp(    \cal{O}_{t-h}^{g_{y}^{h}} )   }{   \pp(\cal{O}_t)  }    \pp\l(    X_h \in dy ; \cal{O}_h \r)
 \lesssim  
           \int_{g(h)}^\infty    \bbm{1}_{\{ y \in I_3 \}}       t_0(y)   e^{ - \int_1^{t_0(y)} \ov(g(s))ds}  \pp\l(    X_h \in dy ; \cal{O}_h \r)   .    
\e{equation}
\b{align*}
%
\intertext{Now (by choice of $A$ sufficiently large if necessary) we have $t_0(y)=f(Ay) \vee (1 + 2/A)=f(Ay)$ for all $y>g(h+1)>g(1)$. Then writing $\zeta(x):=  x e^{-\int_1^{x} \ov(g(u))du}$, noting $\zeta(\cdot)$ is differentiable and $\zeta(1)=1$,  }
\l(   \ref{eqn}\r)    &=           \int_{g(h)}^\infty    \bbm{1}_{\{ y \in I_3 \}}     f(Ay) e^{-  \int_1^{f(Ay)} \ov(g(s))ds}    \pp\l(    X_h \in dy ; \cal{O}_h \r)        
\\
 &=           \int_{  g(h+1)}^{\infty}      f(Ay) e^{-  \int_1^{f(Ay)} \ov(g(s))ds}    \pp\l(    X_h \in dy ; \cal{O}_h \r)        
 \\
 &\leq           \int_{  g(h+1)}^{\infty}      f(Ay) e^{-  \int_1^{f(Ay)} \ov(g(s))ds}    \pp\l(    X_h \in dy  \r)     =           \int_{  g(h+1)}^{\infty}    \zeta(  f(Ay) )   \pp\l(    X_h \in dy  \r)  
  \\
 &=           \int_{  g(h+1)}^{\infty}  \l(  \int_1^{f(Ay)}  \zeta'(  x)dx + 1\r)   \pp\l(    X_h \in dy  \r)  
\\
   &=   \pp\l( X_h \geq g(h+1)\r) +    \int_{  g(h+1)}^{\infty}   \int_1^{f(Ay)}  \zeta'(  x)dx   \pp\l(    X_h \in dy  \r)  . 
 \end{align*}
Changing the order of integration,  and then applying the result that for each fixed $c,h>0$, $\pp(X_h \geq z) \asymp \ov(z)$, uniformly in $z>c$, since $\ov$ is \mbox{\hyperref[rv]{regularly varying}}, see   \cite[Theorem 1 (iii)]{egv79},  
\begin{align*}
    \int_{ 1}^{\infty}   \zeta'(  x)   \int_{ g(h+1)\vee \f{g(x)}{A}}^{\infty}  \hspace{-5pt}    \pp\l(    X_h \in dy  \r)    dx
&\leq   
  \int_{ 1}^{\infty}   \zeta'(  x)   \int_{ \f{g(x)}{A}}^{\infty}     \pp\l(    X_h \in dy  \r)    dx
 \\
 & =   
       \int_{ 1}^{\infty}   \zeta'(  x)     \pp\l(    X_h \geq \f{g(x)}{A}  \r)    dx
 \\
 &\asymp   
       \int_{ 1}^{\infty}   \zeta'(  x)  \ov\l(\f{g(x)}{A}\r)    dx
\lesssim    
   \int_{ 1}^{\infty}   \zeta'(  x)  \ov\l(g(x)\r)    dx
 \\
 &=    \int_{ 1}^{\infty}   \hspace{-5pt}   \big[ e^{-\int_1^{x} \ov(g(u))du}  -  x \ov(g(x))e^{-\int_1^{x} \ov(g(u))du} \big] \ov\l(g(x)\r)    dx   
 \\
 &\leq    \int_{ 1}^{\infty}  e^{-\int_1^{x} \ov(g(u))du}     \ov\l(g(x)\r)    dx
 \\
  &  = \int_1^\infty    \f{d}{dx}\l(     -e^{- \int_1^{x} \ov(g(u))du   }     \r) dx <\infty,
 \end{align*}
and therefore $\l(   \ref{eqn}\r) <\infty$, so  the dominated convergence theorem applies  on $I_3$, and the order of limits and integration can be swapped between $\l(\ref{intt1}\r)$ and $\l(\ref{intt2}\r)$, as required.

 \e{proof}

\b{proof}[Proof of Proposition  \hbox{\ref{recurrentmess}}]
     Recall    $\pi_h(X)$ is the sample  path of $X$ up to time $h$, $F_1$ is a  functional on the excursion process, 
and $\cal{B}\subseteq \cal{O}_h$, $\cal{B} \in \cal{F}_h$, where $(\cal{F}_u)_{u\geq0}$ is $X$'s   filtration.  Disintegrating on the value of $\pi_h(X)\in\cal{B}$,   
  \b{equation} \label{Erecccc}
  \mathbb{E}\l[           F_1(\pi_h((\eps_s)_{s\geq0}))     \   \bbm{1}_{ \{ \pi_h(X) \in \cal{B}  \}  }  \   |  \  \cal{O}_t  \ \r]   
  \e{equation}
  \[
  =    \underset{\nu\in\cal{B}}\int          \mathbb{E} \big[       F_1(\pi_h((\eps_s)_{s\geq0}))      \big|      \cal{O}_t   ;  \pi_h(X) = \nu          \big]      \pp\big(   \pi_h(X) \in d\nu            | \cal{O}_t   \big).
  \]
Given a fixed path  $\pi_h(X)= \nu$, $\pi_h((\eps_s)_{s\geq0})$ depends only on $\nu$, so  $\pi_h((\eps_s)_{s\geq0})$ is conditionally independent   of   $\cal{O}_t$, and then as $\lim_{t\to\infty} \pp(   \pi_h(X) \in d\nu            | \cal{O}_t   ) = \bb{Q}( \pi_h(X) \in d\nu   )$ by Theorem \mbox{\ref{recurrentthm}}, 
  \b{align} 
  \nonumber
 \lim_{t\to\infty} \l(\ref{Erecccc}\r)
  &\overset{\phantom{\ref{recurrentthm}}}=   \lim_{t\to\infty}
     \int_{\cal{B}}         \mathbb{E} \l[      F_1(\pi_h((\eps_s)_{s\geq0} ))       \big|    \pi_h(X) = \nu          \r]    \pp\l(   \pi_h(X) \in d\nu        | \cal{O}_t   \r)
  \\
  \label{pointingagainatyou}
  &\overset{\ref{recurrentthm}}= 
    \int_{\cal{B}}         \mathbb{E} \l[      F_1(\pi_h((\eps_s)_{s\geq0} ))       \big|    \pi_h(X) = \nu          \r] \bb{Q}( \pi_h(X) \in d\nu   )
    \\
    \nonumber
  &\overset{\phantom{\ref{recurrentthm}}}=\hspace{-1pt}:\hspace{1pt}
  \mathbb{E}_{\mathbb{Q}''}  \l[     F_1(\pi_h((\eps_s)_{s\geq0} ) )         \    \bbm{1}_{\{ \pi_h(X) \in \cal{B} \} }   \r],
  \end{align}
  where we can swap the order of limits and integration since $F_1$ is bounded. 
Taking $F_1\equiv 1$, it follows immediately that $\bb{Q}''(\cdot)$ and $\bb{Q}(\cdot)$ agree on $\cal{D}[0,\infty)$.
    The weak convergence as $t\to\infty$ of $((M_t)_{t\geq0}, (X_s)_{s \geq0})$ under  $\pp(\cdot|\cal{O}_t)$  to $((M_t)_{ t\geq0}, (X_s)_{ s \geq0})$ under $\mathbb{Q}''( \cdot )$   then follows immediately from the fact (see e.g.\  \cite[Ex.\ IV.6.3]{b98} or \cite[p4113]{ks13}) that for all $x>0$, $(M_t)_{t\geq0}$  is uniquely determined by $(\eps_s)_{s\geq0}$ and $(X_s)_{s\geq0}$, and both of $(\eps_s)_{s\geq0}$ and $(X_s)_{s\geq0}$ have weak limits as determined in $\l(\ref{pointingagainatyou}\r)$. That is, we construct $ M$ pathwise via its excursions using $(X_u)_{u\geq0}$ to determine the timing and length of each excursion,  where we sample the excursions of $ M$  using the excursion measure conditional on each excursion length.  The fact that $M$ is recurrent under $\bb{Q}''(\cdot)$ follows immediately from this construction, since $X$ does not explode to infinity by Proposition \mbox{\ref{lemma5}}.
   
  \e{proof}

  \subsubsection{Proof of Theorem \mbox{\ref{entrrep}}}

 \n  Lemmas \mbox{\ref{l6}}, \mbox{\ref{qh}}, \mbox{\ref{bigjump}}, \hbox{\ref{1fyh}},    proven in Sections \mbox{\ref{l6proof}}, \mbox{\ref{lemmasproofs}}, are needed for the proof of  Theorem \mbox{\ref{entrrep}}.

      \b{lemma}   \label{l6}      For each subordinator and $g=f^{-1}$   in case (\hbox{\hyperref[case1a]{ia}}),          there exist $K,h_0>0$ such that for all   $h> h_0$, with $\ov(y)  =  y^{-\alpha} L(y) $, uniformly in~$y  >  K$, 
  \b{equation}
  \label{l6statement}  
   \pp\l(   X_h \in g(h) dy ;   \cal{O}_h    \r)     \asymp     y^{-1-\alpha}    \f{L(g(h)y)}{L(g(h))}  \pp(\cal{O}_h) dy.
         \e{equation}
    \e{lemma}

 \b{lemma} \label{qh}    For $\d>0$ small enough that $0<f(0)<f(\d)<1$, uniformly for all $h>0$ and $y>g(h+f(\d))$,
  \b{equation*} \label{qbound}
   q_h(y)  \asymp     \Phi_y^h(f(Ay))    \exp\l(       -       \int_{
   1}^{f(Ay)} \ov(g(s))   ds \r).
   \e{equation*}
  \e{lemma}
  
     \b{lemma} \label{bigjump}      
   For a subordinator  and $g=f^{-1}$  in case (\hbox{\hyperref[case1a]{ia}}),   let $S_{\Delta_1^{g(h)}}$ denote the size of    its first jump of size greater than $g(h)$. Then there is $h_0>0$ such that uniformly in  $h>h_0$ and $v>1$,
   \[
   \pp\l(   S_{\Delta_1^{g(h)}} \in g(h)dv  \r)      =     \f{ \Pi(g(h) dv)   }{  \ov(g(h))   }   \asymp    \f{  L(g(h)v) }{L(g(h))}  v^{-1-\alpha} dv.
   \]
  In particular there is $x_0\in (0,\infty)$ so that for all $x>x_0$, with $\Pi(dx)  =  u(x)dx$,
   \begin{equation} \label{2.10.3}
         u(x)   \asymp  x^{-1} \ov(x)    =  L(x) x^{-1-\alpha}    .  
   \end{equation}
   \e{lemma}

  \b{lemma} \label{1fyh}  Recall the notation $\l(\ref{phiyh}\r),\l(\ref{t0y}\r)$. If  $h>0$,       $y>g(h)$, and $t\geq f(Ay)$,  for $A>3\vee (B-1)$,  then $\Phi_y^h(t) \geq  f(y)- h$.
  
  \e{lemma}


 \b{proof}[Proof of Theorem \ref{entrrep}] By Proposition \mbox{\ref{lemma5}},   $ \mathbb{Q}( X_h  \in(   g(h),\infty))  =  1$,   $h  \geq  0$. By Definition \mbox{\ref{entdefn}},
 \b{equation}
 \label{wiff}
 w\in R_g       \iff  \lim_{h\to\infty}  \mathbb{Q}\l(  X_h \in (g(h),   w(h) g(h))    \r)   =0  ,
 \e{equation}
  since  $\mathbb{Q}''(\cdot)$ and $\mathbb{Q}(\cdot)$ agree on the space $\cal{D}[0,\infty)$ containing $X$.  Now, by Theorem \mbox{\ref{recurrentthm}}, 
  \b{equation}    \label{integ}
  \lim_{h\to\infty}  \mathbb{Q}\l(  X_h \in (g(h),   w(h) g(h))    \r)   =    \lim_{h\to\infty} \int_{g(h)}^{w(h)g(h)}  q_h (y) \pp\l( X_h \in dy ; \cal{O}_h    \r) .
  \e{equation} 
  We begin by showing that if  $\lim_{h\to\infty}   \int_h^{  f( w(h) g(h) )  }     \ov(g(s)) ds =0$, then $w\in R_g$.   
   \p{Proof of Sufficient Condition}    Let $\lim_{h\to\infty}   \int_h^{  f( w(h) g(h) )  }     \ov(g(s)) ds =0$.    To show   $w\in R_g$, we will show  that the limit of the integral in $\l(\ref{integ}\r)$~is~zero on~each~of
   \[
     \label{R123}   [g(h),g(h+1)] \cup  [g(h+1),Kg(h)]   \cup    [g(h+1)\vee Kg(h),w(h)g(h)] =: R_1 \cup R_2 \cup R_3
     \]
       separately, where $K$ is the constant as in Lemma \mbox{\ref{l6}}.  
 Note that if $g(h+1)\! >\! K g(h)$,  we need only consider  $R_1 \cup R_3$. Since $g$ is non-decreasing, we only need to consider the value $K$ if $K>1$.

  \p{Proof for $R_1$}  By Theorem \mbox{\ref{recurrentthm}}, $q_h(y)$ is non-decreasing in $y$, so
  \b{align} 
  \label{h+1}
 &  \lim_{h\to\infty} \int_{g(h)}^{ g(h+1)}  q_h (y) \pp\l( X_h \in dy ; \cal{O}_h    \r) 
\\
 \nonumber
   & \leq     \lim_{h\to\infty}    q_h(g(h+1))  \int_{g(h)}^{ g(h+1)}  \pp\l( X_h \in dy ; \cal{O}_h    \r) 
  \leq      \lim_{h\to\infty} q_h(g(h+1))   \pp\l(  \cal{O}_h    \r) .
     \intertext{Applying Lemma \mbox{\ref{qh}}, then applying Corollary \mbox{\ref{lemma1}},}
     \nonumber
     \l(\ref{h+1}\r) &\lesssim     \lim_{h\to\infty}   \Phi_{g(h+1)}^h(f(Ag(h+1)))   e^{      -       \int_{1}^{f(Ag(h+1))} \ov(g(s))   ds      }        \pp\l(  \cal{O}_h    \r)
     \\
     \nonumber
     &=   \lim_{h\to\infty}   \Phi_{g(h+1)}^h(f(Ag(h+1)))  e^{      -       \int_{1}^{f(Ag(h+1))} \ov(g(s))   ds      }         \Phi(h)      \ov(g(h))   .
     \intertext{Now, $ \Phi_{g(h+1)}^h(f(Ag(h+1)))  \leq   f(Ag(h+1)) \lesssim    h$, and   $f(Ag(h+1)) \geq   h$, as  $f=g^{-1}$ is \mbox{\hyperref[orv]{$\cal{O}$-regularly varying}}, increasing, and $A>1$.  By $\l(\ref{phiasymp}\r)$ and $\l(\ref{condt}\r)$,}
     \nonumber
          \l(\ref{h+1}\r) &\lesssim  \lim_{h\to\infty}    h e^{   -       \int_{1}^{f(Ag(h+1))} \ov(g(s))   ds    }     \Phi(h)       \ov(g(h))   
  \\
  \nonumber
  &\leq       \lim_{h\to\infty}    h e^{   -       \int_{1}^{h} \ov(g(s))   ds    }     \Phi(h)       \ov(g(h))    
          \overset{\l(\ref{phiasymp}\r)}\lesssim    \lim_{h\to\infty}   h  \ov(g(h))  \overset{\l(\ref{condt}\r)}=0.
     \end{align}

                      \p{Proof for    \mbox{\hyperref[R123]{$R_2$}}}  Recall that we only need to consider $R_2$ when $g(h+1)< Kg(h)$, in which case $K$ must satisfy $K>1$.  By Theorem \mbox{\ref{recurrentthm}}, $q_h(y)$ is non-decreasing in $y$, so
  \b{align}
 \label{limint}
   &  \lim_{h\to\infty} \int_{g(h+1)}^{Kg(h)}  q_h (y   ) \pp\l( X_h \in dy ; \cal{O}_h    \r)
\\
 \nonumber
  &\leq    \lim_{h\to\infty}  q_h (Kg(h) )   \int_{g(h+1)}^{Kg(h)}   \pp\l( X_h \in dy ; \cal{O}_h    \r)
\leq     \lim_{h\to\infty}  q_h (Kg(h) )   \pp\l( \cal{O}_h    \r).
\intertext{Applying Corollary \mbox{\ref{lemma1}}, then   Lemma \mbox{\ref{qh}} 
 ($g(h+f(\d))<g(h+1)< Kg(h)$ for $\d$ in  Lemma \mbox{\ref{qh}}),
 }
 \nonumber  \l( \ref{limint}\r)    &\overset{\ref{lemma1}}\leq     \lim_{h\to\infty}  q_h (Kg(h) )       \Phi(h)\ov(g(h)).
\\
  \nonumber
   &\overset{\ref{qh}}\lesssim      \lim_{h\to\infty}  \Phi_{Kg(h)}^h ( f(AKg(h))  )    e^{   - \int_1^{ f(AKg(h))   }    \ov(g(s))ds  }     \Phi(h)\ov(g(h)).
%
\intertext{Observe that $\Phi_{Kg(h)}^h(t)\leq t$  for all $t>0$, by $\l(\ref{phiyh}\r)$, and  $f(AKg(h)) \geq f(g(h))=h$ since $f$ increasing and $A,K>1$.      Moreover, as $f$ is  $\cal{O}$-regularly varying at $\infty$,  $f(AKg(h)) \lesssim f(g(h))=h$, as $h\to\infty$,  so}
\nonumber
  \l( \ref{limint}\r)   &\leq  \lim_{h\to\infty}       f(AKg(h))     \exp\l(    - \int_1^{  f(AKg(h))  }    \ov(g(s))ds  \r)     \Phi(h)\ov(g(h)).
\\
 \nonumber
  &\leq  \lim_{h\to\infty}     f(AKg(h))  \exp\l(    - \int_1^{  h }    \ov(g(s))ds  \r)     \Phi(h)\ov(g(h))
  \\
  \nonumber
  &\lesssim  \lim_{h\to\infty}   h   \exp\l(    - \int_1^{  h }    \ov(g(s))ds  \r)     \Phi(h)\ov(g(h)).
\end{align}
By $\l(\ref{phiasymp}\r)$,  $\Phi(h)  \asymp   \exp(  \int_1^h \ov(g(s))ds )$ as $h\to\infty$, so as $\lim_{h\to\infty} h \ov(g(h))  =0$ by $\l(\ref{condt}\r)$, we conclude  that  $\l(\ref{limint}\r)=0$, so the integral over $R_2$ is zero, and  thus
\b{equation}      
\label{R1R2}
    \lim_{h\to\infty}     \int_{R_1\cup R_2}        q_h (y) \pp\l( X_h \in dy ; \cal{O}_h    \r)  = 0 .
 \e{equation}

\p{Proof for   \mbox{\hyperref[R123]{$R_3$}}}   Now we wish to show convergence to zero of 
  \b{equation} \label{kw}
    \int_{g(h+1)\vee Kg(h)}^{w(h)g(h)}  q_h (y) \pp\l( X_h \in dy ; \cal{O}_h    \r)    =      \int_{\f{g(h+1)}{g(h)}\vee K }^{w(h) }  q_h (g(h)v) \pp\l( X_h \in g(h)dv ; \cal{O}_h    \r) .
  \e{equation}
 Applying Lemma \mbox{\ref{l6}}, then changing variables back to $y= g(h)v$,   recalling that $\ov(g(h))= g(h)^{-\alpha} L(g(h))$,    then applying Corollary \mbox{\ref{lemma1}}, as $h\to\infty$,
  \b{align}
\nonumber
\l(\ref{kw}\r) &\overset{\ref{l6}}\asymp     \int_{\f{g(h+1)}{g(h)}\vee K  }^{w(h) }  q_h (g(h)v)  v^{-1-\alpha}  \f{ L(g(h)v)}{L(g(h))}   \pp(\cal{O}_h)    dv
\\
\nonumber
  &\overset{\phantom{\ref{l6}}}=    \f{ \pp\l(\cal{O}_h\r) g(h)^{\alpha}  }{ L(g(h))   }    \int_{\f{g(h+1)}{g(h)}\vee K  }^{w(h) }  q_h (g(h)v)  g(h)^{-\alpha}   v^{-1-\alpha} L(g(h)v) dv
\\
\nonumber
  &\overset{\phantom{\ref{l6}}}=      \f{ \pp\l(\cal{O}_h\r)  }{ \ov(g(h))   }  \hspace{-2.5pt}   \int_{g(h+1)\vee  K g(h)}^{w(h)g(h) }      q_h (y)    y^{-1-\alpha} L(y) dy
%
%
%
  \overset{\ref{lemma1}}\sim     \Phi(h)  \hspace{-5pt}  \int_{g(h+1)\vee  Kg(h) }^{w(h)g(h) }      q_h (y)    y^{-1-\alpha} L(y) dy.
\intertext{Changing variables \hspace{-2pt}  ($ u     =   Ay$), applying Lemma \mbox{\ref{qbound}} and then the uniform convergence theorem \cite[Theorem 1.2.1]{bgt89},   as $L$ is \mbox{\hyperref[sv]{slowly varying}} at $   \infty       $ and  $\Phi_y^h(f(u))    \leq    f(u)$,    it follows that as $h  \to  \infty,$}
\nonumber
   \l(\ref{kw}\r) &\overset{\ref{qbound}}\asymp       \Phi(h)   \int_{g(h+1)\vee  Kg(h) }^{w(h) g(h)}   \Phi_y^h(f(Ay))    \exp\l(       -       \int_{1}^{f(Ay)} \ov(g(s))   ds \r)   y^{-1-\alpha}L(y) dy
%
  %
%
%
\\
\nonumber
  &\overset{\phantom{\ref{qbound}}}\leq        \Phi(h)   \int_{ Kg(h) }^{w(h) g(h)}   \Phi_y^h(f(Ay))    \exp\l(       -       \int_{1}^{f(Ay)} \ov(g(s))   ds \r)   y^{-1-\alpha}L(y) dy
\\
\nonumber
  &\overset{\phantom{\ref{qbound}}}\lesssim       \Phi(h)   \int_{AKg(h) }^{Aw(h) g(h)}   \Phi_y^h(f(u))    \exp\l(       -       \int_{1}^{f(u)} \ov(g(s))   ds \r)   u^{-1-\alpha}L(u) du
\\
\nonumber
   &\overset{\phantom{\ref{qbound}}}\leq       \Phi(h)   \int_{AKg(h) }^{Aw(h) g(h)}  f(u)    \exp\l(       -       \int_{1}^{f(u)} \ov(g(s))   ds \r)   u^{-1-\alpha}L(u) du.
   \intertext{Since $A,K>1$,  we can split up the integral as follows, and we will deal with each term separately:}
 %
\nonumber
   \l(\ref{kw}\r)  &\lesssim     \Phi(h)   \int_{g(h) }^{w(h) g(h)}  f(u)    \exp\l(       -       \int_{1}^{f(u)} \ov(g(s))   ds \r)   u^{-1-\alpha}L(u) du
\\
\nonumber
  &+    \Phi(h)   \int_{w(h)g(h) }^{Aw(h) g(h)}  f(u)    \exp\l(       -       \int_{1}^{f(u)} \ov(g(s))   ds \r)   u^{-1-\alpha}L(u) du
\\
\label{J12}
  &=: J_1(h) +   J_2(h).
  \end{align}

  \p{Proof for $J_2(h)$}  \hspace{-0.3cm}  As     $f(u)  \hspace{-0.025cm}   \geq   \hspace{-0.025cm}    f(w(h)g(h))  \hspace{-0.025cm}    \geq   \hspace{-0.025cm}   f(g(h))  \hspace{-0.025cm}   =  \hspace{-0.025cm}   h$ for $u  \hspace{-0.025cm}   \geq  \hspace{-0.025cm}    w(h)g(h)$,~by~$\l(\ref{phiasymp}\r)$,  
 \b{align}
 \nonumber
  J_2(h) &\overset{\phantom{(\ref{phiasymp})}}=   \Phi(h)   \int_{w(h)g(h) }^{Aw(h) g(h)}  f(u)   e^{   -       \int_{1}^{f(u)} \ov(g(s))   ds  }     u^{-1-\alpha}L(u) du
 \\
 \nonumber
 &\overset{\phantom{(\ref{phiasymp})}}=       \Phi(h)   e^{     -       \int_{1}^{h} \ov(g(s))   ds }     \int_{w(h)g(h) }^{Aw(h) g(h)}  f(u)  e^{      -       \int_{h}^{f(u)} \ov(g(s))   ds  }    u^{-1-\alpha}L(u) du 
 \\
 \nonumber
  &\overset{(\ref{phiasymp})}\lesssim   \hspace{-5pt} \int_{w(h)g(h) }^{Aw(h) g(h)}    \hspace{-5pt}  f(u)  e^{      -       \int_{h}^{f(u)} \ov(g(s))   ds  }    u^{-1-\alpha}L(u) du 
 \overset{\phantom{(\ref{phiasymp})}}\leq  \hspace{-5pt}     \int_{w(h)g(h) }^{Aw(h) g(h)}    \hspace{-5pt} f(u)      u^{-1-\alpha}L(u) du.
  \intertext{Since $f$, $f'$ are \mbox{\hyperref[orv]{$\cal{O}$-regularly varying}} at $\infty$, one can verify that, uniformly for all sufficiently large $u$, $f(u)/u  \asymp  f'(u)$,  see   \cite[Prop 2.10.3]{bgt89}.  Recall  that in case (\hbox{\hyperref[case1]{i}}), $ u f'(u)    \ov(u)$ decreases to 0 as $u\to\infty$, so   as $h\to\infty$,}
\nonumber
    J_2(h) &\lesssim             \int_{w(h)g(h) }^{Aw(h) g(h)}  f'(u)      u^{-\alpha}L(u) du
  =     \int_{w(h)g(h) }^{Aw(h) g(h)}   u  f'(u)    \ov(u)  u^{-1}  du.
\\
\nonumber
  &\leq           w(h)g(h)  \      f'\big(w(h)g(h)\big)    \    \ov\big(w(h)g(h)\big)       \int_{w(h)g(h) }^{Aw(h) g(h)} u^{-1} du
\\
\label{j2}
  &= o(1) \times \int_{w(h)g(h) }^{Aw(h) g(h)} u^{-1} du = o(1) \times \log(A)= o(1),
  \end{align}
so $\lim_{h\to\infty} J_2 (h)  =0$, and $J_2(h)$ never contributes. Now we consider   $J_1(h)$.

  \p{Proof for \mbox{\hyperref[J12]{$J_1(h)$}}} First, changing variables from $s$ to $ v:=g(s)$, so that $s=f(v)$,  
  \b{align}
\nonumber
  J_1(h) &=      \Phi(h)   \int_{g(h) }^{w(h) g(h)}  f(u)    \exp\l(       -       \int_{1}^{f(u)} \ov(g(s))   ds \r)   u^{-1-\alpha}L(u) du
\\
\nonumber
  &=     \Phi(h)   \int_{g(h) }^{w(h) g(h)}  f(u)    \exp\l(       -       \int_{g(1)}^{u} \ov(v)  f'(v) dv \r)   u^{-1-\alpha}L(u) du.
  \intertext{Recall $u^{-\alpha} L(u)   =   \ov(u)$, and $f(u)   \asymp     u  f'(u)$ uniformly as $u  \to  \infty$,  so as $h  \to  \infty$, }
\nonumber
    J_1(h) &\asymp     \Phi(h)   \int_{g(h) }^{w(h) g(h)}  f'(u)    \exp\l(       -       \int_{g(1)}^{u} \ov(v)  f'(v) dv \r)   u^{-\alpha}L(u) du
\\
\nonumber
  &=       \Phi(h)   \int_{g(h) }^{w(h) g(h)}  f'(u) \ov(u)    \exp\l(       -       \int_{g(1)}^{u} \ov(v)  f'(v) dv \r)   du.
%
 \intertext{Changing variables from $u$ to $z:= e^{-\int_{g(1)}^u \ov(v) f'(v) dv}$ and applying $\l(\ref{phiasymp}\r)$, it follows that as $h\to\infty$,   }
\nonumber
    J_1(h)   &\overset{\phantom{\l(\ref{phiasymp}\r)}}\asymp     
  \Phi(h)  \l[   e^{-\int_{g(1)}^{g(h)}     \ov(v)f'(v) dv  }   -  e^{-\int_{g(1)}^{w(h)g(h)}   \ov(v)f'(v) dv  }        \r]
\\
  &\overset{\phantom{\l(\ref{phiasymp}\r)}}=       \Phi(h)  \l[   e^{-\int_{1}^{h}     \ov(g(s))ds  }   -  e^{-\int_{1}^{f(w(h)g(h))}   \ov(g(s))  ds }        \r]
%
%
     \label{j1} 
     \overset{\l(\ref{phiasymp}\r)}\asymp         1-  e^{   -\int_{h}^{f(w(h)g(h))}   \ov(g(s))  ds }.
  \end{align}
    Thus by $\l(\ref{R1R2}\r)$, $\l(\ref{j2}\r)$ and $\l(\ref{j1}\r)$,    whenever   $\lim_{h\to\infty}\int_{h}^{f(w(h)g(h))}   \ov(g(s))  ds  =0$,  $w$ is in   the entropic repulsion envelope  $R_g$, as required for the sufficient condition.

  Now we will prove that if $w\in R_g$, then $ \lim_{h\to\infty}   \int_h^{f(w(h)g(h))}   \ov(g(s)) ds =0$.

 \p{Proof of Necessary Condition}  
 Let $w\in R_g$. Then by $\l(\ref{qqq}\r)$,     
 \b{align}
 %
 \nonumber
  0 &\overset{\phantom{\ref{l6}}}=   \lim_{h\to\infty}  \mathbb{Q}\l(  X_h \in (g(h),   w(h) g(h))    \r)
 \\
  \nonumber
  &
  \overset{\phantom{\ref{l6}}}=    \lim_{h\to\infty}      \int_{ g(h)}^{ w(h) g(h)}   q_h(  y)    \pp\l(    X_h \in     dy ; \cal{O}_h   \r)  
  %
  \overset{\phantom{\ref{l6}}}=    \lim_{h\to\infty}      \int_{ K  g(h) }^{ w(h)  g(h) }   q_h(  y)    \pp\l(    X_h \in    dy ; \cal{O}_h   \r)    ,
    \intertext{since the limit of  the integral over  $R_1 \cup R_2=(g(h),Kg(h))$ is always zero by $\l(\ref{R1R2}\r)$, regardless of $\lim_{h\to\infty}   \int_h^{  f( w(h) g(h) )  }     \ov(g(s)) ds =0$. Changing variables to $v= y / g(h) $ and  applying Lemma \mbox{\ref{l6}},      }
      \nonumber
0     &\overset{\phantom{\ref{l6}}}=    \hspace{-4pt}    \lim_{h\to\infty}      \int_{ K    }^{ w(h)    }  \hspace{-9pt}   q_h(  g(h)v)    \pp\l(    X_h \hspace{-2pt} \in \hspace{-2pt}  g(h)dv  ; \cal{O}_h   \r)       
   \overset{\ref{l6}}=    \lim_{h\to\infty}      \pp(\cal{O}_h)     \hspace{-4pt}    \int_{K   }^{w(h)  }  \hspace{-9pt}  q_h( g(h)v )        v^{-1-\alpha}  \f{ L(g(h)v) }{ L(g(h))}  dv.
   \intertext{Changing variables   to $y=g(h)v$ and recallling that $\ov(g(h))= g(h)^{-\alpha} L(g(h))$, then applying Corollary  \mbox{\ref{lemma1}}, and Lemmas \mbox{\ref{qh}} and \mbox{\ref{1fyh}},    }  
%
%
\nonumber
0 &\overset{\phantom{\ref{1fyh}}}=      \lim_{h\to\infty}    \f{ \pp(\cal{O}_h) }{ \ov(g(h))}    \int_{K g(h)}^{w(h)g(h)}     \hspace{-7pt} q_h( y )       y^{-1-\alpha}  L(y)      dy  
  \overset{\hspace{2pt} \ref{lemma1} \hspace{2pt} }=      \lim_{h\to\infty}      \Phi(h)       \int_{K g(h)}^{w(h)g(h)}      \hspace{-5pt} q_h( y )               y^{-1-\alpha}  L(y)   dy 
\\
\nonumber
 &\overset{  \hspace{2pt}  \ref{qh}  \hspace{2pt}  }=      \lim_{h\to\infty}   \Phi(h)      \int_{K g(h)}^{w(h)g(h)}   \Phi_y^h(f(Ay))   e^{ - \int_1^{f(Ay)}  \ov(g(s))ds   }               y^{-1-\alpha}  L(y)     dy
\\
\nonumber
 &\overset{\ref{1fyh}}\geq      \lim_{h\to\infty} \Bigg[   \Phi(h)      \int_{K g(h)}^{w(h)g(h)}   f(y)  e^{-   \int_1^{f(Ay)}   \ov(g(s))ds   }    y^{-1-\alpha}  L(y) dy
\\
\nonumber
  &\hspace{1.75cm}    -     h \Phi(h)      \int_{K g(h)}^{w(h)g(h)}     e^{-   \int_1^{f(Ay)}   \ov(g(s))ds   }    y^{-1-\alpha}  L(y) dy   \Bigg]
\\
\nonumber
   &\overset{ \hspace{2pt}  \phantom{\ref{1fyh}}}{=:}                  \lim_{h\to\infty} \l[ I_1 - I_2 \r].
  \end{align}
  First we consider $I_2$.  Note that $AK>1$, so  $f(Ay) \geq f(AKg(h)) \geq h$ for all $ y\geq    K g(h)$. Then since $\Phi(h) \asymp  e^{-   \int_1^{h}   \ov(g(s))ds   } $  by  $\l(\ref{phiasymp}\r)$, 
  \b{align*}
   \lim_{h\to\infty}    \l| I_2  \r|    &\overset{\phantom{\l(\ref{phiasymp}\r)}}\leq  \hspace{-2pt}   \lim_{h\to\infty} \hspace{-2pt}  h \Phi(h)  \hspace{-3pt}     \int_{K g(h)}^{w(h)g(h)}  \hspace{-7pt}   e^{-   \int_1^{h}   \ov(g(s))ds   }    y^{-1-\alpha}  L(y) dy \hspace{-2pt}
  \lesssim   \hspace{-2pt}    \lim_{h\to\infty}\hspace{-2pt}  h  \hspace{-4pt}   \int_{K g(h)}^{w(h)g(h)}  \hspace{-7pt}       y^{-1-\alpha}  L(y) dy.
  \intertext{By Lemma \mbox{\ref{bigjump}},    $   y^{-1-\alpha} L(y)dy     \asymp    \Pi(dy)$, so as $\ov$ is \mbox{\hyperref[rv]{regularly varying}} at $\infty$,~by~$\l(\ref{condt}\r)$,  }
%
%
 %
%
%
%
     \lim_{h\to\infty}    \l| I_2  \r|   &\overset{\ref{bigjump}}\lesssim        \lim_{h\to\infty}  h     \int_{K g(h)}^{w(h)g(h)}    \Pi(dy) 
    \overset{\phantom{\ref{bigjump}}}\leq     \lim_{h\to\infty}  h   \ov(Kg(h))    \lesssim    \lim_{h\to\infty}  h \ov(g(h)) \overset{\l(\ref{condt}\r)}=0,
  \end{align*}
    so $I_2=0$, and thus $   \lim_{h\to\infty}  I_1 \leq 0$. As $I_1$ is   non-negative,    $ \lim_{h\to\infty}  I_1 = 0$. Now, changing variables  to $v:=Ay$, as $f$ is \mbox{\hyperref[orv]{$\cal{O}$-regularly varying}}  at $\infty$ and $L$     is \mbox{\hyperref[sv]{slowly varying}} at $\infty$,    by the uniform convergence theorem~\cite[Theorem~1.2.1]{bgt89}, 
\b{align*}
0 =  \lim_{h\to\infty}      I_1
&=    \lim_{h\to\infty}   \Phi(h)      \int_{AK g(h)}^{Aw(h)g(h)}   f\l(\f{v}{A}\r)  e^{-   \int_1^{f(v)}   \ov(g(s))ds   }    v^{-1-\alpha} A^{\alpha}  L\l(\f{v}{A}\r) dv
  \\
&\asymp   \lim_{h\to\infty}   \Phi(h)      \int_{AK g(h)}^{Aw(h)g(h)}   f(v)  e^{-   \int_1^{f(v)}   \ov(g(s))ds   }    v^{-1-\alpha}  L(v) dv.
\intertext{Recall    $v^{-\alpha} L(v) = \ov(v)$, and     $f(v) \asymp   v f'(v)$ uniformly for all  large enough $v$, because   $f$, $f'$ are  $\cal{O}$-regularly varying at $\infty$, see   \cite[Prop 2.10.3]{bgt89}. Then}
   0  &=        \lim_{h\to\infty}   \Phi(h)      \int_{AK g(h)}^{Aw(h)g(h)}   f(v)    e^{-   \int_1^{f(v)}   \ov(g(s))ds   }      v^{-1}  \ov(v) dv
 \\
 &\asymp      \lim_{h\to\infty}   \Phi(h)      \int_{AK g(h)}^{Aw(h)g(h)}     f'(v) \ov(v)    e^{-   \int_1^{f(v)}   \ov(g(s))ds   }     dv.
  \intertext{Now,   one can verify  that $P(v):= \int_{g(1)}^v \ov(u) f'(u)du =  \int_1^{f(v)}  \ov(g(s))ds$ by changing variables from $u$ to $s  =  f(u)$. Then as   $A  >  3$ and $P'(v)  =    \ov(v) f'(v)  \geq~\hspace{-5pt}0$, }
 0  &=          \lim_{h\to\infty}   \Phi(h)     \int_{  AK g(h)}^{A w(h)g(h)}   P'(v)  e^{-  P(v)  }   dv 
  \\
       &\geq          \lim_{h\to\infty}    \l[   \Phi(h)  \hspace{-3pt}   \int_{  g(h)}^{w(h)g(h)}    \hspace{-8pt}  P'(v)  e^{-  P(v)  }   dv     -    \Phi(h)    \hspace{-3pt}   \int_{g(h)}^{AK g(h)}    \hspace{-8pt}  P'(v)  e^{-  P(v)  }   dv   \r]
 \\
 &=:  \lim_{h\to\infty} \l[K_1 -K_2  \r] . 
  \end{align*}
  Now, recall that by $\l(\ref{j1}\r)$,   as $h\to\infty$, we have
  \b{equation} \label{Kexp}
       K_1 \asymp  J_1 \asymp   \l( 1- e^{  -      \int_h^{f(w(h)g(h))}   \ov(g(s))ds}   \r) .
    \e{equation}
So if we  prove $\lim_{h\to\infty} K_1 =0$, then    $ \lim_{h\to\infty}   \int_h^{f(w(h)g(h))}   \ov(g(s))ds =0$,
 and the proof is complete. 
  As $K_1$ is always non-negative,  it suffices to   prove that $\lim_{h\to\infty}     K_1    \leq       0$. 
  To prove this, we will show  that   $\lim_{h\to\infty}   | K_2|        =     0 $.
  Since $g=f^{-1}$, note   $ f(v) > h$ for $v>g(h)$, then as $\Phi(h) \asymp e^{\int_1^h \ov(g(s))ds}$ by $\l(\ref{phiasymp}\r)$, 
  \b{align*}
 \lim_{h\to\infty}   |K_2|    &\overset{\phantom{\l(\ref{phiasymp}\r)}}=     \lim_{h\to\infty}   \Phi(h)  \int_{g(h)}^{AK g(h)}    \ov(v) f'(v)   e^{-   \int_1^{f(v)}   \ov(g(s))ds   }      dv
\\
  &\overset{\l(\ref{phiasymp}\r)}\asymp     \lim_{h\to\infty}        \int_{g(h)}^{AK g(h)}    \hspace{-5pt}    \ov(v) f'(v)  e^{-   \int_h^{f(v)}   \ov(g(s))ds   }     dv
    \leq        \lim_{h\to\infty}           \int_{g(h)}^{AK g(h)}   \hspace{-5pt} v   f'(v)  \ov(v)   v^{-1}     dv.
   \intertext{Recall that by assumption, $v f'(v) \ov(v)$ decreases to 0 as $v\to\infty$, and hence}
    \lim_{h\to\infty}   |K_2|   &\overset{\phantom{\l(\ref{phiasymp}\r)}}\lesssim  \lim_{h\to\infty}  g(h)  \  f'(g(h))  \    \ov(g(h))      \int_{g(h)}^{AKg(h)}   v^{-1}dv
   \\
     &\overset{\phantom{\l(\ref{phiasymp}\r)}}= \lim_{h\to\infty}      g(h)  \  f'(g(h))  \    \ov(g(h))    \times \log(AK) =0.
  \end{align*}


    \end{proof}

\subsubsection{Proof of Corollary \mbox{\ref{entrrepstable}}}
     
     \b{proof}[Proof of Corollary \mbox{\ref{entrrepstable}}]
    We need to verify that   a stable subordinator of index $\alpha\in(0,1)$ satisfies  $\l(\ref{l6cond}\r)$,   so Theorem \mbox{\ref{entrrep}} applies.    For $t>0$ and $x>g(t)+x_0$, by the scaling property of stable subordinators (see \cite[p227]{b98}),
 \b{equation} \label{corollproofbit}
  f_t(x) =    t^{-\f{1}{\alpha}}   f_1\l( \f{x}{t^{\f{1}{\alpha}}}   \r).
 \e{equation}
     
     \n Now consider the result (see  \cite[Theorem 1.12]{n18})  that for a stable subordinator of index $\alpha\in(0,1)$, $f_1(v) \sim c_\alpha  v^{-1-\alpha}$ as $v\to\infty$, for  $c_\alpha>0$ constant. In~particular, for all large enough $v$, $f_1(v)$ is arbitrarily close to $c_\alpha v^{-1-\alpha}$. Taking e.g.\ $a'=2c_\alpha$, it follows that there exist    $a',C \in  (0,\infty)$   such that for all $v>C$, $    f_1(v) \leq a' v^{-1-\alpha}$.
     
    As $\Pi(dv)=u(v)dv = c v^{-1-\alpha}dv$ for a constant $c>0$,   if we can   show that  $x/ t^{1/\alpha}  \geq C$  for all $t>0$, $x>g(t)+x_0$, with a suitable choice of  $x_0>0$,~then~by~$\l(\ref{corollproofbit}\r)$,
      \[
       f_t(x)   =  t^{-\f{1}{\alpha}}   f_1\l( \f{x}{t^{\f{1}{\alpha}}}   \r)    \leq  a'c t  x^{-1-\alpha}  = a t u(x),
      \]
  for $a:=a'c$,   so   condition $\l(\ref{l6cond}\r)$ is satisfied, and the proof will be complete.
 Indeed, by $\l(\ref{condt}\r)$, we have $\lim_{t\to\infty} t\ov(g(t))= \lim_{t\to\infty} tg(t)^{-\alpha}  =0$,     so there exists  $D\in(0,\infty)$    such that  for all $t>D$,        $tg(t)^{-\alpha}   \leq   C^{-\alpha}$,   so    $t^{1/\alpha} \leq  C^{-1} g(t)$, and hence for all $t>D$,
\[
\f{x}{t^{\f{1}{\alpha}} }   \geq   \f{g(t) + x_0}{t^{\f{1}{\alpha}}}  \geq   \f{g(t)}{t^{\f{1}{\alpha}}}  \geq  C   .  
 \]

 \n On the other hand, if $t\leq D$, then
 $  x/t^{1/\alpha} \geq  (g(t)+  x_0 ) /  D^{1/\alpha}   \geq   x_0  /  D^{1/\alpha}  $, and
 (choosing $x_0$ large enough that 
$x_0 / D^{1/\alpha} > C$ if necessary),
 we conclude that  $x/t^{1/\alpha}>C$   for all      $t>0$, $ x> g(t)+x_0$, so    \cite[Theorem 1.12]{n18}  applies to $\l(\ref{corollproofbit}\r)$. It follows that   condition   $\l(\ref{l6cond}\r)$ is satisfied, and so  Theorem \mbox{\ref{entrrep}} applies,~as~required.

     \e{proof}

%

\section{Proof of Lemma \hbox{\ref{lemma2}}}   \label{lemma2proof}

 In order to prove  Lemma \mbox{\ref{lemma2}}, we require  
  Lemmas  \hbox{\ref{keylemma}} and \hbox{\ref{1.1}}.  The proof of Lemma \mbox{\ref{keylemma}} is  provided in Section~\hbox{\ref{lemmasproofs}}, whereas the proof of Lemma \mbox{\ref{1.1}} is provided immediately below.

   \b{lemma} \label{keylemma}
   Let $(X_t)_{t\geq0}$ be a subordinator satisfying the assumptions in case (\hbox{\hyperref[case1]{i}}) or (\hbox{\hyperref[case2]{ii}}). Then there exists a constant $C>0$, which depends only on the law of $X$, such that for all $t>0$, $A(t)\in(1,\infty)$,    $B(t)>0$,  and  $ H(t) \in (0,1) $,    
       \b{equation}   \label{keyeqn}
   \pp\l( X_t^{ (0,A(t))}    >  B(t) \r) \leq    \exp  \l(    C  t \log\l(   \f{1}{H(t)}   \r)   H(t)^{-\f{A(t)}{B(t)}}    \ov(A(t)) \f{A(t)}{B(t)}   \r)    H(t).
   \e{equation}
                 \e{lemma}    
   

   \b{lemma} \label{1.1}   Recalling $\l(\ref{gyh}\r)$ and $\l(\ref{t0y}\r)$, if $t   >   t_0(y)$, $h\geq0$, and $y\geq0$,   then $g_{y}^{h}(t)     \geq      \l(    1    -    1/A    \r)g(t)$.
  
  \e{lemma}

     \b{proof}[Proof of Lemma \hbox{\ref{1.1}}]  By   $\l(\ref{t0y}\r)$, $t  >   t_0(y)  \geq   f(Ay).$  The result holds trivially when $y=0$. Otherwise, as  $g  =  f^{-1}$~is~increasing,
  \b{align*}
  g_{y}^{h}(t)  
   =     \l(    \f{g(t+h)}{g(t)}   - \f{y}{g(t)}     \r)g(t)  
   \geq      \l(    1  - \f{y}{g(t)}     \r)g(t)  
   \geq      \l(    1  - \f{y}{g(f(Ay))}     \r)g(t) 
   =     (    1  - A^{-1}    )g(t).
  \end{align*}
\e{proof}

  \n Now, in addition to showing that for each $A>3$, the inequalities $\l(\ref{1}\r)$ and $\l(\ref{2}\r)$  hold uniformly in  $h>0$, $y>g(h)$, and  $ t> t_0(y)$, we shall show that the inequalities $\l(\ref{1}\r)$ and $\l(\ref{2}\r)$  hold when $y=h=0$ for $t>t_0(0)>0$,  with $\rho(t)$ in place of $\rho_y^h(t)$.  For brevity, let us introduce the following notation:  
  \begin{equation}
  \label{calS}
  \cal{S} := \l\{  (h,y,t)\in \bb{R}^3 : h>0, y>g(h), t>t_0(y)    \r\} \cup \l\{ (0,0,t)\in\bb{R}^3 :  t>t_0(0)  \r\}.
  \end{equation}
 \n  
  Lemma \hbox{\ref{lemma2}} shall be proven by  splitting up $\rho_y^h(t)$ into   smaller pieces, and then showing that the inequalities $\l(\ref{1}\r)$ and $\l(\ref{2}\r)$ hold for each piece separately. 

  \b{proof}[Proof of Lemma  \hbox{\ref{lemma2}}] Lemma \mbox{\ref{lemma2}} 
is simpler to prove in 
case (\hbox{\hyperref[case2]{ii}}) than   case (\hbox{\hyperref[case1]{i}}) thanks to the     condition $\l(\ref{cond2}\r)$. We  thus omit the   proof  in~case~(\hbox{\hyperref[case2]{ii}}).   
\n Firstly, since $g      =     f^{-1}$ is continuous (and hence so is $g_y^h$), for all $h,y\geq0$ and $x,t  >    0$,
  \b{equation}
  \label{s-}
  \pp\big(\cal{O}^{g_y^h}_{t, X^{(0,x)}}  \big)     =      \pp\big(\cal{O}^{g_y^h}_{t-, X^{(0,x)}}\big)   ,
  \end{equation} 
 where  $\cal{O}^{g_y^h}_{t-, X^{(0,x)}} := \bigcap_{u<t} \cal{O}^{g_y^h}_{u, X^{(0,x)}} $. Moreover, by   $\l(\ref{truncoo}\r)$,  for all $h,y\geq0$, $x>0$, and $t\geq s  >    0$,
    \b{equation}
  \label{subevent}
  \pp(\cal{O}^{g_y^h}_{t, X^{(0,x)}}  )     \leq   \pp\big( \cal{O}^{g_y^h}_{s,X^{(0,x)}}  ; X_t^{(0,x)} \geq g_y^h(t)  \big)    \leq \pp\big( X_t^{(0,x)} > g_y^h(t) \big).
  \e{equation}  
  Now, we partition and disintegrate on the value of $\Delta_1^{g_y^h(t)}$, which is exponentially distributed with rate $\ov(g_y^h(t))$. It follows by $\l(\ref{ogyh}\r)$, $\l(\ref{s-}\r)$ and $\l(\ref{truncoo}\r)$, since $s<t$,  that  for all $h,y\geq0$ and $t  >    0$,
  \begin{equation}
  \label{ttos}  \pp( \cal{O}^{g_y^h}_{t} |  \Delta_1^{g_y^h(t)}  \! = \! s  ) \!  =\!    \pp( \cal{O}^{g_y^h}_{s} |  \Delta_1^{g_y^h(t)}   \!= \! s  )  \! =\!     \pp( \cal{O}^{g_y^h}_{s-} |  \Delta_1^{g_y^h(t)}  \! = \! s  )  \! = \!  \pp(\cal{O}^{g_y^h}_{s-,X^{(0,g_y^h(t))}}    )\!  =  \! \pp(\cal{O}^{g_y^h}_{s,X^{(0,g_y^h(t))}}    )  ,
  \end{equation} where the last equality holds since the small and large jumps are independent, so it follows that 
  \begin{align}
  \nonumber
  \pp\big( \cal{O}^{g_y^h}_{t}  \big) &=   \pp\big( \cal{O}^{g_y^h}_{t} ;  \Delta_1^{g_y^h(t)} \leq t  \big)  +   \pp\big( \cal{O}^{g_y^h}_{t}  ;  \Delta_1^{g_y^h(t)} >t   \big) 
%
 %
 %
 \\
  \nonumber
  &=  \ov(g_y^h(t))    \int_0^t      \pp( \cal{O}^{g_y^h}_{t} |  \Delta_1^{g_y^h(t)} =s  )   e^{ -  \ov(g_y^h(t)) s} ds +  \pp( \cal{O}^{g_y^h}_{t}  ;  \Delta_1^{g_y^h(t)} >t   ) 
 \\
  \label{oneone}
  &=   \ov(g_y^h(t))    \int_0^t   \pp ( \   \cal{O}^{g_y^h}_{s,X^{(0,g_y^h(t))}}    \    )   \  e^{ -  \ov(g_y^h(t)) s}   ds   +   \pp( \cal{O}^{g_y^h}_{t}  ;  \Delta_1^{g_y^h(t)} >t   )  .
  %
  %
\intertext{Now, observe    that by the definition   $\l(\ref{truncoo}\r)$, $  \pp( \cal{O}^{g_y^h}_{s} |  \Delta_1^{g_y^h(t)} >s  ) = \pp(\cal{O}^{g_y^h}_{s,X^{(0,g_y^h(t))}}    )  $, so for all $h,y\geq0$,}  
  \nonumber
\pp(\cal{O}^{g_y^h}_s )
&= \pp(   \cal{O}^{g_y^h}_s ; \Delta_1^{g_y^h(t)} \leq s ) +  \pp(   \cal{O}^{g_y^h}_s ; \Delta_1^{g_y^h(t)} > s ) 
\\
  \nonumber
&=  \pp(   \cal{O}^{g_y^h}_s ; \Delta_1^{g_y^h(t)} \leq s ) +  \pp(   \cal{O}^{g_y^h}_s | \Delta_1^{g_y^h(t)} > s ) \pp(\Delta_1^{g_y^h(t)} > s ) 
\\
  \label{twotwo}
&= \pp(   \cal{O}^{g_y^h}_s ; \Delta_1^{g_y^h(t)} \leq s ) + \pp\l(\cal{O}^{g_y^h}_{s,X^{(0,g_y^h(t))}}    \r)   e^{ -  \ov(g_y^h(t)) s}.
\intertext{Disintegrating  on  $\Delta_1^{g_y^h(t)} $,   recalling the notation introduced in $\l(\ref{D}\r)$ and $\l(\ref{phiyh}\r)$, by  $\l(\ref{oneone}\r)$, $\l(\ref{twotwo}\r)$,~and~$\l(\ref{s-}\r)$,}
%
\nonumber
  \pp\l( \cal{O}^{g_y^h}_{t}  \r)   &\overset{\phantom{\l(\ref{s-}\r)}}= \ov(g_y^h(t)) \int_0^t    \big[   \mathbb{P}   (\cal{O}^{g_y^h}_{s}) - \pp(   \cal{O}^{g_y^h}_s ; \Delta_1^{g_y^h(t)} \leq s ) \big] ds    +  \pp( \cal{O}^{g_y^h}_{t}  ;  \Delta_1^{g_y^h(t)} >t   )  
\\
  \label{final123}
    &\overset{\l(\ref{s-}\r)}= \hspace{-1pt} \ov(g_y^h(t))   \Phi(t)
\hspace{-1pt}-\hspace{-1pt} \ov(g_y^h(t))^2   \int_0^t      \int_0^s    \hspace{-1pt}   \mathbb{P}   (\cal{O}^{g_y^h}_{v,X^{(0,  g_y^h(t))}  })   e^{-\ov(g_y^h(t))v}   dv
ds    \hspace{-1pt}  +\hspace{-1pt}  \pp( \cal{O}^{g_y^h}_{t}  ;  \Delta_1^{g_y^h(t)}\hspace{-2pt} >\hspace{-1pt}t   )  
\\
%
  %
%
%
%
%
  %
  %
  %
%
\nonumber 
&\leq   \ov(g_{y}^{h}(t))     \Phi_y^h(t)      + \pp\l(\oo_{t,X^{(  0, g_{y}^{h}(t) ) } }^{g_{y}^{h}} \r)
- \ov(g_{y}^{h}(t))^2 \int_0^t\int_0^s \mathbb{P}   (\cal{O}_{v}^{g_{y}^{h}})   e^{-\ov(g_{y}^{h}(t))v}   dv ds  
\\     
  \label{triv}
&\leq   \ov(g_{y}^{h}(t))     \Phi_y^h(t)      + \pp\l(\oo_{t,X^{(  0, g_{y}^{h}(t) ) } }^{g_{y}^{h}} \r)    .
  \end{align}
%
Recall the notation $\l(\ref{D}\r)$, $\l(\ref{D()}\r)$, and  $\l(\ref{rhoyh}\r)$. By  $\l(\ref{triv}\r)$, partitioning on the value of $\D^{\big(\f{g_{y}^{h}(t)}{\log(t)} , g_y^h(t)\big)  } $, 
\b{align}
\nonumber
\rho_y^h(t) : \hspace{-6pt} &\overset{\phantom{\l(\ref{triv}\r)}}{=} \f{ \pp(  \cal{O}_t^{g_{y}^{h}}   )  }{  \Phi_y^h(t)   } -   \ov(g_{y}^{h}(t))
   \overset{\l(\ref{triv}\r)}\leq    \f{1}{\Phi_y^h(t)}   \pp\l(\oo_{t,X^{(  0, g_{y}^{h}(t) ) }  }^{g_{y}^{h}} \r)    
    \\
    \nonumber
      &\overset{\phantom{\l(\ref{triv}\r)}}=     \f{1}{\Phi_y^h(t)}  \Bigg[   \pp     \Big(      \oo_{t,X^{(0,g_{y}^{h}(t))} }^{ g_{y}^{h} } ; \D^{\f{g_{y}^{h}(t)}{\log(t)}   }      > t    \Big)  +    \pp   \Big(   \oo_{t,X^{(0,g_{y}^{h}(t))} }^{ g_{y}^{h} } ; \D^{\big(\f{g_{y}^{h}(t)}{\log(t)} , g_y^h(t)\big)  }      \leq        t   \Big) \Bigg]
\\
\label{ab}
     &\overset{\phantom{\l(\ref{triv}\r)}}{=} \hspace{-3pt}:    \f{1}{\Phi_y^h(t)}  \l[(a)+(b) \r].
     \end{align}

\n So   to prove  $\l(\ref{1}\r)$, we  need to prove, uniformly in $(h,y,t) \in \mbox{\hyperref[calS]{$\cal{S}$}}$,
\b{align}
\label{1'}
(a)+(b) &\lesssim   \f{\Phi_y^h(t)}{t\log(t)^{1+\eps}} \l( 1 + \f{1}{f(y)-h} \r) .
  \intertext{For $\l(\ref{2}\r)$, we  need        suitable $u$ so that uniformly in $(h,y,t) \in \mbox{\hyperref[calS]{$\cal{S}$}}$,}
\label{2'}
(a)+(b) &\leq  \Phi_y^h(t) u(t) \ov(g(t))  \l( 1 + \f{1}{f(y)-h} \r) .
\end{align} 

   \p{Proof for (\mbox{\hyperref[ab]{$a$}})
   }  
   Recall the notation $\l(\ref{truncoo}\r)$ and  $\l(\ref{ogyh}\r)$.  By $\l(\ref{subevent}\r)$,   
   \b{align}
   \label{fsoi}
   (\mbox{\hyperref[ab]{$a$}})     &\overset{\phantom{\ref{1fyh}}}= \pp     \Big(      \oo_{t,X^{(0,g_{y}^{h}(t))} }^{ g_{y}^{h} } ; \D^{\f{g_{y}^{h}(t)}{\log(t)}   }      > t    \Big)
  \overset{\l(\ref{subevent}\r)}\leq      \pp\Big(   X_t^{(0,\f{g_{y}^{h}(t)}{\log(t)})} >  g_{y}^{h}(t) \Big).
  \end{align}
 To bound $\l(\ref{fsoi}\r)$, we shall bound  $\pp(   X_t^{(0,\f{g_{y}^{h}(t)}{\log(t)})} > K g_{y}^{h}(t) ) $, with $K\in (0,1]$, giving a more general bound which  will also be used later in the proof of Lemma \mbox{\ref{lemma2}}. By Lemma \mbox{\ref{1.1}},  for all $(h,y,t) \in \mbox{\hyperref[calS]{$\cal{S}$}}$,  $g_y^h(t)/\log(t) \geq (1-A^{-1})g(t)/\log(t)> 1$,    so by  Lemma \mbox{\ref{1fyh}},  applying Lemma \hbox{\ref{keylemma}}   with $H(t) =  t^{-n}  $, for   $n>1$,   it follows that  uniformly in $(h,y,t) \in \mbox{\hyperref[calS]{$\cal{S}$}}$, 
     \begin{align}
       \nonumber
    \pp\Big(   X_t^{(0,\f{g_{y}^{h}(t)}{\log(t)})} > K  g_{y}^{h}(t) \Big)  
         &\overset{  \ref{1fyh}  }\leq \f{\Phi_y^h(t)}{f(y)-h}   \pp\Big(   X_t^{(0,\f{g_{y}^{h}(t)}{\log(t)})} > K  g_{y}^{h}(t) \Big)   \overset{  \ref{keylemma}  }\leq  
      \f{\Phi_y^h(t)}{f(y)-h}        \exp\l(   (*)    \r)   t^{-n}     ,
   \\
      \label{taref1}
   (*)  &\overset{\phantom{\ref{1fyh}}}\lesssim         t      \log(t^n)  \     t^{ n   \f{g_{y}^{h}(t)}{K g_{y}^{h}(t)\log(t)}    }    \  \ov\l( \f{g_{y}^{h}(t)}{\log(t)} \r)      \f{g_{y}^{h}(t)}{K g_{y}^{h}(t)\log(t)} 
   \\
   \nonumber
         &  \overset{\phantom{\ref{1fyh}}}= \f{n}{K} e^{\f{n}{K}} t \ov\l( \f{g_{y}^{h}(t)}{\log(t)} \r)   \lesssim 1  ,  
   \end{align}
    since   $\lim_{t\to\infty}   t \ov\l( g_{y}^{h}(t) / \log(t)\r)=0$, uniformly in $h,y$, by Lemma \mbox{\ref{1.1}} and $\l(\ref{condt}\r)$, and so  $\lim_{t\to\infty}(*)=0$. Thus        
  $(\mbox{\hyperref[ab]{$a$}})   \lesssim t^{-n}   \Phi_y^h(t)     / (f(y)-h  )\leq    t^{-1}\log(t)^{-1-\eps}    \Phi_y^h(t)     / (f(y)-h  ) $,  uniformly in $(h,y,t) \in \mbox{\hyperref[calS]{$\cal{S}$}}$,  as required for~$\l(\ref{1'}\r)$.

 To show $(\mbox{\hyperref[ab]{$a$}}) \leq \Phi_y^h(t) u(t) \ov(g(t)) /(f(y)-h) $,  by   Lemma \hbox{\ref{keylemma}} with $H(t) =  \ov(g(t))^2$,  applying Lemma \mbox{\ref{1.1}}, as $\ov$ is \mbox{\hyperref[rv]{regularly varying}} at $\infty$ with index $-\alpha\in(-1,0)$, uniformly in $(h,y,t) \in \mbox{\hyperref[calS]{$\cal{S}$}}$,   
\b{align}
\label{aaa}
 \pp\Big(   X_t^{(0,\f{g_{y}^{h}(t)}{\log(t)})} > K  g_{y}^{h}(t) \Big) 
   &\overset{\ref{keylemma}}\leq       
       \exp\l(  (*)    \r)  \ov(g(t))^2,
\\
\nonumber
(*) &\overset{\phantom{\ref{keylemma}}}\lesssim    t           \log\l(\ov(g(t))^{-2}\r)       \f{   }{ }    \ov(g(t))^{-\f{2}{K\log(t)}}    \ov\l( \f{g_{y}^{h}(t)}{\log(t)}  \r)         \f{ g_{y}^{h}(t)   }{   K g_{y}^{h}(t)       \log(t)     }     
\\
\nonumber
  &\overset{\ref{1.1}}\lesssim   t     \ov(g(t))^{-\f{2}{K\log(t)}}     \f{       \log\l(\ov(g(t))^{-2}\r) }{  \log(t)     }   \ov\l( \f{g(t)}{\log(t)}  \r)
  \\
\nonumber
  &\overset{\phantom{\ref{1.1}}}=  2 t \f{\log\l(\f{1}{\ov(g(t))}\r)}{\log(t)}     
  e^{\f{2}{K\log(t)} \log\l(\f{1}{\ov(g(t))}\r)  }  \ov\l( \f{g(t)}{\log(t)}  \r) .
 \intertext{Now,  using the inequality $x\leq e^x$, bounding $\ov(g(t)/\log(t))/\ov(g(t))\leq t$,  it follows that }
  \nonumber
  (*)  &\lesssim       t   \ov\l( \f{g(t)}{\log(t)}  \r) 
  e^{\f{2+K}{K\log(t)} \log\l(\f{1}{\ov(g(t))}\r)  }  
 \\
 \nonumber
  &=        t   \ov\l( \f{g(t)}{\log(t)}  \r) 
  e^{\f{2+K}{K\log(t)} \log\l(\f{t \ov\l(\f{g(t)}{\log(t)}\r)}{\ov(g(t))}\r)  + \f{2+K}{K\log(t)} \log\l(\f{1}{t\ov\l(\f{g(t)}{\log(t)}\r)}\r)      }  
\\
\label{finalbitpart}
  &\lesssim        t   \ov\l( \f{g(t)}{\log(t)}  \r) 
  e^{  \f{2+K}{k\log(t)} \log\l(\f{1}{t\ov\l(\f{g(t)}{\log(t)}\r)}\r)      }    =   
  e^{  \l[ \f{2+K}{K\log(t)}  -1 \r] \log\l(\f{1}{t\ov\l(\f{g(t)}{\log(t)}\r)}\r)      }   , 
 \end{align}
which is thus bounded uniformly in $(h,y,t) \in \mbox{\hyperref[calS]{$\cal{S}$}}$ because $\lim_{t\to\infty} t \ov(g(t)/\log(t))=0 $ by $\l(\ref{condt}\r)$. By $\l(\ref{aaa}\r)$ and Lemma \mbox{\ref{1fyh}}, it follows that   $(\mbox{\hyperref[ab]{$a$}})\leq \Phi_y^h(t)  u(t) \ov(g(t)) / (f(y)-h) $     uniformly in $(h,y,t) \in \mbox{\hyperref[calS]{$\cal{S}$}}$, for suitable $u$, as required for $\l(\ref{2'}\r)$.

   \p{Partitioning (\mbox{\hyperref[ab]{$b$}})
   }
 Now we   partition $(\mbox{\hyperref[ab]{$b$}})$. 
  Let $\Delta_m^{(a,b)}$ denote the time of our subordinator's $m$th jump of size larger in $(a,b)$, as in $\l(\ref{D()}\r)$.           With $\beta$ as in $\l(\ref{condt}\r)$,  for $m\geq1$ such that $m> \beta/(\beta-1)$ and  $m > 1/(\alpha(\beta-1)) $, for $c>0$ such that  $1    -(   m    -    1) c >0$,  for all $t$ large enough that $1/\log(t)   \leq c$,
  \b{align}
   \nonumber   
   (\mbox{\hyperref[ab]{$b$}}) &= \pp  \hspace{1pt}   \Big(  \oo_{t,X^{(0,g_{y}^{h}(t))} }^{ g_{y}^{h} } ; \D^{\big(\f{g_{y}^{h}(t)}{\log(t)} , g_y^h(t)\big)  }     \leq t \Big)           \phantom{\Bigg(}
 \\   
 \nonumber
   &= \pp \hspace{1pt}  \Big(  \oo_{t,X^{(0,g_{y}^{h}(t))} }^{ g_{y}^{h} } ; \D^{\big(\f{g_{y}^{h}(t)}{\log(t)} , cg_y^h(t)\big)  }    \leq t   ; \Delta_m^{\big(\f{g_{y}^{h}(t)}{\log(t)} , cg_y^h(t)\big)  }     > t     ;   \D^{    c g_{y}^{h}(t)   }     > t   \Big)         \phantom{\Bigg(}
\\
\nonumber
   &+
     \pp  \hspace{1pt}   \Big(  \oo_{t,X^{(0,g_{y}^{h}(t))} }^{ g_{y}^{h} } ;\Delta_m^{\big(\f{g_{y}^{h}(t)}{\log(t)} , cg_y^h(t)\big)  }     \leq  t        ;   \D^{    c g_{y}^{h}(t)   }     > t    \Big)            \phantom{\Bigg(}
   \\
   \nonumber
     &+
   \pp\l(  \oo_{t,X^{(0,g_{y}^{h}(t))} }^{ g_{y}^{h} } ;     \D^{ (   c g_{y}^{h}(t)  , g_y^h(t)) }     \leq t  \r)        \phantom{\Bigg(}
\\
\label{2ABC}
   &=: (2A) + (2B) + (2C). 
   \end{align}     

   \p{Proof for (\mbox{\hyperref[2ABC]{$2A$}})
   } Disintegrating on the value of $ \overline{\Delta}_1 := \D^{\big(\f{g_{y}^{h}(t)}{\log(t)} , cg_y^h(t)\big)  } $, the time of the first jump whose size lies between $g_{y}^{h}(t)/\log(t) $ and $c g_{y}^{h}(t)$ as defined in $\l(\ref{D()}\r)$, which is exponentially distributed,
    \b{align*}
   (\mbox{\hyperref[2ABC]{$2A$}}) &= \pp \hspace{1pt} \Big(  \oo_{t,X^{(0,g_{y}^{h}(t))} }^{ g_{y}^{h} } ; \D^{\big(\f{g_{y}^{h}(t)}{\log(t)} , cg_y^h(t)\big)  }     \leq t   ; \Delta_m^{\big(\f{g_{y}^{h}(t)}{\log(t)} ,c g_y^h(t)\big)  }     > t     ;   \D^{    c g_{y}^{h}(t)   }     > t   \Big) 
\\
   &= \int_0^t     \pp\Big(  \oo_{t,X^{(0,g_{y}^{h}(t))} }^{ g_{y}^{h} }  ; \Delta_m^{\big(\f{g_{y}^{h}(t)}{\log(t)} , cg_y^h(t)\big)  }    > t     ;   \D^{    c g_{y}^{h}(t)   }     > t  \ \Big| \    \overline{\Delta}_1      = s  \Big)    \pp\Big(   \overline{\Delta}_1      \in ds  \Big)
 \\
   &\leq  \ov\l(  \f{g_{y}^{h}(t)}{\log(t)}       \r)   \hspace{-1pt}   \int_0^t     \pp\Big(  \oo_{t,X^{(0,c g_{y}^{h}(t))} }^{ g_{y}^{h} } 
    ; \Delta_m^{\big(\f{g_{y}^{h}(t)}{\log(t)} , cg_y^h(t)\big)  }  > t   
       \ \Big| \     \overline{\Delta}_1     = s  \Big)   ds.
   \end{align*}        
   Now,  with $\overline{\Delta}_k$ denoting the time of the $k$th jump of size between $g_{y}^{h}(t)/\log(t) $ and $c g_{y}^{h}(t)$,  
\b{align*}
&     \pp\Big(  \oo_{t,X^{(0,c g_{y}^{h}(t))} }^{ g_{y}^{h} }  ; \Delta_m^{\big(\f{g_{y}^{h}(t)}{\log(t)} , cg_y^h(t)\big)  }    > t      \ \Big| \     \overline{\Delta}_1    = s  \Big)  
\\
    &=  \sum_{k=1}^{m-1}      \pp\Big(    \oo_{t,X^{(0,c g_{y}^{h}(t))} }^{ g_{y}^{h} }         \Big|      \overline{\Delta}_1         =      s  ; \overline{\Delta}_{k+1}          >       t     ;  \overline{\Delta}_k        \leq         t     \Big)  \pp\Big(    \overline{\Delta}_{k+1}     >      t     ;  \overline{\Delta}_k     \leq      t    \ \Big| \     \overline{\Delta}_1 = s   \Big)
\\
    &\leq    \sum_{k=1}^{m-1}  \pp\Big(  \oo_{t,X^{(0,c g_{y}^{h}(t))} }^{ g_{y}^{h} }        \ \Big| \       \overline{\Delta}_1         =      s  ; \overline{\Delta}_{k+1}          >       t     ;  \overline{\Delta}_k        \leq         t   \Big)  ,
  %
  %
  %
  \end{align*}
  and then it follows that
      \b{align}   
      \nonumber
   (\mbox{\hyperref[2ABC]{$2A$}}) &\overset{\phantom{\l(\ref{subevent}\r)}}\leq         \sum_{k=1}^{m-1}      \ov      \l(     \f{g_{y}^{h}(t)}{\log(t)}           \r)       \int_0^t          \pp\Big(  \oo_{t,X^{(0,c g_{y}^{h}(t))} }^{ g_{y}^{h} }        \ \Big| \       \overline{\Delta}_1         =      s  ; \overline{\Delta}_{k+1}          >       t     ;  \overline{\Delta}_k        \leq         t   \Big)  ds    .
  %
\intertext{Now, by $\l(\ref{subevent}\r)$, given that by time $t$ there are $k$  jumps of size in $[  g_{y}^{h}(t)/\log(t),c  g_{y}^{h}(t)]$,}     
\nonumber
 (\mbox{\hyperref[2ABC]{$2A$}})  &\overset{\l(\ref{subevent}\r)}\leq         \sum_{k=1}^{m-1}    \ov     \l(            \f{g_{y}^{h}(t)}{\log(t)}           \r)         \int_0^t           \pp        \Big(         X_t^{(0, c  g_{y}^{h}(t)   )}     >         g_{y}^{h}(t)   \Big|      \overline{\Delta}_1         =      s  ; \overline{\Delta}_{k+1}          >       t     ;  \overline{\Delta}_k        \leq         t      \Big)     ds  
 \\
 \nonumber
 &\overset{\phantom{\l(\ref{subevent}\r)}}\leq        \sum_{k=1}^{m-1}      \ov\l(  \f{g_{y}^{h}(t)}{\log(t)}       \r)      \int_0^t          \pp\Big(       X_{t}^{\big(0, \f{g_{y}^{h}(t)}{\log(t)}\big)} > (1-kc) g_{y}^{h}(t)    \Big)   ds
 \\
 \nonumber
   &\overset{\phantom{\l(\ref{subevent}\r)}}=      \sum_{k=1}^{m-1}      \ov\l(  \f{g_{y}^{h}(t)}{\log(t)}       \r)       t    \    \pp\Big(    X_t^{\big(0,  \f{g_{y}^{h}(t)}{\log(t)}   \big)}      >     (1-kc) g_{y}^{h}(t)   \Big)   
   \\
   \label{shortcutt}
   &\overset{\phantom{\l(\ref{subevent}\r)}}\leq   ( m-1)      \ov\l(  \f{g_{y}^{h}(t)}{\log(t)}       \r)       t    \    \pp\Big(    X_t^{\big(0,  \f{g_{y}^{h}(t)}{\log(t)}   \big)}      >     (1-(m-1)c) g_{y}^{h}(t)   \Big)   .
\intertext{Now,  $\lim_{t\to\infty}  t   \ov(  g_{y}^{h}(t)/\log(t)       )   \leq \lim_{t\to\infty}  t   \ov(  (1-1/A)g(t)/\log(t)       )    =0  $ by $\l(\ref{condt} \r)$, uniformly in $h,y$ by  Lemma \mbox{\ref{1.1}}.  Then by $\l(\ref{taref1}\r)$ and $\l(\ref{shortcutt}\r)$, it follows that   }   
   \label{2Alabel}
      (\mbox{\hyperref[2ABC]{$2A$}}) &\lesssim \pp\Big(    X^{\big(0,  \f{g_{y}^{h}(t)}{\log(t)}   \big)}      >     (1-(m-1)c) g_{y}^{h}(t)   \Big)  \overset{(\ref{taref1})}\lesssim     \f{  \Phi_y^h(t)     }{ f(y)-h }    t^{-n}   \lesssim  \f{  \Phi_y^h(t)     }{ f(y)-h }  \f{1}{  t \log(t)^{1+\eps}},        
\end{align}
uniformly in $(h,y,t) \in \mbox{\hyperref[calS]{$\cal{S}$}}$,  as required for $\l(\ref{1'}\r)$. 
  
   For the remaining bound on $(\mbox{\hyperref[2ABC]{$2A$}})$, applying $\l(\ref{aaa}\r)$ and $\l(\ref{shortcutt}\r)$, recalling  $(\ast)\lesssim 1$ by $\l(\ref{finalbitpart}\r)$,  it follows that $(\mbox{\hyperref[2ABC]{$2A$}}) \leq  \ov(g(t))^2    \leq  \Phi_y^h(t) u(t) \ov(g(t))  / (f(y)-h) $  for suitable $u$,   uniformly in $(h,y,t) \in \mbox{\hyperref[calS]{$\cal{S}$}}$,  as required for $\l(\ref{2'}\r)$.

 \p{Proof for (\mbox{\hyperref[2ABC]{$2B$}})}   Disintegrating on the value of $\Delta_1^{\big(\f{g_{y}^{h}(t)}{\log(t)} , cg_y^h(t)\big)  }    $, which   is exponentially distributed with parameter $\ov\l(g_{y}^{h}(t) /\log(t)\r) - \ov(cg_y^h(t)) \leq \ov\l(g_{y}^{h}(t) /\log(t)\r)$,  
 \b{align}
 \nonumber
 (\mbox{\hyperref[2ABC]{$2B$}}) &=    \pp\Big(  \oo_{t,X^{(0,g_{y}^{h}(t))} }^{ g_{y}^{h} } ;\Delta_m^{  \big(\f{g_{y}^{h}(t)}{\log(t)} , cg_y^h(t)\big)     }     \leq  t        ;   \D^{    c g_{y}^{h}(t)   }     > t    \Big)
\\
\nonumber
 &\leq      \ov\l(\f{g_{y}^{h}(t) }{\log(t)} \r)  \int_0^t          \pp\Big(  \oo_{t,X^{(0,g_{y}^{h}(t))} }^{ g_{y}^{h} } ;\Delta_m^{\big(\f{g_{y}^{h}(t)}{\log(t)} , cg_y^h(t)\big)    }     \leq  t        ; \D^{    c g_{y}^{h}(t)   }     > t  \  \big| \    \Delta_1^{\big(\f{g_{y}^{h}(t)}{\log(t)} , cg_y^h(t)\big)   }   =s  \Big)      ds   
\\
\nonumber
 &\leq    \ov\l(\f{g_{y}^{h}(t) }{\log(t)} \r) \int_0^t          \pp\Big(  \oo_{s }^{ g_{y}^{h} } ;\Delta_m^{\big(\f{g_{y}^{h}(t)}{\log(t)} , cg_y^h(t)\big)   }     \leq  t        \    \big|  \    \Delta_1^{\big(\f{g_{y}^{h}(t)}{\log(t)} , cg_y^h(t)\big)  }   =s  \Big)     ds   . 
 \intertext{By  $\l(\ref{s-}\r)$,    $\cal{O}_{s}^{g_y^h} = \cal{O}_{s-  }^{g_y^h}  $. Given $\Delta_1^{ \big(\f{g_{y}^{h}(t)}{\log(t)} , cg_y^h(t)\big)  } =s$,        $ \cal{O}_{s-}^{g_y^h} $  is independent of 
 $\Delta_m^{\big(\f{g_{y}^{h}(t)}{\log(t)} , cg_y^h(t)\big)  }   $, so}
%
\nonumber
  (\mbox{\hyperref[2ABC]{$2B$}})  &\leq     \ov\l(\f{g_{y}^{h}(t) }{\log(t)} \r) \int_0^t          \pp\Big(  \oo_{s  }^{ g_{y}^{h} }  \Big)   \pp\Big(     \Delta_m^{\big(\f{g_{y}^{h}(t)}{\log(t)} , cg_y^h(t)\big)  }     \leq  t             \   \big|  \   \Delta_1^{ \big(\f{g_{y}^{h}(t)}{\log(t)} , cg_y^h(t)\big)   }   =s   \Big)       ds
\\
\nonumber
  &\leq     \ov\l(\f{g_{y}^{h}(t) }{\log(t)} \r)     \pp\Big(     \Delta_{m-1}^{\f{g_{y}^{h}(t) }{\log(t)}   }     \leq  t    \Big)        \int_0^t          \pp\Big(  \oo_{s }^{ g_{y}^{h} }  \Big)   ds.  
\\ 
\intertext{Now,  since $ \Delta_1^{\f{g_{y}^{h}(t) }{\log(t)}  }   $  is exponentially distributed with parameter $\ov\l(g_{y}^{h}(t) /\log(t)\r)$, 
\[
   \pp\Big(     \Delta_{m-1}^{\f{g_{y}^{h}(t) }{\log(t)}   }     \leq  t    \Big)        \leq   \pp\Big(     \Delta_{1}^{  \f{g_{y}^{h}(t) }{\log(t)}   }     \leq  t    \Big)^{m-1}      =    \Big(     1 - e^{- t \ov\big(  \f{g_{y}^{h}(t) }{\log(t)}    \big)     } \Big)^{m-1}       \leq        t^{m-1} \ov\l(  \f{g_{y}^{h}(t) }{\log(t)}  \r)^{     m-1}  
   ,\]  
so recalling the notation in $\l(\ref{phiyh}\r)$, by Lemma \mbox{\ref{1.1}}, uniformly in $(h,y,t) \in \mbox{\hyperref[calS]{$\cal{S}$}}$,}         
  \\ \label{2Bstep}
%
 (\mbox{\hyperref[2ABC]{$2B$}})
 &\overset{\phantom{\ref{1.1}}}\leq        \ov\l(\f{g_{y}^{h}(t) }{\log(t)} \r)^{m}   t^{m-1}   \Phi_y^h(t)
%
%
 \overset{\ref{1.1}}\lesssim     \ov\l(\f{g(t) }{\log(t)} \r)^m t^{m-1}  \Phi_y^h(t).
  \end{align}
 Recall $\ov(x)=x^{-\alpha} L(x)$ for $L$  \mbox{\hyperref[sv]{slowly varying}}  at $\infty$, so by Potter's theorem \cite[Theorem 1.5.6]{bgt89}, for arbitrarily small $\delta>0$,  uniformly in $(h,y,t) \in \mbox{\hyperref[calS]{$\cal{S}$}}$,
 \b{equation} \label{log1}
   \ov\l(\f{g(t) }{\log(t)} \r)    = \log(t)^{\alpha} g(t)^{-\alpha} L(g(t))  \f{   L\l(\f{g(t) }{\log(t)} \r)  }{ L(g(t))  }  \lesssim   \log(t)^{\alpha + \delta}  \ov(g(t)).
 \end{equation}
 Similarly, defining   $g_\beta(t):=g(t)/\log(t)^\beta$, for $\beta>(1+\alpha)/(2\alpha+\alpha^2)$ as in $\l(\ref{condt}\r)$, uniformly in $(h,y,t) \in \mbox{\hyperref[calS]{$\cal{S}$}}$,
 \b{equation}  \label{log2}
    \ov\l(\f{g(t) }{\log(t)} \r)   =   \log(t)^{\alpha} g(t)^{-\alpha}   L\l(  g_\beta(t)   \r)  \f{   L\l(\f{g(t) }{\log(t)} \r)  }{ L(g_\beta(t))  }  \lesssim     \log(t)^{\alpha(1-\beta) + \delta(\beta-1)}\ov(g_\beta(t)).
 \end{equation}
Applying $\l(\ref{log1}\r)$ to $\ov(g(t))$ and $\l(\ref{log2}\r)$ to $\ov(g(t))^{m-1}$, then by $\l(\ref{condt}\r)$, uniformly in $(h,y,t) \in \mbox{\hyperref[calS]{$\cal{S}$}}$,
 \begin{align}
 \nonumber
  (\mbox{\hyperref[2ABC]{$2B$}})
 &\overset{\phantom{\ref{1.1}}}\lesssim  \log(t)^{ m\alpha + m \delta  - (m-1)\alpha \beta + \delta ( \beta-1) (m-1) }      \ov(g(t)) \ov(g_\beta(t))^{m-1}    t^{m-1}    \Phi_y^h(t)  
 \\
 \nonumber
 &\overset{(\ref{condt})}\leq    \log(t)^{ m\alpha + m \delta  - (m-1)\alpha \beta + \delta ( \beta-1) (m-1) }      \ov(g(t))  \Phi_y^h(t).
 \intertext{Now,   for $\beta$ as in $\l(\ref{condt}\r)$,  $m>\beta/(\beta-1)$, so  $m\alpha -   (m-1)\alpha\beta <0$,   choosing $\delta>0$ small enough, we conclude by  Lemma \mbox{\ref{1fyh}} that for suitable  $u$,  $ (\mbox{\hyperref[2ABC]{$2B$}})\leq  u(t) \ov(g(t)) \Phi_y^h(t) / (f(y)-h )  $, uniformly in $(h,y,t) \in \mbox{\hyperref[calS]{$\cal{S}$}}$,  as required for $\l(\ref{2'}\r)$.}
  \intertext{Now we prove  that $(\mbox{\hyperref[2ABC]{$2B$}})\lesssim        t^{-1}\log(t)^{-1-\eps}    \Phi_y^h(t)     / (f(y)-h  )   $. With $g_\beta(t)=g(t)/\log(t)^\beta$,  by $\l(\ref{2Bstep}\r)$ and $\l(\ref{log2}\r)$,   for arbitrarily small $\d>0$,  uniformly in $(h,y,t) \in \mbox{\hyperref[calS]{$\cal{S}$}}$,}
\nonumber
(\mbox{\hyperref[2ABC]{$2B$}}) &\leq        \ov\l(\f{g(t) }{\log(t)} \r)^m t^{m-1}  \Phi_y^h(t) 
\lesssim     \log(t)^{ m \alpha (1-\beta) + m \delta }   \ov(g_\beta(t))^m t^{m-1}  \Phi_y^h(t) .
 \end{align}
By $\l(\ref{condt}\r)$, $\ov(g_\beta(t))\lesssim t^{-1}$, and so
\(
(\mbox{\hyperref[2ABC]{$2B$}})\lesssim    t^{-1} \log(t)^{m\alpha(1-\beta) + m\delta  }  \Phi_y^h(t) 
\).
 Finally, applying Lemma \mbox{\ref{1fyh}}, our choice of $m$ ensures $m\alpha( 1- \beta) <-1$, so  choosing $\d$  small enough,  there exists $\eps>0$ such that   $(\mbox{\hyperref[2ABC]{$2B$}})\lesssim         t^{-1}\log(t)^{-1-\eps}    \Phi_y^h(t)     / (f(y)-h  ) $,  uniformly in $(h,y,t) \in \mbox{\hyperref[calS]{$\cal{S}$}}$, as~required~for~$\l(\ref{1'}\r)$.

\p{Partitioning $(\mbox{\hyperref[2ABC]{$2C$}}) 
$}
  Define $p^*(t) :=1- \log(t)^{-\gamma}   $ for  $\gamma:= (1-\alpha)/(2+\alpha)$, and   let $\Delta_2^{(a,b)}$ denote the time of our subordinator's second jump of size  in $(a,b)$. Then we partition:
   \b{align}      
   \nonumber
     (\mbox{\hyperref[2ABC]{$2C$}}) &=    \pp\l(  \oo_{t,X^{(0,g_{y}^{h}(t))} }^{ g_{y}^{h} } ;     \D^{  (  c g_{y}^{h}(t) , g_{y}^{h}(t) )  }     \leq t  \r)
   \\
   \nonumber
     &=  \pp\l(  \oo_{t,X^{(0,g_{y}^{h}(t))} }^{ g_{y}^{h} } ;     \D^{    (c g_{y}^{h}(t)  , g_{y}^{h}(t) )  }     \leq t       ;    \D^{   ( p^*(t) g_{y}^{h}(t) ,  g_{y}^{h}(t) )  }     \leq t             \r)
    \\
    \nonumber
     &+  \pp\l(  \oo_{t,X^{(0,g_{y}^{h}(t))} }^{ g_{y}^{h} }
      ;     \Delta_2^{   ( c g_{y}^{h}(t)  , g_{y}^{h}(t) )  }     \leq t       
     ;      \D^{   ( p^*(t) g_{y}^{h}(t)  , g_{y}^{h}(t) )  }     > t          \r)
   \\
   \nonumber
     &+  \pp\l(  \oo_{t,X^{(0,g_{y}^{h}(t))} }^{ g_{y}^{h} } 
     ;     \D^{   ( c g_{y}^{h}(t)  , g_{y}^{h}(t) )  }     \leq t       
      ;     \Delta_2^{  (  c g_{y}^{h}(t)  , g_{y}^{h}(t) )  }     > t
       ;         \D^{   ( p^*(t) g_{y}^{h}(t) , g_{y}^{h}(t) )   }     > t               \r)
  \\
  \label{2Cabc}
   &=:  (2Ca) + (2Cb) + (2Cc).
   \end{align}
   

%
%
%
%
%
%
%
%
%
  \p{Proof for (\mbox{\hyperref[2Cabc]{$2Ca$}})}  As $c\in(0,1)$ is fixed,   $c < 1-\log(t)^{-\gamma}=p^*(t)$ for all large enough $t$,  and so
    \b{align}
    \nonumber
    (\mbox{\hyperref[2Cabc]{$2Ca$}}) 
   &=     \pp\l(  \oo_{t,X^{(0,g_{y}^{h}(t))} }^{ g_{y}^{h} } ;         \D^{   ( p^*(t) g_{y}^{h}(t) ,  g_{y}^{h}(t) )  }     \leq t             \r).
 \intertext{Disintegrating on the value of $   \D^{    ( p^*(t) g_{y}^{h}(t) ,  g_{y}^{h}(t) )  }  $, then by $\l(\ref{subevent}\r)$ and the independence as in $\l(\ref{s-}\r)$,}
\nonumber
 (\mbox{\hyperref[2Cabc]{$2Ca$}}) &\overset{\phantom{\l(\ref{s-}\r)}}\leq   \l[   \ov\l( p^*(t) g_{y}^{h}(t)  \r)    - \ov\l( g_{y}^{h}(t)    \r)   \r]   \int_0^t    \pp\l(  \oo_{t,X^{(0,g_{y}^{h}(t))} }^{ g_{y}^{h} } |        \D^{   ( p^*(t) g_{y}^{h}(t) ,  g_{y}^{h}(t) )   }     =s             \r)      ds
\\
\nonumber
 &\overset{\l(\ref{s-}\r)}\leq    \l[   \ov\l( p^*(t) g_{y}^{h}(t)  \r)    - \ov\l( g_{y}^{h}(t)    \r)   \r]   \int_0^t    \pp\l(  \oo_{s,X^{(0,p^*(t) g_{y}^{h}(t))} }^{ g_{y}^{h} }       \r)      ds 
\\
\nonumber
   &\overset{\phantom{\l(\ref{s-}\r)}}\leq  \l[   \ov\l( p^*(t) g_{y}^{h}(t)  \r)    - \ov\l( g_{y}^{h}(t)    \r)   \r]  \Phi_y^h(t) =  \ov\l( p^*(t) g_{y}^{h}(t)  \r) \l(  1    -  \f{\ov\l( g_{y}^{h}(t)    \r) }{ \ov\l( p^*(t) g_{y}^{h}(t)  \r)}  \r)  \Phi_y^h(t) .
\intertext{By Lemma  \hbox{\ref{1.1}}, $\ov(p^*(t)  g_{y}^{h}(t))\lesssim \ov(p^*(t)  g(t))\lesssim \ov(g(t))$  uniformly in $(h,y,t) \in \mbox{\hyperref[calS]{$\cal{S}$}}$, so}
%
%
 \label{2Cabit}
  (\mbox{\hyperref[2Cabc]{$2Ca$}}) &\lesssim      \ov\l( g(t)  \r) \l(  1    -  \f{\ov\l( g_{y}^{h}(t)    \r) }{ \ov\l( p^*(t) g_{y}^{h}(t)  \r)}  \r)  \Phi_y^h(t) .     
%
%
%
%
 \intertext{As $\lim_{t\to\infty}p^*(t)   =  1$ and $\ov$ is  \mbox{\hyperref[crv]{CRV}} at $\infty$,   
 \(
    \lim_{t\to\infty}      \ov( p^*(t) g_{y}^{h}(t)  )    / \ov( g_{y}^{h}(t)  )        =  1.
 \)
 Now, $g_y^h(t) \gtrsim g(t)$ uniformly in $y,h$ by Lemma \mbox{\ref{1.1}},  so by $\l(\ref{2Cabit}\r)$ and Lemma \mbox{\ref{1fyh}},  $(\mbox{\hyperref[2Cabc]{$2Ca$}})\leq    \Phi_y^h(t)  u(t) \ov(g(t)) / (f(y)-h)  $ for suitable $u$,  uniformly in $(h,y,t) \in \mbox{\hyperref[calS]{$\cal{S}$}}$,  as~required~for~$\l(\ref{2'}\r)$. } 
\intertext{Now we prove     $(\mbox{\hyperref[2Cabc]{$2Ca$}})\lesssim    t^{-1} \log(t)^{-1-\eps}   \Phi_y^h(t) / (f(y)-h)  $. 
 As $\ov$ is \mbox{\hyperref[rv]{regularly varying}}   at $\infty$,  $\ov(x)=x^{-\alpha} L(x)$ for $L$  \mbox{\hyperref[sv]{slowly varying}}~at~$\infty$, where   $x^NL(x)$ is non-decreasing. Then by $\l(\ref{2Cabit}\r)$,     uniformly in $(h,y,t) \in \mbox{\hyperref[calS]{$\cal{S}$}}$, }
\nonumber
 (\mbox{\hyperref[2Cabc]{$2Ca$}}) &\lesssim       \ov\l( g(t)  \r) \l(  1    -  \f{ p^*(t)^\alpha L\l( g_{y}^{h}(t)    \r)  }{ L\l( p^*(t) g_{y}^{h}(t)  \r)}  \r)  \Phi_y^h(t) 
  =    \ov\l( g(t)  \r) \l(  1    -  \f{ p^*(t)^{\alpha+N} g_{y}^{h}(t)^N L\l( g_{y}^{h}(t)    \r)  }{ p^*(t)^{N} g_{y}^{h}(t)^N  L\l( p^*(t) g_{y}^{h}(t)  \r)}  \r)  \Phi_y^h(t)   
 \\
 \nonumber
 &\leq 
  \ov\l( g(t)  \r) \l( \hspace{17.5pt} 1    -   p^*(t)^{\alpha+N}    \hspace{17.5pt}  \r)  \Phi_y^h(t)   
     = 
    \ov\l( g(t)  \r) \l( \hspace{23.5pt} 1    -   (1-\log(t)^{-\gamma})^{\alpha+N}   \hspace{23.5pt}   \r)  \Phi_y^h(t)   
    \\
    \nonumber
    &\leq (\alpha+N) \ov\l( g(t)  \r)  \log(t)^{-\gamma}  \Phi_y^h(t)  . 
    %
    %
%
%
   %
  %
 %
 \intertext{For $\beta$ as in $\l(\ref{condt}\r)$, $g_\beta(t):=g(t)/\log(t)^\beta$,  by  Potter's theorem and     $\l(\ref{condt}\r)$,  for   arbitrarily small $\tau>0$,  \(
 \ov(g(t))\lesssim    \ov( g_{\beta}(t))  \log(t)^{-\alpha \beta  + \beta \tau}  \leq    t^{-1} \log(t)^{-\alpha \beta  + \beta \tau} 
\)    uniformly in $(h,y,t) \in \mbox{\hyperref[calS]{$\cal{S}$}}$,
so that}
\nonumber
 (\mbox{\hyperref[2Cabc]{$2Ca$}}) &\lesssim     t^{-1}     \log(t)^{-\gamma -\alpha \beta  + \beta \tau} \Phi_y^h(t).
\end{align}
Since $\gamma = (1-\alpha)/(2+\alpha) >1-\alpha\beta$, we may choose $\tau $ sufficiently  small    that $-\gamma -\alpha \beta  + \beta \tau<-1-\eps<-1 $.  Then it follows that
 $(\mbox{\hyperref[2Cabc]{$2Ca$}})\lesssim           t^{-1}\log(t)^{-1-\eps}    \Phi_y^h(t)    $,  uniformly in $(h,y,t) \in \mbox{\hyperref[calS]{$\cal{S}$}}$, as required for $\l(\ref{1'}\r)$.

  \p{Proof for (\mbox{\hyperref[2Cabc]{$2Cb$}})}   
  Disintegrating on  the value of $\hat{\Delta}_1^c:= \D^{ ( c g_{y}^{h}(t) , g_{y}^{h}(t) ) }$,  
  by~$\l(\ref{s-}\r)$,  
\begin{align}
  \nonumber
(\mbox{\hyperref[2Cabc]{$2Cb$}}) &\overset{\phantom{\l(\ref{s-}\r)}}=     \pp\l(  \oo_{t,X^{(0,g_{y}^{h}(t))} }^{ g_{y}^{h} } ;     \Delta_2^{  (   c g_{y}^{h}(t)  , g_{y}^{h}(t) ) }  \leq t   ;   \D^{  (  p^*(t) g_{y}^{h}(t) , g_{y}^{h}(t))  }     > t                       \r) 
\\
 \nonumber
  &\overset{\phantom{\l(\ref{s-}\r)}}\leq     \pp\l(  \oo_{t,X^{(0,g_{y}^{h}(t))} }^{ g_{y}^{h} } ;     \Delta_2^{  (   c g_{y}^{h}(t)  , g_{y}^{h}(t) )   }     \leq t              \r)
\\
  \nonumber
  &\overset{\phantom{\l(\ref{s-}\r)}}\leq   \ov\l(   c g_{y}^{h}(t)    \r)  \int_0^t      \pp     \l(      \oo_{   t,X^{(0,g_{y}^{h}(t))} }^{ g_{y}^{h} } ;     \Delta_2^{   ( c g_{y}^{h}(t) , g_{y}^{h}(t) )   }      \leq   t        \   \Big| \  \hat{\Delta}_1^c       =s   \r)    ds
\\
  \nonumber
  &\overset{\phantom{\l(\ref{s-}\r)}}\leq       \ov    \l(   c g_{y}^{h}(t)    \r)     \int_0^t      \pp   \l(      \oo_{   s,X^{(0,g_{y}^{h}(t))} }^{ g_{y}^{h} }    \ \Big| \   \hat{\Delta}_1^c    = s   \r)    \pp  \l(    \Delta_2^{    (   c g_{y}^{h}(t) , g_{y}^{h}(t) )  }      \leq   t   \    \Big|   \ \hat{\Delta}_1^c       =    s        \r)    ds
\\
  \nonumber
  &\overset{\l(\ref{s-}\r)}\leq    \ov\l(      c g_{y}^{h}(t)    \r)    \pp\l(     \D^{      c g_{y}^{h}(t)   }     \leq t     \r)  \int_0^t   \pp\l(    \oo_{s  }^{ g_{y}^{h} } \r)ds.
 %
 %
  \intertext{Recall for $L$    \mbox{\hyperref[sv]{slowly varying}}  at $\infty$, $\ov(x)   =x^{-\alpha} L(x)$.   Now,    $ \pp\big(      \D^{      c g_{y}^{h}(t)   }      \leq    t   \big)       =     1    -     e^{-t \ov(   c g_{y}^{h}(t)   )  }       \leq       t \ov(   c g_{y}^{h}(t))$,  so by   $\l(\ref{phiyh}\r)$ and  Lemma \mbox{\ref{1.1}}, uniformly in  $(h,y,t) \in \mbox{\hyperref[calS]{$\cal{S}$}}$,}
   %
   (\mbox{\hyperref[2Cabc]{$2Cb$}}) &\overset{\l(\ref{phiyh}\r)}\leq  \ov\l(      c g_{y}^{h}(t)    \r)^2   t       \Phi_y^h(t)
   \label{2cb}
     \overset{\ref{1.1}}\lesssim     \ov(g(t))^2 t   \Phi_y^h(t) 
  \overset{\phantom{\ref{1.1}}}=   g(t)^{-2\alpha}  L(g(t))^2   t \Phi_y^h(t)  
  \\
  \nonumber
   &\overset{\phantom{\ref{1.1}}}=    \log(t)^{ -2\alpha \beta }     \l(\f{g(t)}{\log(t)^{\beta}}\r)^{-2\alpha} L\l(\f{g(t)}{\log(t)^{\beta}}\r)^2     \f{ L(g(t))^2  }{ L\l(\f{g(t)}{\log(t)^{\beta}}\r)^2 }  t \Phi_y^h(t) 
 \\
  \nonumber
   &\overset{\phantom{\ref{1.1}}}=     \log(t)^{ -2\alpha \beta }    \ov\l(\f{g(t)}{\log(t)^{\beta}}   \r)^2         \f{ L(g(t))^2  }{ L\l(\f{g(t)}{\log(t)^{\beta}}\r)^2 }  t \Phi_y^h(t) .
   \intertext{By Potter's theorem  \cite[Theorem 1.5.6]{bgt89}, for arbitrarily small $\d  >  0$, uniformly in $t$,}
  \nonumber
     (\mbox{\hyperref[2Cabc]{$2Cb$}}) &\overset{\phantom{\ref{1.1}}}\lesssim \log(t)^{ -2\alpha \beta }    \ov\l(\f{g(t)}{\log(t)^{\beta}}   \r)^2       \log(t)^{2 \beta\d}         t \Phi_y^h(t) .
   \intertext{It follows by $\l(\ref{condt}\r)$ that       $ \ov\l(g(t)   /   \log(t)^{\beta}   \r)^2 \lesssim  t^{-2}$ uniformly in $t$, and hence}
  \nonumber
   (\mbox{\hyperref[2Cabc]{$2Cb$}}) &\overset{\phantom{\ref{1.1}}}\lesssim  \log(t)^{ -2\alpha \beta }         \log(t)^{2 \beta\d}         t^{-1} \Phi_y^h(t) 
   .
   \intertext{Now,    $2\alpha\beta > 1$ by $\l(\ref{condt}\r)$. Taking $\d$  small enough that $2 \alpha \beta - 2   \beta\d \geq 1+\eps>1$, it follows by Lemma \mbox{\ref{1fyh}} that $(\mbox{\hyperref[2Cabc]{$2Cb$}})    \lesssim       t^{-1}     \log(t)^{-1-\eps} \Phi_y^h(t) /(f(y)-h) $, uniformly in $(h,y,t) \in \mbox{\hyperref[calS]{$\cal{S}$}}$, as required~for~$\l(\ref{1'}\r)$.}
    \intertext{Now we show   $(\mbox{\hyperref[2Cabc]{$2Cb$}})\leq \Phi_y^h(t) u(t) \ov(g(t)) / (f(y)-h)  $. 
    By $\l(\ref{2cb}\r)$, since   $\lim_{t\to\infty} t\ov(g(t))=0$ by $\l(\ref{condt}\r)$, uniformly in $(h,y,t) \in \mbox{\hyperref[calS]{$\cal{S}$}}$,}
  \nonumber
   (\mbox{\hyperref[2Cabc]{$2Cb$}}) &\lesssim    \ov(g(t))^2 t   \Phi_y^h(t) = o(1) \times  \ov(g(t))   \Phi_y^h(t) ,
   \end{align}
  so by Lemma \mbox{\ref{1fyh}}, $(\mbox{\hyperref[2Cabc]{$2Cb$}})   \leq    \Phi_y^h(t) u(t) \ov(g(t)) /( f(y)-h) $   for suitable $u$, uniformly in $(h,y,t) \in \mbox{\hyperref[calS]{$\cal{S}$}}$, as required~for~$\l(\ref{2'}\r)$.

\p{Partitioning (\mbox{\hyperref[2Cabc]{$2Cc$}})} Disintegrating on the value of  $ \D^{c,p^*} := \D^{     ( c g_{y}^{h}(t)  ,     p^*(t) g_{y}^{h}(t) )     }$,   by $\l(\ref{s-}\r)$, $\l(\ref{subevent}\r)$, and Lemma \mbox{\ref{1.1}}, uniformly in $(h,y,t) \in \mbox{\hyperref[calS]{$\cal{S}$}}$,
\b{align}
\nonumber
(\mbox{\hyperref[2Cabc]{$2Cc$}}) &\overset{\phantom{\ref{1.1}}}=  \pp\l(  \oo_{t,X^{(0,g_{y}^{h}(t))} }^{ g_{y}^{h} } ;     \D^{    (   c g_{y}^{h}(t)  ,     g_{y}^{h}(t) )     }     \leq t        ;     \Delta_2^{   (    c g_{y}^{h}(t)  ,     g_{y}^{h}(t) )     }     > t ;         \D^{    (p^*(t) g_{y}^{h}(t)  ,     g_{y}^{h}(t) )     }     > t               \r)
\\
\nonumber
    &\overset{\phantom{\ref{1.1}}}\leq      \ov(       c g_{y}^{h}(  t  )    )      \int_0^t     \pp    \Big(    \oo_{    t,X^{(0,  g_{y}^{h}(t))} }^{ g_{y}^{h} }   ;      \Delta_2^{   (    c g_{y}^{h}(t)  ,     g_{y}^{h}(t) )      }       >   t      ;       \D^{   ( p^*(t) g_{y}^{h}(t)  ,     g_{y}^{h}(t) )      }       >     t         \big|       \D^{     c,p^*     }     =  s     \Big)       ds
\\
\nonumber
&\overset{\l(\ref{subevent}\r) }\leq   \ov(    c g_{y}^{h}(  t  )   )    \int_0^t      \pp    \Big(    \oo_{    s,X^{(0,   c g_{y}^{h}(t))} }^{ g_{y}^{h} }    ;      \hat{X}_{t-s}^{(0,   c g_{y}^{h}(t))}       +     X_{s-}^{(0,   c g_{y}^{h}(t))}     >      (  1 -  p^*(  t   )  )g_{y}^{h}(  t )      \big|       \D^{     c,p^*     }     =  s         \Big)  ds
\\
\nonumber
&\overset{\l(\ref{s-}\r) }=   \ov(    c g_{y}^{h}(  t  )   )    \int_0^t      \pp    \Big(    \oo_{    s,X^{(0,   c g_{y}^{h}(t))} }^{ g_{y}^{h} }    ;      \hat{X}_{t-s}^{(0,   c g_{y}^{h}(t))}       +     X_{s-}^{(0,   c g_{y}^{h}(t))}     >      (  1 -  p^*(  t   )  )g_{y}^{h}(  t )                   \Big)  ds
\\
\nonumber  
&\overset{\ref{1.1}}\lesssim   \ov(    g(  t  ) )   \int_0^t    \pp    \Big(    \oo_{    s,X^{(0,   c g_{y}^{h}(t))} }^{ g_{y}^{h} }    ;    \hat{X}_{t-s}^{(0,   c g_{y}^{h}(t))}       +     X_{s-}^{(0,   c g_{y}^{h}(t))}      >      (  1  -  p^*(  t   )  )g_{y}^{h}(  t  )           \Big)ds,
\intertext{where   $\hat{X}$ is an independent copy of $X$, and we use that the jump at time $s$ has size at most $p^*(t) g_y^h(t)$. Recall that $1-p^*(t) = \log(t)^{-\gamma}$. Then partitioning according to $ \big\{X_{s-}^{(0,   c g_{y}^{h}(t))}    >   g_{y}^{h}(t)/(2\log(t)^\gamma)    \big\}$,}
\nonumber
(\mbox{\hyperref[2Cabc]{$2Cc$}})&\lesssim       \ov(g(t)) \int_0^t  \pp\l(  \oo_{s,X^{(0,   c g_{y}^{h}(t))} }^{ g_{y}^{h} }  \ ; \   \hat{X}_{t-s}^{(0,   c g_{y}^{h}(t))}       >   \f{ g_{y}^{h}(t)          }{   2 \log(t)^\gamma  }           \r)   ds
\\
\nonumber
&+   \ov(g(t)) \int_0^t  \pp\l(  \oo_{s,X^{(0,   c g_{y}^{h}(t))} }^{ g_{y}^{h} }  \ ; \      X_{s-}^{(0,   c g_{y}^{h}(t))}    >   \f{ g_{y}^{h}(t)          }{   2 \log(t)^\gamma  }   \r)   ds
\\
\label{SS*}
&=:       (S) + (S^*).
\end{align}
Next we will bound $(S)$, then later we will split up $(S^*)$ into more pieces.

  \p{Proof for (\mbox{\hyperref[SS*]{$S$}})} 
As $\hat{X}$ is an  independent copy of $X$, we can write
  \b{align}
  (\mbox{\hyperref[SS*]{$S$}}) 
%
%
%
\nonumber
  &=      \ov(g(t)) \int_0^t  \pp\l(  \oo_{s,X^{(0,   c g_{y}^{h}(t))} }^{ g_{y}^{h} }  \r)  \pp\l(   X_{t-s}^{(0,   c g_{y}^{h}(t))}       >   \f{ g_{y}^{h}(t)          }{   2 \log(t)^\gamma  }           \r)   ds
\\
\label{S0}
  &\leq   \ov(g(t))          \pp\l( X_{t}^{(0,   c g_{y}^{h}(t))}       >   \f{ g_{y}^{h}(t)     }{   2 \log(t)^\gamma  }           \r)     \Phi_y^h(t).
 \intertext{Since $\ov$ is \mbox{\hyperref[rv]{regularly varying}} at $\infty$,   applying Potter's theorem  \cite[Theorem 1.5.6]{bgt89}, with $\beta$ as in $\l(\ref{condt}\r)$, to $\ov(g(t))/\ov(g(t)/\log(t)^\beta)$,       for arbitrarily small $\tau>0$, it follows by $\l(\ref{condt}\r)$ that uniformly in $(h,y,t) \in \mbox{\hyperref[calS]{$\cal{S}$}}$,}
\label{S1}
 (\mbox{\hyperref[SS*]{$S$}}) &\lesssim \ov\l(\f{g(t)}{\log(t)^\beta}\r)  \log(t)^{-\alpha \beta + \tau \beta     }    \pp\l( X_{t}^{(0,   c g_{y}^{h}(t))}       >   \f{ g_{y}^{h}(t)     }{   2 \log(t)^\gamma  }           \r)     \Phi_y^h(t) 
 \\
 \label{annoyingbit}
    &\lesssim        t^{-1}   \log(t)^{-\alpha \beta + \tau \beta     }    \pp\l( X_{t}^{(0,   c g_{y}^{h}(t))}       >   \f{ g_{y}^{h}(t)     }{   2 \log(t)^\gamma  }           \r)     \Phi_y^h(t).
\end{align}
Now, we will show that there exists $\eps>0$ such that uniformly in $(h,y,t) \in \mbox{\hyperref[calS]{$\cal{S}$}}$, 
\b{equation}
\label{gamma'bound}
 \pp\l( X_{t}^{(0,   c g_{y}^{h}(t))}       >   \f{ g_{y}^{h}(t)     }{   2 \log(t)^\gamma  }           \r)  
  \lesssim  
    \log(t)^{-1-\eps + \alpha\beta - \tau\beta}.
\e{equation}
Then it follows from $\l(\ref{annoyingbit}\r)$, $\l(\ref{gamma'bound}\r)$, and Lemma \mbox{\ref{1fyh}} that \[
  (\mbox{\hyperref[SS*]{$S$}})\lesssim     \f{\Phi_y^h(t)}{t   \log(t)^{1+ \alpha\beta } }    \leq   \f{\Phi_y^h(t)}{t   \log(t)^{1+ \eps} }  \leq   \f{\Phi_y^h(t)}{t   \log(t)^{1+ \eps} }    \l( 1 + \f{1}{f(y)-h} \r) ,
\]
 uniformly in $(h,y,t) \in \mbox{\hyperref[calS]{$\cal{S}$}}$, as required for $\l(\ref{1'}\r)$. 
 Moreover, it follows by $\l(\ref{S0}\r)$ and  $\l(\ref{gamma'bound}\r)$ that for suitable $u$, 
 $(\mbox{\hyperref[SS*]{$S$}})\leq  \Phi_y^h(t) u(t) \ov(g(t)) / (f(y)-h) $, uniformly  in $(h,y,t) \in \mbox{\hyperref[calS]{$\cal{S}$}}$, as required for $\l(\ref{2'}\r)$.

\n Now, to prove $\l(\ref{gamma'bound}\r)$, we set  $M:= \gamma+ 2/(2+\alpha)= (3-\alpha)/(2+\alpha) $, and partition as follows 
%
%
%
\[
 \pp\big( X_{t}^{(0,   c g_{y}^{h}(t))}       >    g_{y}^{h}(t)    /   2 \log(t)^\gamma              \big) 
\]
\b{equation}
\label{twoprinces}
 = \! \pp\Big( X_{t}^{(0,   c g_{y}^{h}(t))}   \!    >\!   \f{ g_{y}^{h}(t)     }{   2 \log(t)^\gamma  }      ; \Delta_1^{\f{g_y^h(t)}{\log(t)^{M}}} \!\leq\! t      \Big)   
 \!  +\!
    \pp\Big( X_{t}^{(0,   c g_{y}^{h}(t))}     \!  > \!  \f{ g_{y}^{h}(t)     }{   2 \log(t)^\gamma  }    ;   \Delta_1^{\f{g_y^h(t)}{\log(t)^{M}}} \!>\! t     \Big)  =: (Q)\! +\! (Q') .
\e{equation}
Then as $\Delta_1^{\f{g_y^h(t)}{\log(t)^{M}}}$ is exponentially distributed with rate $\ov(g_y^h(t)/\log(t)^{M})$, we can bound
\[
(Q)  \leq  \pp\Big( \Delta_1^{\f{g_y^h(t)}{\log(t)^{M}}} \leq t      \Big) \leq 1- e^{-t\ov(g_y^h(t)/\log(t)^{M})} \leq t \ov\l(\f{g_y^h(t)}{\log(t)^{M}}\r).
\]
By Lemma \mbox{\ref{1.1}} and $\l(\ref{condt}\r)$, applying Potter's Theorem to $\ov(g(t)/\log(t)^M)/\ov(g(t)/\log(t)^\beta)$, it follows that for arbitrarily small $\kappa>0$, uniformly in $(h,y,t) \in \mbox{\hyperref[calS]{$\cal{S}$}}$, 
\[
(Q)
 \overset{\ref{1.1}}\lesssim 
   t \ov\l(\f{g(t)}{\log(t)^{M}}\r)   
 \lesssim  
  t \ov\l(\f{g(t)}{\log(t)^\beta}\r)   \log(t)^{-\alpha (\beta-M) + \kappa( \beta-M)}  
   \overset{(\ref{condt})}\lesssim 
     \log(t)^{-\alpha (\beta-M) + \kappa( \beta-M)}  ,
\]
and then in order for $\l(\ref{gamma'bound}\r)$ to hold,  we  need $-\alpha (\beta-M)+ \kappa( \beta-M)\leq-1-\eps+\alpha\beta-\tau\beta$, so taking $\kappa,\tau,\eps$ small enough, we need
$
-\alpha (\beta-M) < -1 + \alpha\beta,
$ that is, $(3-\alpha) / (2+\alpha) =M < 2\beta -1/\alpha$. 
This is indeed true since $\beta> (1+2\alpha)/(2\alpha+\alpha^2)$, from which it follows that
\[
2\beta -\f{1}{\alpha}  > 2\f{1+2\alpha}{2\alpha+\alpha^2}  -\f{1}{\alpha} = \f{2+4\alpha -2 -\alpha}{2\alpha + \alpha^2} = \f{3}{2+\alpha} > \f{3-\alpha}{2+\alpha} = M,
\]
and the desired bound for $(Q)$ holds.
To bound $(Q')$, we provide only a proof for $\alpha\beta\leq1$, as a similar, simpler argument works when $\alpha\beta>1$. By Lemma \mbox{\ref{keylemma}} with $H(t)=  \log(t)^{-1-\eps + \alpha\beta - \tau\beta}$, 
\b{equation}
\label{(*)bitproof}
(Q')
\leq 
  \pp\Big( X_{t}^{(0,   \f{g_y^h(t)}{\log(t)^{M}} )}       >   \f{ g_{y}^{h}(t)     }{   2 \log(t)^\gamma  }       \Big) 
  \overset{\ref{keylemma}}\leq \exp((*)) \log(t)^{-1-\eps + \alpha\beta - \tau\beta},
\e{equation}
\[
(*) \lesssim  t \log(\log(t)^{1+\eps - \alpha\beta + \tau\beta}) \log(t)^{[1+\eps - \alpha\beta + \tau\beta] 2 \log(t)^{\gamma-M}}   \ov\l(\f{g_y^h(t)}{\log(t)^M}   \r)   2 \log(t)^{\gamma-M}.
\]
Now, for $\eta >0$  small enough that $\eta<M-\gamma$, observe that since $M>\gamma$, uniformly in $t>t_0(y)$,
\[
(*) \lesssim  t \log(t)^{\f{\eta}{2}}  \log(t)^{\f{\eta}{2}}  \ov\l(\f{g_y^h(t)}{\log(t)^M}   \r)   \log(t)^{\gamma-M} \leq t \ov\l(\f{g_y^h(t)}{\log(t)^M}   \r) .
\]
Now, recall  $M= (3-\alpha)/(2+\alpha) $, $\beta>(1+2\alpha)/(2\alpha+\alpha^2)$, and $\alpha<1$,  
from which one can verify that $M<\beta$. Then it follows by Lemma \mbox{\ref{1.1}} and $\l(\ref{condt}\r)$ that uniformly in $(h,y,t) \in \mbox{\hyperref[calS]{$\cal{S}$}}$, 
\[
(*) \overset{\ref{1.1}}\lesssim t \ov\l(\f{g(t)}{\log(t)^M}   \r) \leq t \ov\l(\f{g(t)}{\log(t)^\beta}   \r) \overset{(\ref{condt})}= o(1),
\]
so the desired bounds for $(Q)$ and $(Q')$ are proven, and the proof of $\l(\ref{gamma'bound}\r)$ is complete.

\p{Partitioning  (\mbox{\hyperref[SS*]{$S^*$}})} Now we   partition $(\mbox{\hyperref[SS*]{$S^*$}})$. For $\beta$ as in $\l(\ref{condt}\r)$, $\gamma= (1-\alpha)/(2+\alpha)$, and $\d:=1+\gamma$, write $g_\d(t):= g(t)/\log(t)^\d$. Recall the notation $\l(\ref{D}\r)$. Then 
\b{align}
\nonumber
(\mbox{\hyperref[SS*]{$S^*$}}) &=     \ov(g(t)) \int_0^t  \pp\l(  \oo_{s,X^{(0,   c g_{y}^{h}(t))} }^{ g_{y}^{h} }  \ ; \      X_{s-}^{(0,   c g_{y}^{h}(t))}    >   \f{ g_{y}^{h}(t)          }{   2 \log(t)^\gamma  }   \r)   ds
\\
\nonumber
&=   \ov(g(t)) \int_0^t  \pp\l(  \oo_{s,X^{(0,   c g_{y}^{h}(t))} }^{ g_{y}^{h} }  \ ; \      X_{s-}^{(0,   c g_{y}^{h}(t))}    >   \f{ g_{y}^{h}(t)          }{   2 \log(t)^\gamma  }     \ ; \    \D^{g_\d(t)} \leq s        \r)   ds
\\
\nonumber
&+  \ov(g(t)) \int_0^t  \pp\l(  \oo_{s,X^{(0,   c g_{y}^{h}(t))} }^{ g_{y}^{h} }  \ ; \      X_{s-}^{(0,   c g_{y}^{h}(t))}    >   \f{ g_{y}^{h}(t)          }{   2 \log(t)^\gamma  }     \ ; \      \D^{g_\d(t)} >s     \r)   ds
\\
\label{S1*S2*}
&=: (S_1^*) + (S_2^*).
\end{align}

      \p{Proof for (\mbox{\hyperref[S1*S2*]{$S_1^*$}})}
       Disintegrating on the value of $ \D^{g_{\d}(t)} $, by $\l(\ref{s-}\r)$~and~$\l(\ref{phiyh}\r)$,
      \begin{align}
      \nonumber
      ( \mbox{\hyperref[S1*S2*]{$S_1^*$}})  &\overset{\phantom{\l(\ref{phiyh}\r)}}\leq      \ov( g(   t   )   )   \ov (    g_{\d}(  t   )   )     \int_0^t    \int_0^s     \pp \Big(   \oo_{   s,X^{(0,   c g_{y}^{h}(t))} }^{ g_{y}^{h} }    ;        X_{s-}^{(0,   c g_{y}^{h}(t))}     >     \f{ g_{y}^{h}(t)          }{   2 \log(t)^\gamma  }     \big|   \D^{g_{\d}(t)}  =  v         \Big)    dv  ds
  \\
  \nonumber
        &\overset{\phantom{\l(\ref{phiyh}\r)}}\leq    \ov( g(   t   )   )   \ov (    g_{\d}(   t   )   )     \int_0^t   \int_0^t     \pp \Big(   \oo_{   v,X^{(0,   g_{\d}(   t   )  )} }^{ g_{y}^{h} }    ;        X_{t}^{(0,   c g_{y}^{h}(t))}     >     \f{ g_{y}^{h}(t)          }{   2 \log(t)^\gamma  }     \big|   \D^{g_{\d}(t)}  =  v         \Big)    dv  ds
    \\
      &\overset{ \l(\ref{s-}\r)}\leq     t     \ov( g(   t   )   )   \ov (    g_{\d}(  t  )   )    \int_0^t    \pp \Big(   \oo_{   v,X^{(0,   g_{\d}(   t   ) )} }^{ g_{y}^{h} }  \Big)    dv   
    \label{s1*}
       \overset{\l(\ref{phiyh}\r)}\leq      t \ov(g(t))  \ov(g_{\d}(t))         \Phi_y^h(t).
      \intertext{As $\ov$ is \mbox{\hyperref[rv]{regularly varying}} with index $-\alpha$,  for $g_\beta(t)  :=   g(t)/\log(t)^\beta$ and   $\beta$ as in $\l(\ref{condt}\r)$, applying Potter's theorem  \cite[Theorem 1.5.6]{bgt89} to $\ov(g(t))/\ov(g_\beta(t))$ and $\ov(g_\delta(t))/\ov(g_\beta(t))$,  for $\tau>0$, uniformly in $t$,}
      \nonumber
     ( \mbox{\hyperref[S1*S2*]{$S_1^*$}}) &\lesssim t   \ov(g_{\beta}(t))^2  \log(t)^{   -\alpha \beta    -\alpha(\beta-\delta)     +\beta\tau   + (\beta-\d)\tau }    \Phi_y^h(t),
          \intertext{Now, by $\l(\ref{condt}\r)$, $\lim_{t\to\infty} t \ov(g_\beta(t))   =  0$,~so~uniformly in $t$,}
    \nonumber
    ( \mbox{\hyperref[S1*S2*]{$S_1^*$}})   &\leq   t^{-1}   \log(t)^{   -\alpha \beta    -\alpha(\beta-\delta)     +\beta\tau   + (\beta-\d)\tau }   \Phi_y^h(t).
\end{align}

\n  Now, one can verify that  $-\alpha\beta - \alpha(\beta- \d )< -1$, using  that $\d=1+\gamma = 1+ (1-\alpha)/(2+\alpha)$ and $\beta>(1+2\alpha)/(2\alpha+\alpha^2)$. So taking   $\tau$ small enough, we  conclude by Lemma \mbox{\ref{1fyh}} that     $( \mbox{\hyperref[S1*S2*]{$S_1^*$}}) \lesssim   t^{-1}\log(t)^{-1-\eps}  \Phi_y^h(t) / (f(y)-h) $, uniformly in $(h,y,t) \in \mbox{\hyperref[calS]{$\cal{S}$}}$, as required for $\l(\ref{1'}\r)$.

  To show   $( \mbox{\hyperref[S1*S2*]{$S_1^*$}}) \leq  \Phi_y^h(t) u(t) \ov(g(t)) /(f(y)-h) $,  one can verify that for each $\alpha\in(0,1)$ and for $\beta$ as in $\l(\ref{condt}\r)$,      $\d =  1+ (1-\alpha)/(2+\alpha) <    (1+2\alpha)/(2\alpha+\alpha^2)   < \beta$. Thus   $\lim_{t\to\infty} t\ov(  g_\d(t))=0$ by $\l(\ref{condt}\r)$. Then by $\l(\ref{s1*}\r)$ and Lemma \mbox{\ref{1fyh}}, for suitable $u$, uniformly in $y>g(h),h>0$ as $t  \to  \infty$, 
\(
      ( \mbox{\hyperref[S1*S2*]{$S_1^*$}})        \lesssim   t \ov(g(t))  \ov(g_{\d}(t))         \Phi_y^h(t)   \leq    u(t) \ov(g(t)) \Phi_y^h(t) / (f(y)-h) $,  as required~for~$\l(\ref{2'}\r)$.

  \p{Proof for (\mbox{\hyperref[S1*S2*]{$S_2^*$}})} 
  Note for $g_\d(t):=g(t)/\log(t)^\d$,  by Lemma \mbox{\ref{1.1}},  for all $h>0$, $y>g(h)$, and for all large enough $t$, $g_\d(t) \leq c g_y^h(t)$, so   
  \b{align}
  \nonumber
  (\mbox{\hyperref[S1*S2*]{$S_2^*$}})  &=      \ov(g(t)) \int_0^t  \pp\l(  \oo_{s,X^{(0,   c g_{y}^{h}(t))} }^{ g_{y}^{h} }  \ ; \      X_{s-}^{(0,   c g_{y}^{h}(t))}    >   \f{ g_{y}^{h}(t)          }{   2 \log(t)^\gamma  }     \ ; \      \D^{g_{\d}(t)} >s     \r)   ds
\\
\nonumber
 &=      \ov(g(t)) \int_0^t  \pp\l(  \oo_{s,X^{(0,  g_{\d}(t))} }^{ g_{y}^{h} }  \ ; \      X_{s-}^{(0,  g_{\d}(t))}    >   \f{ g_{y}^{h}(t)          }{   2 \log(t)^\gamma  }     \ ; \      \D^{g_{\d}(t)} >s     \r)   ds
 \\
\label{Horatio'scomment}
 &\leq     \ov(g(t)) \int_0^t  \pp\l(  X_{s-}^{(0,  g_{\d}(t))}    >   \f{ g_{y}^{h}(t)          }{   2 \log(t)^\gamma  }        \r)   ds
  \leq  t \ov(g(t))   \pp\l(     X_{t}^{(0,  g_{\d}(t))}    >   \f{ g_{y}^{h}(t)          }{   2 \log(t)^\gamma  }          \r)   .
%
%
  \intertext{For  $g_\beta(t):=g(t)/\log(t)^\beta$, with $\beta$ as in $\l(\ref{condt}\r)$,   applying Potter's theorem \cite[Theorem 1.5.6]{bgt89} to $\ov(g(t))$,  for arbitrarily small $\tau>0$,  by Lemma \mbox{\ref{1.1}},}
\nonumber
      (\mbox{\hyperref[S1*S2*]{$S_2^*$}})  &\overset{\phantom{\ref{1.1}}}\lesssim     t \ov(g_\beta(t))   \log(t)^{  -\alpha\beta + \tau \beta }        \pp\l(     X_{t}^{(0,  g_{\d}(t))}    >   \f{ g_{y}^{h}(t)          }{   2 \log(t)^\gamma  }          \r)    
     \\
      \nonumber  
      &\overset{\ref{1.1}}\leq       t \ov(g_\beta(t))   \log(t)^{  -\alpha\beta + \tau \beta }        \pp\l(     X_{t}^{(0,  g_{\d}(t))}    >   \f{(1-A^{-1}) g (t)          }{   2 \log(t)^\gamma  }          \r)    .
  \intertext{Applying Lemma \mbox{\ref{keylemma}} with   $H(t) = 1/(t  \log(t)^{1 + \eps   -\alpha\beta    })$, $\eps>\tau\beta $,   then applying $\l(\ref{condt}\r)$ and Lemma \mbox{\ref{1fyh}}, uniformly in $(h,y,t) \in \mbox{\hyperref[calS]{$\cal{S}$}}$ by Lemma \mbox{\ref{1.1}},   }
  \nonumber
   (\mbox{\hyperref[S1*S2*]{$S_2^*$}})  &\overset{\hspace{1.5pt}\ref{keylemma}  \hspace{1.5pt}  }\lesssim           t \ov(g_\beta(t))   \log(t)^{  -\alpha\beta + \tau \beta }  \exp ((*)) \f{1}{t  \log(t)^{1 + \eps   -\alpha\beta    }}
   \overset{ \hspace{1.5pt} \l(\ref{condt}\r)  \hspace{1.5pt}   }=  o(1) \times  \exp ((*)) \f{1}{t  \log(t)^{1 + \eps   -\tau\beta    }}
   \\
   \label{S2*bit}
   &\overset{\ref{1fyh}}\leq      o(1) \times       \exp ((*)) \f{1}{t  \log(t)^{1 + \eps   -\tau\beta    }}    \f{    \Phi_y^h(t) }{    f(y)-h},
   \\
   \nonumber
    (*) &\overset{\phantom{\ref{1fyh}}}\lesssim   t  \log\l(t    \log(t)^{1+\eps -\alpha\beta}    \r)    \l(t    \log(t)^{1+\eps -\alpha\beta}    \r)^{   \f{2  \log(t)^{\gamma-\d}}{1-A^{-1}}  }   \ov(g_\d(t))   \log(t)^{\gamma - \d   }.
  \intertext{Now, if $\lim_{t\to\infty} (*)     =0$, then $(\mbox{\hyperref[S1*S2*]{$S_2^*$}})     \lesssim      t^{-1}        \log(t)^{-1-\eps}           \Phi_y^h(t)  /       (   f( y)         -      h) $. Indeed,    
    }
 \nonumber
   (*) &\lesssim   t  \log\l(t    \log(t)^{1+\eps -\alpha\beta}    \r)    \l(t    \log(t)^{1+\eps -\alpha\beta}    \r)^{  \f{2  \log(t)^{\gamma-\d} }{ 1- A^{-1} }   }   \ov(g_\d(t))   \log(t)^{\gamma - \d   }
\\
\nonumber
  &=    t  \ov(g_\d(t))      t^{ \f{ 2\log(t)^{\gamma-\d}  }{ 1- A^{-1} }     }       \log\big(t    \log(t)^{1+\eps -\alpha\beta}    \big)     \log(t)^{  \l( 1+\eps -\alpha\beta   \r) \f{2\log(t)^{\gamma-\d}}{ 1- A^{-1} }    +\gamma-\d   }.
   \intertext{Now, $\d  >  \gamma$, so  $\lim_{t\to\infty}\log(t)^{\gamma-\d}    =  0$, and for arbitrarily small $\kappa  >  0$,     uniformly in $t$,
   }
   \nonumber
    (*) &\lesssim  t  \ov(g_\d(t))     t^{ \f{ 2\log(t)^{\gamma-\d}  }{ 1- A^{-1} }     }         \log\big(t    \log(t)^{1+\eps -\alpha\beta}    \big)     \log(t)^{  \l( 1+\eps -\alpha\beta   \r) \kappa   +\gamma-\d   }.
\intertext{Note that $\log(t\log(t)^{1+\eps - \alpha\beta})\lesssim\log(t)$, uniformly in $t$. Then applying Potter's theorem    \cite[Theorem 1.5.6]{bgt89}    
   to  $  \ov(g_\d(t)) $, for $\beta$ as in $\l(\ref{condt}\r)$ and arbitrarily small $c  >  0$, uniformly in $t$,  }
\nonumber
 (*) &\lesssim     t  \ov(g_\beta(t))         t^{ \f{ 2\log(t)^{\gamma-\d}  }{ 1- A^{-1} }     }         \log(t)^{ 1 +  \l( 1+\eps -\alpha\beta   \r)\kappa -(\beta-\d)\alpha   + (\beta-\delta)c  +\gamma-\d    }.
  \intertext{Recalling     $\gamma   -  \d     =     -   1$,  \(
  t^{ \f{ 2\log(t)^{\gamma-\d}  }{ 1- A^{-1} }     }      =      e^{  \f{  2\log(t)^{1+\gamma-\d}  }{ 1- A^{-1}  }   }   \hspace{-1pt}= \hspace{-1pt} e^{\f{2}{1-A^{-1}}}.
\) Then since  $\lim_{t\to\infty} t \ov(g_\beta(t)) = 0  $ by $\l(\ref{condt}\r)$,  using that $1+\gamma-\d=0$, as $t\to\infty$,}
%
  (*) &\lesssim       \log(t)^{ 1 +  \l( 1+\eps -\alpha\beta   \r)\kappa -(\beta-\d)\alpha   + (\beta-\delta)c  +\gamma-\d    }
  %
    \label{777}
  =  \log(t)^{    \l( 1+\eps -\alpha\beta   \r)\kappa   -(\beta-\d)\alpha   + (\beta-\delta)c       }.
  %
 %
  %
\end{align}
 Now, $\d<\beta$, so  $-(\beta-\d)\alpha<0$. Choosing $\kappa$, $c$ small enough that the   exponent in $\l(\ref{777}\r)$  is negative,  $\lim_{t\to\infty}(*)=0$. Then by $\l(\ref{S2*bit}\r)$,   uniformly in $(h,y,t) \in \mbox{\hyperref[calS]{$\cal{S}$}}$, $(\mbox{\hyperref[S1*S2*]{$S_2^*$}})\lesssim   t^{-1} \log(t)^{-1-\eps}  \Phi_y^h(t)/(f(y)-h) $,  as required for $\l(\ref{1'}\r)$.

  To prove $(\mbox{\hyperref[S1*S2*]{$S_2^*$}})\leq  \Phi_y^h(t)  u(t) \ov(g(t))  / (f(y)-h)  $, 
    applying Lemma \mbox{\ref{1.1}} and  Lemma \mbox{\ref{keylemma}} with $H(t)= 1/ (t \log(\log(t)))   $ to $\l(\ref{Horatio'scomment}\r)$,   uniformly in $(h,y,t) \in \mbox{\hyperref[calS]{$\cal{S}$}}$,
  \b{align*}
  (\mbox{\hyperref[S1*S2*]{$S_2^*$}})  &\overset{\l(\ref{Horatio'scomment}\r)}\lesssim  t \ov(g(t))    \pp\l(     X_{t}^{(0,  g_{\d}(t))}    >   \f{  (1-A^{-1}) g(t)          }{   2 \log(t)^\gamma  }          \r) 
  \overset{\hspace{1.5pt}\ref{keylemma}\hspace{1.5pt}}\leq        \f{\ov(g(t))}{\log(\log(t))}   \exp \l( (*)      \r)      ,
\\
  (*)  &\lesssim       t \log\l(    t \log(\log(t)) \r)      \l(   t \log(\log(t))  \r)^{      \f{ 2\log(t)^{\gamma-\d}   }{   (1-A^{-1})}     }    \ov(g_\d(t))      \log(t)^{\gamma-\d}        .
  \intertext{Recall $1+\gamma- \d=0$. Noting that $    t   \log(\log(t))   \lesssim   t^2   $  uniformly in $t>t_0(y)>0$,  it follows that }
 (*)   
  &\lesssim      t      \log(t)^{1+\gamma-\d}     t^{   \f{   4   \log(t)^{\gamma-\d}    }{ 1-A^{-1}}   }      \ov(g_\d(t))        
    =       t       e^{  \f{   4 (1-A^{-1})  \log(t)^{1+\gamma-\d}   }{ 1-A^{-1}}   }      \ov(g_\d(t))        
    =      t        e^{  \f{   4 (1-A^{-1})     }{ 1-A^{-1}}   }      \ov(g_\d(t))        .
  %
  \end{align*}
  Now, since $ \d< \beta$, by $\l(\ref{condt}\r)$,  $   \lim_{t\to\infty}t\ov(g_\d(t)) = 0$, and      hence  $\lim_{t\to\infty}(*)=0$, so by Lemma   for suitable $u$,  $(\mbox{\hyperref[S1*S2*]{$S_2^*$}})\leq   \Phi_y^h(t) u(t) \ov(g(t)) / (f(y)-h)  $,  uniformly in $(h,y,t) \in \mbox{\hyperref[calS]{$\cal{S}$}}$, as required for $\l(\ref{2'}\r)$.

 \end{proof}

\section{Proof of Lemma  \hbox{\ref{l6}}} \label{l6proof}    

  \b{proof}[Proof of Lemma \mbox{\ref{l6}}]    
  %
%
%
    Let $T_{g(h)}$ denote the time when $X$  first passes above  $g(h)$, and let $S_{\Delta_1^{g(h)}} 
    $ be the size of $X$'s first jump of size larger than $g(h)$. For each~$y  >  K$,   with $K>0$ a large, fixed constant, 
  \b{align}
  \nonumber
      \pp\l(   X_h \in g(h) dy ;   \cal{O}_h    \r)   
 &=         \pp\l(   X_h \in g(h) dy ;      X_{T_{g(h)}} \leq  \f{g(h)y}{2}      ;  \cal{O}_h    \r)  
\\
\nonumber
 &+      \pp\l(   X_h \in g(h) dy ;      X_{T_{g(h)}} >  \f{g(h)y}{2}     ;   S_{\Delta_1^{g(h)}}  < \f{g(h)y}{2}      ;\cal{O}_h    \r)    
 \\
 \nonumber
 &+      \pp\l(   X_h \in g(h) dy ;      X_{T_{g(h)}} >  \f{g(h)y}{2}      ;  S_{\Delta_1^{g(h)}}  \geq \f{g(h)y}{2}      ; \cal{O}_h    \r) 
 \\
 \label{starting}
 &=:   \sigma_h^1(dy)   +  \sigma_h^2(dy)    + \sigma_h^3(dy).   
   \end{align} 
We     will bound $ \pp\l(   X_h \in g(h) dy ;   \cal{O}_h    \r)   $ by bounding  these 3 terms separately.  

  \p{Upper Bound for \mbox{\hyperref[starting]{$\sigma_h^1(dy)$}}}  We shall disintegrate on the values of $T_{g(h)}$ and $X_{T_{g(h)}}$. Observe by $\l(\ref{o}\r)$ that
  \(
\pp( \cal{O}_h  ; T_{g(h)} \in ds  ) 
 = \pp ( \cal{O}_{ s }  ;  T_{g(h)} \in ds   )
  \),
so that we can  apply  $\l(\ref{s-}\r)$ and the independent increments property, with the notation $\pp(X_t\in dx)=f_t(x)dx  $,  to yield
  \b{align}
\nonumber
  \mbox{\hyperref[starting]{$\sigma_h^1(dy)$}} &=   \int_{s=0}^{h} \int_{w=1}^{\f{y}{2}}        \mathbb{P}\l(   X_h  \in    g(h)  dy ;   X_{T_{g(h)}}    \in g(h) dw ;   T_{g(h)} \in ds  ;  \cal{O}_h  \r)
\\
\nonumber
    &=  \int_{0}^{h} \int_{1}^{\f{y}{2}}       f_{h-s} (g(h)  (y-w))g(h)  dy \pp\l(    X_{T_{g(h)}}    \in g(h) dw ;   T_{g(h)}  \in ds   ;  \cal{O}_s      \r)
    \\
\nonumber
    &=  \int_{0}^{h} \int_{1}^{\f{y}{2}}       f_{h-s} (g(h)  (y-w))g(h)  dy \pp\l(    X_{T_{g(h)}}    \in g(h) dw ;   T_{g(h)}  \in ds   ;  \cal{O}_h      \r).
    \end{align}
Now, $g(h)  (y-w) > g(h)y/2>g(h)+x_0\geq  g(h-s) + x_0$, for all large enough $h$, with   $x_0$ as  in Assumption \mbox{\ref{case1a}}, so   $\l(\ref{l6cond}\r)$ applies to $  f_{h-s} (g(h)  (y-w))$.   Applying $\l(\ref{2.10.3}\r)$, since $ y-w\geq y/2$ and $L$ is \mbox{\hyperref[sv]{slowly varying}} at $\infty$, uniformly in $y>K$ by  \cite[Theorem 1.2.1]{bgt89},  as $h\to\infty$,
\begin{align}
  \nonumber
     \mbox{\hyperref[starting]{$\sigma_h^1(dy)$}}       &\overset{(\ref{l6cond})}\lesssim    \int_{0}^{h}   \int_{1}^{\f{y}{2}}   (h  -  s) u(g(h)(y  -  w)) g(h) dy \pp  \l(    X_{T_{g(h)}}    \in g(h) dw ;    T_{g(h)} \in ds;  \cal{O}_h    \r)
\\
 \nonumber
    &\overset{(\ref{2.10.3})}\lesssim   \int_{0}^{h}   \int_{1}^{\f{y}{2}}  \f{(h  -  s)}{ g(h)^{\alpha} }  (y  -  w)^{-1-\alpha} L(g(h)(y-w)) dy 
  \pp  \l(     X_{T_{g(h)}}      \in g(h) dw ;    T_{g(h)} \in ds;  \cal{O}_h     \r)
  \\
  \nonumber
   &\overset{\phantom{(\ref{l6cond})}}\lesssim   \int_{0}^{h}       \int_{1}^{\f{y}{2}}  \f{h}{ g(h)^{\alpha}}\f{L(g(h))}{L(g(h))} y^{-1-\alpha} L\l(\f{g(h)y}{2}\r) dy 
  \pp  \l(    X_{T_{g(h)}}    \in g(h) dw ;   T_{g(h)} \in ds  ;  \cal{O}_h  \r)
\\
  \nonumber
  &\overset{\phantom{(\ref{l6cond})}}\lesssim   y^{-1-\alpha}\f{ L\l(g(h)y\r)}{L(g(h))} dy   \f{h L(g(h)) }{ g(h)^{\alpha}}     \int_{0}^{h}   \int_{1}^{\f{y}{2}}      
 \pp  \l(    X_{T_{g(h)}}     \in   g(h) dw ;   T_{g(h)}   \in     ds ;  \cal{O}_h   \r)
\\
  &\overset{\phantom{(\ref{l6cond})}}\leq    y^{-1-\alpha}\f{ L\l(g(h)y\r)}{L(g(h))} dy \    h   \ov(g(h)) \ 
 \pp\l(  \cal{O}_h    \r)
 \label{sigma1bound}
 \overset{(\ref{condt})}=   o(1) \times    y^{-1-\alpha}\f{ L\l(g(h)y\r)}{L(g(h))} \pp(\cal{O}_h) dy    ,
\end{align}
   where, recalling   $ g(h)^{-\alpha} L(g(h))=   \ov(g(h)) $,   the last step follows by~$\l(\ref{condt}\r)$.

  \p{Simplifying the Expressions for \mbox{\hyperref[starting]{$\sigma_h^2(dy)$}} and \mbox{\hyperref[starting]{$\sigma_h^3(dy)$}}}  
     Recall     $\l(\ref{D}\r)$, and      $\mbox{\hyperref[starting]{$\sigma_h^2(dy)$}} + \mbox{\hyperref[starting]{$\sigma_h^3(dy)$}}    = \pp (
   X_h \in g(h) dy ;      X_{T_{g(h)}}   >    g(h)y/2    ;  \cal{O}_h    )$.  
   Choosing $K>4$, we have     $g(h)y/2   >   2  g(h) $  for  $y  >  K$. 
  As $T_{g(h)}$ is the first passage time above  $g(h)$,   if $X_{T_{g(h)}}   >   2g(h)$, then  $X$   crosses $g(h)$ by a jump larger than $g(h)$, so
       since $T_{g(h)}  \leq  \Delta_1^{g(h)}$, $T_{g(h)}  =   \Delta_1^{g(h)}$.

  Then since $X_t<g(h)$ for all $t <  T_{g(h)}$, it follows that $X_{T_{g(h)} -} = X_{\Delta_1^{g(h)}  -}   < g(h)$, as $X$   has c\`adl\`ag sample paths, almost surely.
  Moreover, if $\cal{O}_h$ holds, then $X$   crosses $g(h)$ by time $h$, so  $\Delta_1^{g(h)} = T_{g(h)}   \leq h$.  Thus:
     \begin{equation}
     \label{firstinclusion}
   \l\{     X_{T_{g(h)}}   >   \f{g(h)y}{2}     ;   S_{\Delta_1^{g(h)}}   <  \f{g(h)y}{2}      ;\cal{O}_h     \r\}    \subseteq             \l\{     \Delta_1^{g(h)}     \leq    h ;       X_{\Delta_1^{g(h)}  -}    < g(h)        ;   S_{\Delta_1^{g(h)}}   <     \f{g(h)y}{2}      ;\cal{O}_h       \r\}    ,
  \end{equation}
  and therefore  we can bound $  \sigma_h^2(dy) $ by
  \b{align}
  \nonumber
  \mbox{\hyperref[starting]{$\sigma_h^2(dy)$}}   &=      \pp\l(   X_h \in g(h) dy ;      X_{T_{g(h)}} >  \f{g(h)y}{2}     ;   S_{\Delta_1^{g(h)}}  < \f{g(h)y}{2}      ;\cal{O}_h    \r)  
\\
\label{sigma2}
  &\leq         \pp \l(  \hspace{-2pt} X_h \in g(h) dy ;    \Delta_1^{g(h)}   \leq h ;       X_{\Delta_1^{g(h)}  -}  \hspace{-2pt} <   \hspace{-2pt}  g(h)        ;   S_{\Delta_1^{g(h)}} \hspace{-2pt} < \hspace{-2pt} \f{g(h)y}{2}      ;\cal{O}_h \hspace{-2pt}   \r)  \hspace{-2pt}  .  
  \end{align}
 For \mbox{\hyperref[starting]{$\sigma_h^3(dy)$}},   the converse analogous inclusion to $\l(\ref{firstinclusion}\r)$ holds too, that is, if  $ \cal{O}_h $, $    \Delta_1^{g(h)}   \leq h$,         $  X_{\Delta_1^{g(h)}  -}  \hspace{-0.15cm}    <   \hspace{-0.05cm}   g(h)   $, and $S_{\Delta_1^{g(h)}}  \geq g(h)y/2$        hold, then we have   
%
%
  %
  %
    %
    %
       \[
   X_{T_{g(h)}} =   X_{\Delta_1^{g(h)}}  \geq   X_{\Delta_1^{g(h)}}  -   X_{\Delta_1^{g(h)}-}   = S_{\Delta_1^{g(h)}} >   g(h)y/2  ,
   \]
  and therefore   $\sigma_h^3(dy)  $ satisfies
     \b{align}
    \nonumber    
    \mbox{\hyperref[starting]{$\sigma_h^3(dy)$}}    &=      \pp\l(   X_h \in g(h) dy ;      X_{T_{g(h)}} >  \f{g(h)y}{2}      ;  S_{\Delta_1^{g(h)}}  \geq \f{g(h)y}{2}      ; \cal{O}_h    \r) 
  \\
    \label{contain} 
    &=        \pp    \l(\hspace{-0.1cm}    X_h \in g(h) dy ;       \Delta_1^{g(h)}   \leq h ;       X_{\Delta_1^{g(h)}  -}  \hspace{-0.15cm}    <   \hspace{-0.05cm}   g(h)        ;  S_{\Delta_1^{g(h)}}  \hspace{-0.05cm} \geq \hspace{-0.05cm} \f{g(h)y}{2}      ; \cal{O}_h  \hspace{-0.075cm}  \r) .
  \end{align}

    \p{Upper Bound for \mbox{\hyperref[starting]{$\sigma_h^2(dy)$}}}
     By $\l(\ref{sigma2}\r)$ and $\l(\ref{ttos}\r)$, disintegrating on the values of $\Delta_1^{g(h)}  $, $ X_{\Delta_1^{g(h)}-}  $, and   $ S_{\Delta_1^{g(h)}}   $, 
        by  independence of increments and the Markov property, with $\pp(X_t\in dx)=f_t(x)dx  $, 
  \b{align}
  \nonumber
  \mbox{\hyperref[starting]{$\sigma_h^2(dy)$}} &\overset{\l(\ref{sigma2}\r)}\leq       \pp\l(   X_h \in g(h) dy ;    \Delta_1^{g(h)}   \leq h ;       X_{\Delta_1^{g(h)}  -}   < g(h)        ;   S_{\Delta_1^{g(h)}}  < \f{g(h)y}{2}      ;\cal{O}_h    \r) 
\\
\nonumber
  &\hspace{2pt}\overset{\l(\ref{ttos}\r)}=   \int_{s=0}^{h}     \int_{w=0}^{1}    \int_{    v=0}^{\f{y}{2}}          \pp\big(   X_h \in g(h) dy ;    \Delta_1^{g(h)}   \in ds ;       X_{\Delta_1^{g(h)}  -}   \in g(h) dw         ;  
  \\
  \nonumber
  &\ \hspace{96pt}  S_{\Delta_1^{g(h)}}  \in g(h) dv    ;    \cal{O}_s    \big) 
 \\
 \nonumber
  &\overset{\phantom{\l(\ref{sigma2}\r)}}= \int_{0}^{h}     \int_{0}^{1}    \int_0^{\f{y}{2}}      f_{h-s} ( g(h) (y-w-v) ) g(h)   dy
  \\
  \nonumber   &\    \hspace{2.3cm}  \times  \pp\l(    \Delta_1^{g(h)}   \in ds ;       X_{\Delta_1^{g(h)}  -}   \in g(h) dw         ;   S_{\Delta_1^{g(h)}}  \in g(h) dv    ;\cal{O}_s    \r)
  \\
 \nonumber
  &\hspace{2pt}\overset{\l(\ref{ttos}\r)}= \int_{0}^{h}     \int_{0}^{1}    \int_0^{\f{y}{2}}      f_{h-s} ( g(h) (y-w-v) ) g(h)   dy
  \\
  \nonumber   &\    \hspace{2.3cm}  \times  \pp\l(    \Delta_1^{g(h)}   \in ds ;       X_{\Delta_1^{g(h)}  -}   \in g(h) dw         ;   S_{\Delta_1^{g(h)}}  \in g(h) dv    ;\cal{O}_h    \r).
  \end{align}
Note    $y-w-v > y/3>K/3$ for $w\leq1$,  $ v\leq y/2$.  So  as $h\to\infty$, $g(h)(y-w-v) \geq g(h-s) + x_0$,    so    we can apply  $\l(\ref{l6cond}\r)$ and $\l(\ref{2.10.3}\r)$. Now, $g(h)^{-\alpha} L(g(h))  =  \ov(g(h))$  for   $L$   \mbox{\hyperref[sv]{slowly varying}} at $\infty$,    so by $\l(\ref{condt}\r)$,  uniformly in $y  >  K$   by the uniform convergence theorem  \cite[Theorem 1.2.1]{bgt89},~as~$h  \to  \infty$,     
\begin{align}
\nonumber
  \mbox{\hyperref[starting]{$\sigma_h^2(dy)$}}  &\overset{\l(\ref{l6cond}\r)}\lesssim \int_{0}^{h}     \int_{0}^{1}    \int_0^{\f{y}{2}}    (h-s) u (g(h)(y-w-v))  g(h) dy
  \\  
   \nonumber 
  &\  \hspace{2cm}    \times  \pp\l(    \Delta_1^{g(h)}   \in ds ;       X_{\Delta_1^{g(h)}  -}   \in g(h) dw         ;   S_{\Delta_1^{g(h)}}  \in g(h) dv    ;\cal{O}_h    \r)
%
%
\\
    \nonumber
       &\overset{\l(\ref{2.10.3}\r)}\lesssim   \int_{0}^{h}     \int_{0}^{1}    \int_0^{\f{y}{2}}     \f{   (h-s)      L(g(h)(y-w-v))      }{g(h)^{\alpha}  (y-w-v)^{1+\alpha}}   dy
  \\      \nonumber
    &\  \hspace{2cm}    \times  \pp\l(    \Delta_1^{g(h)}   \in ds ;       X_{\Delta_1^{g(h)}  -}   \in g(h) dw         ;   S_{\Delta_1^{g(h)}}  \in g(h) dv    ;\cal{O}_h    \r)
  \\
    \nonumber
       &\overset{\phantom{\l(\ref{2.10.3}\r)}}\lesssim   \int_{0}^{h}     \int_{0}^{1}    \int_0^{\f{y}{2}}  \f{   (h-s) }{g(h)^{\alpha}  y^{1+\alpha}}  L(g(h)y) dy
  \\      \nonumber
    &\  \hspace{2cm}    \times  \pp\l(    \Delta_1^{g(h)}   \in ds ;       X_{\Delta_1^{g(h)}  -}   \in g(h) dw         ;   S_{\Delta_1^{g(h)}}  \in g(h) dv    ;\cal{O}_h    \r)
     \\       \nonumber
     &\overset{\phantom{\l(\ref{l6cond}\r)}}\leq     \f{   h   L(g(h)) }{g(h)^{\alpha}  }   \int_{0}^{h}     \int_{0}^{1}    \int_0^{\f{y}{2}}  y^{-1-\alpha}    \f{L(g(h)y) }{L(g(h)} dy
  \\      \nonumber  
  &\  \hspace{2cm}    \times  \pp\l(    \Delta_1^{g(h)}   \in ds ;       X_{\Delta_1^{g(h)}  -}   \in g(h) dw         ;   S_{\Delta_1^{g(h)}}  \in g(h) dv    ;\cal{O}_h    \r)
%
 \\
         \nonumber
  &\hspace{2pt}  \overset{\l(\ref{condt}\r)}=    o(1) \times    \int_{0}^{h}     \int_{0}^{1}    \int_0^{\f{y}{2}}       y^{-1-\alpha}    \f{L(g(h)y) }{L(g(h)} dy
  \\   \nonumber
   &\  \hspace{2cm}    \times  \pp\l(    \Delta_1^{g(h)}   \in ds ;       X_{\Delta_1^{g(h)}  -}   \in g(h) dw         ;   S_{\Delta_1^{g(h)}}  \in g(h) dv    ;\cal{O}_h    \r)
   \\
      %
      %
%
%
\label{sigma2bound}
  &\overset{\phantom{\l(\ref{l6cond}\r)}}\leq       o(1) \times  y^{-1-\alpha}  \f{L(g(h)y)}{L(g(h))} \pp(\cal{O}_h)  dy.
  \end{align}     

   \b{comment}
   \p{Summary of Upper Bounds}    We have shown that as $h\to\infty$,
   \[
   \int_K^\infty    q_h(y g(h))     \sigma_h(dy) \lesssim   o(1) \times \pp(\cal{O}_h)   \int_{y=\lfloor{K}\rfloor}^\infty         q_h(yg(h))        y^{-1-\alpha}     dy
   \]
   \[
   +     2   \pp\l(   \Delta_1^{g(h)} \leq h      ;   X_{\Delta_1^{g(h)} -}  < g(h)          ; \cal{O}_h    
       \r)   \int_{y=\lfloor{K}\rfloor}^\infty         q_h(yg(h))        y^{-1-\alpha}  \f{L(g(h)    y)}{L(g(h))}    dy
    \] 
   \[
   \lesssim  \pp\l( \cal{O}_h   \r) \int_{y=\lfloor{K}\rfloor}^\infty         q_h(yg(h))        y^{-1-\alpha}    \f{L(g(h)    y)}{L(g(h))}    dy.
   \]
  \end{comment}

   \p{Upper Bound for \mbox{\hyperref[starting]{$\sigma_h^3(dy)$}}}   Disintegrating on the values of $\Delta_1^{g(h)} $, $ X_{\Delta_1^{g(h)}-} $ and  $S_{\Delta_1^{g(h)}}  $,   then applying $\l(\ref{ttos}\r)$, independence of increments, the Markov property, 
  and Lemma \mbox{\ref{bigjump}}, it follows that uniformly among $y>K$ as $h\to\infty$, with $\pp(X_t\in dx)   =   f_t(x)dx  $,
    \b{align}
        \nonumber
  \mbox{\hyperref[starting]{$\sigma_h^{3}(dy)$}} &\overset{\phantom{\ref{bigjump}}}=     \int_{s=0}^h \int_{w=0}^1 \int_{v=\f{y}{2}}^{y-w}  \pp\Big(   X_h \in g(h) dy ;   \Delta_1^{g(h)}   \in ds  ; 
  \\
        \nonumber
  &\ \hspace{3.4cm} \hspace{3pt}     X_{\Delta_1^{g(h)}  -}   \in  g(h)dw     ;   S_{\Delta_1^{g(h)}}   \in   g(h)dv      ;\cal{O}_h    \Big) 
\\
        \nonumber
  &\hspace{-2pt}\overset{\l(\ref{ttos}\r)}=  \int_0^h \int_0^1 \int_{\f{y}{2}}^{y-w} f_{h-s}  ( g(h) (y-w-v) ) g(h)dy  \pp\l(  S_{\Delta_1^{g(h)}}   \in   g(h)dv     \r)
  \\
        \nonumber
  &\ \hspace{2.5cm} \times \pp\l(    \Delta_1^{g(h)}   \in ds  ;       X_{\Delta_1^{g(h)}  -}   \in  g(h)dw        ;\cal{O}_s    \r) 
%
%
   %
   %
\\
        \nonumber
     &\overset{\ref{bigjump}}\lesssim    \int_0^h \int_0^1 \int_{\f{y}{2}-w}^{y-w}      f_{h-s}  ( g(h) (y-w-v) ) g(h)dy               \  v^{-1-\alpha}  \f{L(g(h)v)}{    L(g(h))} dv
\\
        \nonumber
  &\ \hspace{2.5cm}   \times  
      \pp\l(    \Delta_1^{g(h)}   \in ds  ;       X_{\Delta_1^{g(h)}  -}   \in  g(h)dw        ;\cal{O}_s    \r)
%
%
      %
      %
      \\
        \nonumber
      &\hspace{-2pt}\overset{\l(\ref{ttos}\r)}=     \int_0^h \int_0^1 \int_{\f{y}{2} -w}^{y- w}    f_{h-s}   ( g(h) (y-w-v) ) g(h) dy \   v^{-1-\alpha}   \f{L(g(h)v)}{     L(g(h))} dv
     \\
        \nonumber
     &\ \hspace{2.5cm} 
     \times
      \pp\l(    \Delta_1^{g(h)}   \in ds  ;       X_{\Delta_1^{g(h)}  -}   \in  g(h)dw        ;\cal{O}_h    \r).
\intertext{Now, as $ y/3 \leq y/2-1\leq  y/2 -w \leq v \leq y$, applying the uniform convergence theorem \cite[Theorem 1.2.1]{bgt89} to $L(g(h)v)/L(g(h)y)$,  uniformly in $y>K$ as $h\to\infty$,}
\mbox{\hyperref[starting]{$\sigma_h^{3}(dy)$}} &\lesssim y^{-1-\alpha}   \f{L(g(h)y)}{     L(g(h))} dy    \int_0^h \int_0^1 \int_{\f{y}{2} -w}^{y- w}    f_{h-s}   ( g(h) (y-w-v) ) g(h)   dv
     \\
        \nonumber
     &\ \hspace{4.4cm} 
     \times
      \pp\l(    \Delta_1^{g(h)}   \in ds  ;       X_{\Delta_1^{g(h)}  -}   \in  g(h)dw        ;\cal{O}_h    \r).
   \intertext{Changing variables to $u  =   g(h)(y-w-v)$,    uniformly in~$y  >  K$,~as~$h  \to  \infty$, }
        \nonumber
     \mbox{\hyperref[starting]{$\sigma_h^{3}(dy)$}}   &\lesssim  y^{-1-\alpha}    \f{L(g(h)y)}{      L(g(h))} dy    \int_0^h \int_0^1 \int_{0}^{\f{g(h)y}{2} }   f_{h-s}  ( u ) du           
\\
        \nonumber
     &\ \hspace{4.4cm}     \times
      \pp\l(    \Delta_1^{g(h)}   \in ds  ;       X_{\Delta_1^{g(h)}  -}   \in  g(h)dw        ;\cal{O}_h    \r)
   \\
        \nonumber
   &=  y^{-1-\alpha}   \f{L(g(h)y)}{      L(g(h))} dy    \int_0^h \int_0^1    \pp\l(  X_{h-s} \leq \f{g(h)y}{2} \r)            
\\
        \nonumber
     &\ \hspace{4.4cm}     \times
      \pp\l(    \Delta_1^{g(h)}   \in ds  ;       X_{\Delta_1^{g(h)}  -}   \in  g(h)dw        ;\cal{O}_h    \r)
        \\
        \nonumber
   &\leq  y^{-1-\alpha}   \f{L(g(h)y)}{      L(g(h))} dy    \int_0^h \int_0^1   
      \pp\l(    \Delta_1^{g(h)}   \in ds  ;       X_{\Delta_1^{g(h)}  -}   \in  g(h)dw        ;\cal{O}_h    \r)
     \\
       \label{penultimatebound}
         &=  y^{-1-\alpha}   \f{L(g(h)y)}{      L(g(h))} dy    
      \pp\l(    \Delta_1^{g(h)}   \leq h  ;       X_{\Delta_1^{g(h)}  -}   < g(h)       ;\cal{O}_h    \r)
          \\
            \label{sigma3bound}
         &\leq  y^{-1-\alpha}    \f{L(g(h)y)}{    L(g(h))}     \pp\l(    \cal{O}_h    \r) dy    
    .
\end{align}

    \p{Conclusion of Upper Bound} By  
   $\l(\ref{starting}\r)$,   $\l(\ref{sigma1bound}\r)$, $\l(\ref{sigma2bound}\r)$, and $\l(\ref{sigma3bound}\r)$,          we conclude, as required for the upper bound in $\l(\ref{l6statement}\r)$, that  uniformly in $y>K$,    as $h\to\infty$,    
   \b{equation} \label{sigmaupper}
   \pp(X_h \in g(h) dy    ;   \cal{O}_h ) \lesssim     y^{-1-\alpha}  \f{  L(g(h)y)}{ L(g(h))}   \pp(\cal{O}_h)  dy.
     \end{equation}
 %
  %
 %
 %
  %
  %
 Now we will prove the lower bound  on $  \pp\l(       X_h \in g(h)dy   ; \cal{O}_h  \r) $.

\p{Proof of Lower Bound}   Now, fixing  $y_0>0$,     for  all  $y>K$, as $h\to\infty$,
\begin{align}
\nonumber
     \pp\l(       X_h \in g(h)dy   ; \cal{O}_h  \r)    &\geq   \pp\Big(       X_h \in g(h)dy   ;      \D^{g(h)}   \leq h-1   ;   X_{\Delta_1^{g(h)}-}   < g(h)    ;    \hspace{3cm}  
\\
     \label{lb}
   &\hspace{3.25cm}         \f{g(h)y}{2}  \leq S_{\Delta_1^{g(h)}}\leq    g(h)y-y_0      ;    \cal{O}_h   \hspace{-1pt}  \Big)              .
\end{align}
      %
%
%
%
 %
   %
   %
%
%
  %
Disintegrating on the values of $  \Delta_1^{g(h)} $, $ X_{\Delta_1^{g(h)}-}$, and $S_{\Delta_1^{g(h)}}  $,       applying the Markov property, noting that by $\l(\ref{ttos}\r)$,  $\{ \Delta_1^{g(h)} \hspace{-2pt}  \in  \hspace{-2pt} ds ;  \cal{O}_h \} \hspace{-2pt}=\hspace{-2pt} \{ \Delta_1^{g(h)} \hspace{-2pt}  \in  \hspace{-2pt} ds ;  \cal{O}_{\Delta_1^{g(h)}} \}$  for each $s\leq h$,    
with~$\pp(X_t \in dx) \hspace{-2pt}=\hspace{-2pt} f_t(x)dx$,
\b{align} 
   \nonumber
    \l(\ref{lb}\r)     
%
 %
%
    &\overset{\l(\ref{ttos}\r)}=       \int_{s=0}^{h-1}   \int_{w=0}^1   \int_{v=\f{y}{2}}^{y-\f{y_0}{g(h)}-w }      \pp\Big(    X_{h} \in g(h)dy ;   \Delta_1^{g(h)} \in ds  ;   X_{\Delta_1^{g(h)}-}   \in g(h) dw   ;
    \\
    \nonumber
        &\ \hspace{4.6cm}    S_{\Delta_1^{g(h)}}  \in g(h) dv    ; \cal{O}_{\Delta_1^{g(h)}}      \Big) 
        \\
    \nonumber
     &\overset{\l(\ref{ttos}\r)}=           \int_{0}^{h-1}   \int_{0}^1   \int_{v=\f{y}{2}}^{y-\f{y_0}{g(h)}-w}   f_{h-s}(g(h)(y-w-v))  g(h) dy  \pp\l( S_{\Delta_1^{g(h)}}  \in g(h) dv        \r)
   \\
    \nonumber
    &\ \hspace{3.25cm} \times \pp\l(   \Delta_1^{g(h)} \in ds  ;   X_{\Delta_1^{g(h)}-}   \in g(h) dw   ; \cal{O}_{h}  \r) .
  %
  %
  \intertext{Applying Lemma \mbox{\ref{bigjump}},  noting    $h-s\geq1$,   $y/2 <  2y/3 -w$, $v\asymp y$, and    $   L(g(h)v)\hspace{-1pt}  \asymp \hspace{-1pt}    L(g(h)y)$ uniformly in $y>K$ as $h\to\infty$ by the uniform convergence theorem \cite[Thm 1.2.1]{bgt89}, it follows that uniformly in~$y  >  K$~as~$h  \to  \infty$,  }
   \nonumber
     \l(\ref{lb}\r)   &\overset{\ref{bigjump}}\gtrsim          \int_{0}^{h-1}   \int_{0}^1   \int_{v=\f{y}{2}}^{y-\f{y_0}{g(h)} - w}    f_{h-s} (g(h)(y-w-v) ) g(h)   dy        v^{-1-\alpha}  \f{L(g(h)v)}{L(g(h))}        dv
   \\
    \nonumber
    &\ \hspace{3.25cm} \times \pp\l(   \Delta_1^{g(h)} \in ds  ;   X_{\Delta_1^{g(h)}-}   \in g(h) dw   ; \cal{O}_h  \r)  
  %
  \\
    \nonumber
     &\gtrsim          \int_{0}^{h-1}   \int_{0}^1   \int_{v=\f{2y}{3}-w}^{y-\f{y_0}{g(h)} - w}    f_{h-s} (g(h)(y-w-v) )     g(h)   dy        y^{-1-\alpha}     \f{L(g(h)y)}{L(g(h))}          dv
   \\
    \nonumber
  &\ \hspace{3.25cm} \times \pp\l(   \Delta_1^{g(h)} \in ds  ;   X_{\Delta_1^{g(h)}-}   \in g(h) dw   ; \cal{O}_h  \r) .
       \intertext{Changing   variables   to $u =    g(h) (y -w-v)  $, noting that  $y/3 > 1$ for all  $y>K$ and that $h-s\geq1$, 
    }    
   %
   \nonumber
             \l(\ref{lb}\r)        &\gtrsim  
       y^{-1-\alpha}      \f{L(g(h)y)}{L(g(h))}         dy   \int_{0}^{h-1}   \int_{0}^1          \int_{u=y_0}^{\f{g(h)y}{3 } }      f_{h-s}(u) du 
   \\
   \nonumber
  &\ \hspace{4.35cm} \times \pp\l(   \Delta_1^{g(h)} \in ds  ;   X_{\Delta_1^{g(h)}-}   \in g(h) dw   ; \cal{O}_h  \r)
   \\
  \nonumber   
  &=      y^{-1-\alpha}      \f{L(g(h)y)}{L(g(h))}         dy   \int_{0}^{h-1}   \int_{0}^1    \l[     \pp\l(  X_{h-s}\leq  \f{g(h)y}{3}    \r) -     \pp\l(  X_{h-s}\leq y_0    \r) \r]
   \\
   \nonumber
  &\ \hspace{4.35cm} \times \pp\l(   \Delta_1^{g(h)} \in ds  ;   X_{\Delta_1^{g(h)}-}   \in g(h) dw   ; \cal{O}_h  \r)
   \\
  \nonumber   
  &\geq      y^{-1-\alpha}      \f{L(g(h)y)}{L(g(h))}         dy   \int_{0}^{h-1}   \int_{0}^1   \l[     \pp\l(  X_{h}\leq g(h)    \r) -     \pp\l(  X_{1}\leq y_0    \r)   \r]
   \\
   \nonumber
  &\ \hspace{4.35cm} \times \pp\l(   \Delta_1^{g(h)} \in ds  ;   X_{\Delta_1^{g(h)}-}   \in g(h) dw   ; \cal{O}_h  \r).
  \end{align}
%
  %
Now,  with $X^{(0,g(h))}$ again denoting the process with no jumps bigger than $g(h)$,
  \[
  \pp\l(   X_h \leq    g(h) \r)  =   \pp(   X_h^{(0,g(h))} \leq    g(h) ) \pp( \Delta_1^{g(h)} > h) =   \pp(   X_h^{(0,g(h))} \leq    g(h) ) e^{- h \ov(g(h))}  ,
  \]
  and   since $\lim_{h\to\infty} h \ov(g(h))=0$ by $\l(\ref{condt}\r)$, by Markov's inequality,  as $h\to\infty$,
  \[
 \pp\l(   X_h \leq    g(h) \r)    \overset{\l(\ref{condt}\r)}\sim  \      \pp(   X_h^{(0,g(h))} \leq    g(h) ) \   \geq 1 - \f{ \mathbb{E}[   X_h^{(0,g(h))}  ]   }{   g(h)    }    \geq      1 - \f{ h \int_0^{g(h)} \ov(x)dx   }{   g(h)    }.
  \] 
  Now, by $\l(\ref{condt}\r)$ and Karamata's theorem  \cite[Prop 1.5.8]{bgt89}, as $h\to\infty$,
  \[
   \pp\l(   X_h \leq    g(h) \r)     \gtrsim        1    -  \f{    h g(h) \ov(g(h))    }{   g(h)   }   = \  1 - h \ov(g(h))  \  \overset{\l(\ref{condt}\r)}\sim  1.   \hspace{66pt}  
  \]
    Then as $\pp(X_1\leq y_0) = \textup{constant} <1$, taking $y_0$  large enough that $  \pp\l(   X_h \leq    g(h) \r)   - \pp(X_1\leq y_0)  \gtrsim 1$ uniformly, we get that  uniformly in $y>K$ as $h\to\infty$, 
  \begin{align}
   \nonumber
      \l(\ref{lb}\r)   &\gtrsim         \f{L(g(h)y)}{L(g(h))}         y^{-1-\alpha}      dy   \int_{s=0}^{h-1}   \int_{w=0}^1       \pp\l(   \Delta_1^{g(h)} \in ds  ;   X_{\Delta_1^{g(h)}-}   \in g(h) dw   ; \cal{O}_h  \r)
  %
  %
       \\
       \label{lowerbound}
       &=     \f{L(g(h)y)}{L(g(h))}        y^{-1-\alpha}      dy     \pp\l(    \Delta_1^{g(h)}      \leq h-1    ;   X_{\Delta_1^{g(h)}-}     <g(h)     ; \cal{O}_h    \r)   .    
        \end{align}

    \b{comment}   
       \p{Remark}  The powers highlighted in \textcolor{red}{red} are subject to some dispute/confusion. If we work at the level of densities, ignoring any potential $1-\rho$ term introduced by changing from $d(y-z)$ to $d((1-\rho))$, then we get the powers highlighted in red, and the argument works as above. 
        
         If we do include the disputed $\textcolor{pink}{1-\rho}$ term in the above calculations, then since $\rho<1 - \f{1}{g(h)y} $ and $\lim_{h\to\infty} \f{1}{g(h)}= 0$, we can bound $1-\rho \geq   g(h) y \geq  20 g(h) \geq 20$, for all sufficiently large $h$, and the proof as above works. 
         
         However, if we choose to include the disputed $\textcolor{pink}{1-\rho}$ term here, and also for the upper bound, then we end up with an upper bound which is strictly smaller than our lower bound, which is a contradiction. 
         
         This fact alone confirms that we can ignore the $1-\rho$   term introduced by changing from $d(y-z)$ to $d((1-\rho)y)=(1-\rho)dy$, and so we should write the proof using densities rather than densities multiplied by $dy$.
        \end{comment}

\p{Proof by Contradiction Step}    
Now we assume for a contradiction that    
  \b{equation} 
  \label{contradiction0}
  \liminf_{h\to\infty}     \f{  \pp\l(    \Delta_1^{g(h)} \leq h-1  ;   X_{\Delta_1^{g(h)}-}   <g(h)   ; \cal{O}_h    \r)    }{   \pp(\cal{O}_h)     }      = 0 .
   \end{equation}
   As    $\Delta_1^{g(h)} \hspace{-1.5pt} $ is exponentially distributed with rate $\ov(g(h)) \hspace{-0.5pt} $, by Corollary \mbox{\ref{lemma1}},   as~$h      \to      \infty \hspace{-0.5pt}$,
   \begin{align*}
  &  \pp\l(    \Delta_1^{g(h)} \in[ h-1 ,h] ;   X_{\Delta_1^{g(h)}-}   <g(h)   ; \cal{O}_h    \r)  \leq  \pp\l(    \Delta_1^{g(h)} \in[ h-1 ,h]       \r) 
\\
   & \hspace{-4pt}\leq    \pp\l(    \Delta_1^{g(h)} \leq 1 \r)     \hspace{-1pt}   =    \hspace{-1pt}    1     \hspace{-1pt}    - e^{- \ov(g(h))}   \hspace{-1pt}    \leq     \hspace{-1pt}    \ov(g(h)) \overset{\ref{lemma1}}\sim   \f{\pp(\cal{O}_h )}{   \Phi(h)} \overset{\l(\ref{IfPhi}\r)}= o(1) \times \pp\l( \cal{O}_h\r)   ,
   \end{align*}
   since  $\lim_{h\to\infty} \Phi(h)=\infty$ by $\l(\ref{IfPhi}\r)$, so it follows that $\l(\ref{contradiction0}\r)$ holds if and only if
    \b{equation} 
  \label{contradiction}
  \liminf_{h\to\infty}     \f{  \pp\l(    \Delta_1^{g(h)} \leq h  ;   X_{\Delta_1^{g(h)}-}   <g(h)   ; \cal{O}_h    \r)    }{   \pp(\cal{O}_h)     }      = 0 .
   \end{equation}
By  
  $\l(\ref{starting}\r)$,    $\l(\ref{sigma1bound}\r)$, $\l(\ref{sigma2bound}\r)$,  and $\l(\ref{penultimatebound}\r)$, we get that $\l(\ref{contradiction}\r)$ implies, along a subsequence of $h$, as $h\to\infty$,
  %
 %
 %
  %
  %
  \b{align} 
\nonumber
  \pp  \l(   X_h    \geq      Kg(h) ;    \cal{O}_h    \r) &    =     \int_K^\infty        \pp(  X_h \in g(h) dy ; \cal{O}_h
     )      =       o(1)     \times       \pp(\cal{O}_h)        \int_{ K}^{ \infty}     y^{-1-\alpha  }       \f{L(g(h)y)}{L(g(h))}  dy      .
  \intertext{Changing variables from $y$ to $u=g(h)y$,  }
\nonumber
 \pp  \l(   X_h  \hspace{-2pt}  \geq    \hspace{-2pt}  Kg(h) ;    \cal{O}_h    \r) &  =   o(1) \times   \f{   g(h)^\alpha  \pp(\cal{O}_h)  }{     L(g(h))   }   \int_{Kg(h)}^{\infty}       u^{-1-\alpha  }       L(u)  du.
  \intertext{As $L$ is \mbox{\hyperref[sv]{slowly varying}}, applying the result   \cite[Prop 1.5.10]{bgt89} to $ \int_{Kg(h)}^{\infty}       u^{-1-\alpha  }       L(u)  du$, as $h\to\infty$, }
  %
  \pp\l(   X_h \geq K g(h) ;   \cal{O}_h    \r) &\lesssim     o(1) \times   \f{   g(h)^\alpha  \pp(\cal{O}_h)  }{       L(g(h))   }   (Kg(h))^{-\alpha}           L(Kg(h))  
      \label{finall}
     = o(1) \times  \pp(\cal{O}_h)  .
 \intertext{But considering the subevent $ \{ \Delta_1^{g(h)}=\Delta_1^{Kg(h)}\leq h ; \cal{O}_h\}  \subseteq \{  X_h \geq K g(h) ;   \cal{O}_h     \}$,   disintegrating on the value of $\Delta_1^{g(h)} $, and applying  the Markov property,
}
   \nonumber
   \pp\l(   X_h \geq K g(h) ;   \cal{O}_h    \r) &\geq \pp\l(\Delta_1^{g(h)}=\Delta_1^{Kg(h)}\leq h ; \cal{O}_h\r)
  %
  %
 =     \int_0^h         \pp\l(\Delta_1^{g(h)}=\Delta_1^{Kg(h)}\in ds  ; \cal{O}_h\r)
  \\
   \nonumber
   &=      \int_0^h         \pp\l(\Delta_1^{g(h)} \in ds  ;    S_{\Delta_1^{g(h)}} \geq K g(h)      ; \cal{O}_{s,X^{(0,g(h))}}  \r)
  \\
   \nonumber
   &=\int_0^h        \pp\l(  \cal{O}_{s,X^{(0,g(h))}}     ;    \Delta_1^{g(h)} \in ds   \r)   \pp\l(  S_{\Delta_1^{g(h)}} \geq K g(h)  \r)
  \\
   \nonumber
       &=     \f{  \ov(Kg(h))  }{    \ov(g(h))  }   \int_0^h        \pp\l(  \cal{O}_{s,X^{(0,g(h))}}     ;    \Delta_1^{g(h)} \in ds   \r)  
  \\
   \nonumber
       &=       \f{  \ov(Kg(h))  }{    \ov(g(h))  }   \int_0^h        \pp\l(     \cal{O}_{h}     ;    \Delta_1^{g(h)}    \in    ds      \r)  
   %
  %
     =       \f{  \ov(Kg(h))  }{    \ov(g(h))  }  \pp\l(     \cal{O}_{h}     ;    \Delta_1^{g(h)} \leq h    \r)  .
\intertext{By $\l(\ref{final123}\r)$ with $g_y^h(t)=g(t)$, and by Corollary \mbox{\ref{lemma1}},    $ \pp(  \cal{O}_{h}     ;    \Delta_1^{g(h)}     \leq   h )     \sim   \pp(\cal{O}_h)$ as $h  \to  \infty$, so~as~$h  \to  \infty$,      }
\label{contr2}
       \pp\l(   X_h \geq K g(h) ;   \cal{O}_h    \r) 
        &\geq      \f{  \ov(Kg(h))  }{    \ov(g(h))  }  \pp\l(  \cal{O}_{h}      \r)     
   \sim K^{-\alpha}  \pp(\cal{O}_h),
      \end{align}
   because $\ov$ is \mbox{\hyperref[rv]{regularly varying}} at $\infty$, so    $\l(\ref{finall}\r)$  contradicts $\l(\ref{contr2}\r)$, and therefore   $\liminf_{h\to\infty}     \pp\big(    \Delta_1^{g(h)} \!  \leq \! h  ;     X_{\Delta_1^{g(h)}  -}   \!  < \! g(h)   ; \cal{O}_h    \big)   / \pp(\cal{O}_h ) \!  > \! 0  $. By   $\l(\ref{lowerbound}\r)$, uniformly in $y\!  > \! K$~as~$h\!  \to \! \infty$, 
   \[
 \pp(  X_h \in g(h) dy  ;    \cal{O}_h)  \gtrsim                y^{-1-\alpha}      \f{ L(g(h)y)}{L(g(h))  }  \pp\l(  \cal{O}_h  \r)    dy      ,
    \]
  as required for the lower bound in  $\l(\ref{l6statement}\r)$, so the proof of Lemma \mbox{\ref{l6}}~is~complete.

    \e{proof}

\section{Proofs of Auxiliary Lemmas} \label{lemmasproofs}  

 \b{proof}[Proof of Lemma \mbox{\ref{lowerboundforrho}}] Firstly, we will show that $\rho(t)=o(\ov(g(t)))$ as $t\to\infty$. By $\l(\ref{2}\r)$ with $y=h=0$, it is immediate that $\lim_{s\to\infty}\rho(s)/\ov(g(s))\leq 0$. Now, by $\l(\ref{rho}\r)$ and $\l(\ref{final123}\r)$, 
  \b{align}
  \nonumber
 - \rho(s)  &=      \f{1}{\Phi(s)}    \l[       \ov(g(s))^2   \int_0^s      \int_0^v       \mathbb{P}   (\cal{O}_{w})   e^{-\ov(g(s))w}   dw
dv      - \pp( \oo_{s}  ;  \Delta_1^{g(s)} >s   )        \r] 
\\
\nonumber
&\leq 
      \f{1}{\Phi(s)}        \ov(g(s))^2   \int_0^s      \int_0^v       \mathbb{P}   (\cal{O}_{w})   e^{-\ov(g(s))w}   dw dv    
  \\
  \label{withouttheintegralbit}
&\leq 
     \f{1}{\Phi(s)}    s    \ov(g(s))^2   \int_0^s           \mathbb{P}   (\cal{O}_{w})     dw       =         s    \ov(g(s))^2       .
  \intertext{Now, applying $\l(\ref{condt}\r)$ and $\l(\ref{cond2}\r)$ in cases  (\hbox{\hyperref[case1]{i}}) and (\hbox{\hyperref[case2]{ii}}) respectively, it follows that $\lim_{s\to\infty} -\rho(s)/\ov(g(t)) \leq 0$, and hence $\rho(s)=o(\ov(g(s)))$ as $s\to\infty$.}
  \intertext{Next we will show that $-\int_1^\infty \rho(s)ds<\infty$. Indeed, this follows immediately  in case    (\hbox{\hyperref[case2]{ii}}) by $\l(\ref{withouttheintegralbit}\r)$ and $\l(\ref{cond2}\r)$. In  case (\hbox{\hyperref[case1]{i}}), applying Potter's theorem \mbox{\cite[Theorem 1.5.6]{bgt89}} to $\ov(g(s))/\ov(g(s)/\log(s)^\beta)$, for arbitrarily small $\tau>0$, it follows by $\l(\ref{withouttheintegralbit}\r)$ and $\l(\ref{condt}\r)$ that}
  \nonumber
-\int_1^\infty \rho(s)ds &\lesssim   \int_1^\infty      s^2    \ov\l( \f{g(s)}{\log(s)^\beta}\r)^2 s^{-1} \log(s)^{-2\alpha\beta + 2\tau}    ds  \overset{(\ref{condt})}\lesssim   \int_1^\infty      s^{-1} \log(s)^{-2\alpha\beta + 2\tau}    ds <\infty,
  \end{align}
  where we simply take $\tau$ small enough that $-2\alpha\beta+2\tau <-1$, which is possible since $\alpha\beta>1/2$.

  \end{proof}

%
%
%

  \b{proof}[Proof of Lemma  \hbox{\ref{lemma4}}]    Recall that  $\ov(x) = x^{-\alpha} L(x)$  for $L$  \mbox{\hyperref[sv]{slowly varying}}  at $\infty$, so as  $\ov$ is non-increasing, for~large~$N>0$, using that $t_0(y)\geq f(Ay)$,
  \b{align}
\label{lll}
 &\  \ \ \   \int_{f(Ay)}^t ( \ov(g(s+h) - y) - \ov(g(s))) ds     \leq    \int_{f(Ay)}^t  ( \ov(g(s) - y) - \ov(g(s))) ds
  \\ \nonumber
  &\leq        \int_{f(Ay)}^\infty ( \ov(g(s) - y) - \ov(g(s))) ds
%
%
\hspace{20pt} =  \int_{f(Ay)}^\infty    \l(     \f{  L(g(s) - y)   }{    (g(s) - y)^\alpha }    -     \f{  L(g(s) )   }{    g(s)^\alpha }   \r)   ds 
\\ 
\nonumber
 &= \int_{f(Ay)}^\infty     \f{  L(g(s) - y)  }{    (g(s) - y)^\alpha }        \l( 1     -      \l(\f{g(s)-y}{g(s)}\r)^\alpha  \f{  L(g(s) )  }{    L(g(s) - y) }      \r)      ds
 \\ 
\nonumber
 &= \int_{f(Ay)}^\infty     \f{  L(g(s) - y)  }{    (g(s) - y)^\alpha }        \l( 1     -      \l(\f{g(s)-y}{g(s)}\r)^{\alpha+N}  \f{g(s)^N      L(g(s) )     }{(g(s)-y)^N    L(g(s) - y)   }      \r)      ds.
\intertext{Now, $A>B-1$, and  $x^NL(x)$ is non-decreasing in $x$ for $x>B$ in case (\mbox{\hyperref[case1]{i}}),~so  }
\nonumber
  \l(\ref{lll}\r) &  \leq         
\int_{f(Ay)}^\infty     \f{  L(g(s) - y)  }{    (g(s) - y)^\alpha }        \l( 1     -      \l(\f{g(s)-y}{g(s)}     \r)^{\alpha+N}   \r)   ds. 
\intertext{One can verify   $1-(1-y/g(s))^{\alpha+N} \lesssim y/g(s)$, uniformly in $y  >  0,s  >  f(Ay)$,~so}
\nonumber
    \l(\ref{lll}\r) &\lesssim      \int_{f(Ay)}^\infty     \f{  L(g(s) - y)  }{    (g(s) - y)^\alpha }   \f{y}{g(s)}       ds .
    \intertext{As $g(s)  -  y\hspace{-2pt}\geq \hspace{-2pt}(1  -  A^{-1}) g(s) $ for   $s  >  f(Ay)$, and $\ov(x)\hspace{-3pt} =\hspace{-3pt} x^{-\alpha} L(x)$~is~non-increasing,
    }
  \l(\ref{lll}\r)      &\lesssim    y \int_{f(Ay)}^\infty     \f{  L((1-A^{-1})g(s) )  }{    g(s)^{1+\alpha }}    ds.
\intertext{Applying the   uniform convergence theorem \cite[Theorem 1.2.1]{bgt89} to the slowly varying function $L$, substituting $u=g(s)$, as  $uf'(u)\ov(u)$ is decreasing, we conclude that uniformly in $y>0$ (and so also uniformly in $h>0,y>g(h)$),}
\nonumber
    \l(\ref{lll}\r) &\lesssim 
  y \int_{f(Ay)}^\infty     \f{  L(g(s) )  }{    g(s)^{1+\alpha }}    ds
  =   y \int_{Ay}^\infty    \f{L(u)}{u^{1+\alpha}}   f'(u)du  
   =   y \int_{Ay}^\infty    u^{-2} u  f'(u)  \ov(u)    du    
\\
\nonumber
 & \leq  A  y^2 f'(Ay) \ov(Ay) \int _{Ay}^\infty    u^{-2}   du =\f{1}{2} y f'(Ay) \ov(Ay) \lesssim y f'(y) \ov(y)  .
  \end{align}
  \e{proof}

  \b{proof}[Proof of Lemma \mbox{\ref{qh}}]   
   First recall that by Theorem \mbox{\ref{recurrentthm}},  
 \[
q_h(y)     \hspace{-0.1cm}  =     \hspace{-0.1cm} \f{    \Phi_y^h(t_0(y))}{    \Phi(1)}      \hspace{-0.1cm} \lim_{t\to\infty}     
      \exp\l(    \int_{  t_0(y)}^{t}       \hspace{-0.1cm}    \l(  \ov(g_{y}^{h}(s) )    \hspace{-0.1cm}  +    \hspace{-0.1cm}    \rho_y^h(s)   \r)   ds      
 -    \hspace{-0.1cm}   \int_1^{ t   }    \hspace{-0.1cm}   \l(  \ov(g(s))  \hspace{-0.1cm}  +  \hspace{-0.1cm}  \rho(s)   \r) ds  \hspace{-0.05cm}     \r)  \hspace{-0.1cm}  .
 \]  
  Now,     by $\l(\ref{1}\r)$ in Lemma \mbox{\ref{lemma2}},  uniformly in $h>0,y>g(h)$,
   \[
   \l|  \int_{t_0(y) }^{  \infty  }   \rho_y^h(s) ds   \r|   \lesssim    \int_{t_0(y)}^\infty     \f{1}{   s \log(s)^{1+\eps} }     \l(  1 + \f{1}{ f(y) -h}  \r)  ds   ,
  \]
  and $ 1/ (f(y) -h) \leq 1/f(\d)<\infty$ since $y\geq g(h+f(\d))$,   so the $\rho_y^h$ integral is  bounded uniformly  in $h>0,y>g(h+f(\d))$.  
   By  Remark \mbox{\ref{yh0}},  $\int_1^\infty \rho(s)ds  <  \infty$. 
       For   $y>g(h+f(\d))>\d$,   $f(Ay)>f(A\d)$, so taking $A$ sufficiently large if necessary,  $t_0(y) := f(Ay) \vee f(1  +  2/A)   =   f(Ay)$, then by Lemma \mbox{\ref{lemma4}}, $ \limsup_{t\to\infty}  \int_{f(Ay)}^t    ( \ov(g_{y}^{h}(s))    -   \ov(g(s))) ds   <  \infty$, and so we have  uniformly in $h>0,y>g(h)$,
  \b{equation*} 
   q_h(y)  \lesssim     \Phi_y^h(f(Ay))    \exp\l(       -       \int_{1}^{f(Ay)} \ov(g(s))   ds \r).
  \end{equation*}
  For the converse inequality, as $\ov$ is non-increasing,  for   $y  >  g(h)$ (so~$f(y)  >  h$),
  \begin{align*}
         &  \phantom{\leq} \  \int_{f(Ay)}^t ( \ov(g(s)) - \ov(g(s+h)-y)) ds  
        &&\hspace{-29pt}\leq      \int_{f(Ay)}^t ( \ov(g(s)) - \ov(g(s+h))) ds  
        \\
         &=      \int_{f(Ay)}^t  \ov(g(s)) ds      -        \int_{f(Ay)+h}^{t+h}  \ov(g(s)) ds  
         &&\hspace{-29pt}\leq      \int_{f(Ay)}^{f(Ay)+h}   \ov(g(s)) ds    
         \\
         &\ &&\hspace{-29pt} \leq       h \ov( g( f(Ay))) \leq      h \ov( y)      \leq f(y) \ov(y).
  \end{align*}
   Then as $y>g(h+f(\d))>\d$ and $\lim_{y\to\infty} f(y) \ov(y) =0$ by $\l(\ref{condt}\r)$ (recall $f^{-1}=g$),   we  conclude   
  \b{equation*} 
   q_h(y)  \asymp     \Phi_y^h(f(Ay))    \exp\l(       -       \int_{1}^{f(Ay)} \ov(g(s))   ds \r).
  \end{equation*}
  %
  %
  %
  %
  %
  \end{proof}
  
   \b{proof}[Proof of Lemma \ref{bigjump}]   In case (\hbox{\hyperref[case1a]{ia}}), with $\Pi(dx)=u(x)dx$, 
 %
 $u(x)$  has bounded decrease and bounded increase, and as $\ov$ is \mbox{\hyperref[rv]{regularly varying}} at $\infty$ with index $-\alpha\in(-1,0)$ in case (\hbox{\hyperref[case1]{i}}), it follows that $\ov$ has positive increase and bounded increase (see \cite[p71]{bgt89} for precise definitions of bounded decrease, bounded increase, and positive increase). Thus we can apply  \cite[Prop 2.2.1]{bgt89}, yielding that $xu(x) \asymp \ov(x)$ for all sufficiently large $x$, so
   \[
        \f{ \Pi(g(h) dv)   }{  \ov(g(h))   }    =    \f{ u(g(h)v)     g(h)dv  }{   \ov(g(h))  }   \asymp  \f{  \ov(g(h)v) g(h)   dv}{g(h)v \ov(g(h))    }   =   v^{-1-\alpha}  \f{L(g(h)v)}{L(g(h))}  dv.
   \]
   \end{proof}

  \b{proof}[Proof of Lemma \mbox{\ref{1fyh}}]  For  $t  \geq   f(Ay)$,   $A  >  3\vee  (B-1)$, as $f$ is increasing,
\[
t  
 \geq f(Ay) \geq f(y) \geq f(y)-h .
\]
For $y>0$,       $y>g(h)$, and $s\leq  f(y)-h$, we have $g_{y}^{h}(s) \leq 0$, so  
\[
\pp\Big( \cal{O}_s^{g_{y}^{h}}  \Big) =   \pp\l(   X_u \geq g_{y}^{h}(u),   \forall u\leq s  \r)\geq  \pp\l(   X_u \geq 0, \forall u\leq s  \r) =1,
\]  and  we conclude, as required, that 
\[
\Phi_y^h(t)  =   \int_0^t  \pp\Big( \cal{O}_s^{g_{y}^{h}}  \Big) ds    \geq   \int_0^{ f(y)- h}   \pp\Big( \cal{O}_s^{g_{y}^{h}}  \Big) ds =  f(y)- h.
  \]
  \e{proof}

   \b{proof}[Proof of Lemma  \hbox{\ref{keylemma}}] By Markov's inequality and $\l(\ref{lk}\r)$, with $\lambda = \log(1/H(t)) /  B(t)$,
                         \begin{align}
                       \nonumber
                         & \   \pp\l( X_t^{ (0,A(t))}   > B(t) \r)    &&\hspace{-6.2cm}=  \pp\l(   e^{\lambda X_t^{(0,A(t))}}   \geq  e^{\lambda B(t)}      \r)   
                        \\
                        \nonumber
                             &\leq   \mathbb{E}\l[    e^{\lambda X_t^{(0,A(t))}}   \r]  e^{-\lambda B(t)} &&\hspace{-6.2cm}=   \exp\l(  t \int_0^{A(t)}  \lambda e^{\lambda x} (\ov(x)-\ov(A(t))) dx  \r)   H(t)
                        \\
                        \nonumber
                        &\leq \exp \l(   t \f{\log(1/H(t))}{ B(t)}    e^{ \lambda A(t)} \int_0^{A(t)}  \ov(x) dx         \r)    H(t)
                        \\
                        \label{secondlastinkeylemma}
                        &=  \exp\l(   t \f{\log(1/H(t))}{ B(t)}    H(t)^{ - \f{A(t)}{B(t)}} \int_0^{A(t)}  \ov(x) dx         \r)    H(t).
                        \end{align}
                       Now, by   \cite[Theorem 2.6.1(b)]{bgt89}, which applies as $\ov$ has  \mbox{\hyperref[matus]{lower  index}} $  \beta(\ov)>-1$ in cases (\hbox{\hyperref[case1]{i}}) and (\hbox{\hyperref[case2]{ii}}), there exists  $C>0$    such that for all   $A(t)>1$,
                           \b{equation}
                           \label{lastinkeylemma}
                         \int_0^{A(t)}  \ov(x) dx  \leq    \int_0^1 \ov(x)dx +   C  \ov(A(t))A(t).
                      \e{equation}
                      Now, consider $\mu(\ov):= \liminf_{x\to\infty} \log(\ov(x))/\log(x)$, which by \mbox{\cite[Prop 2.2.5]{bgt89}} satisfies $\mu(\ov) \geq \beta(\ov) >-1$, so  $\liminf_{x\to\infty} \log(\ov(x))/\log(x) >1$, and thus  $\liminf_{x\to\infty} x \ov(x)>0$. So uniformly among $A(t)\!>\!1$, $\int_0^1 \ov(x)dx \! \lesssim \! A(t)\ov(A(t))$, and    $\l(\ref{keyeqn}\r)$ follows from $\l(\ref{secondlastinkeylemma}\r)$ and $\l(\ref{lastinkeylemma}\r)$, as~required.
                        \e{proof}

\section*{Acknowledgements} Thanks to Mladen Savov and Ronnie Loeffen for discussing this work at various stages. Further thanks to an anonymous referee for their feedback on the paper. Funding was provided by the EPSRC.

\bibliographystyle{plain}
\bibliography{llt}

\begin{thebibliography}{10}

\bibitem{aks12}
{Aurzada, F.}, {Kramm, T.}, and {Savov, M.}
\newblock First passage times of {L\'evy} processes over a one-sided moving
  boundary.
\newblock {\em {Markov Process. Related Fields}}, 21(1):1--38, 2015.

\bibitem{b18}
{Barker, A.}
\newblock Fractal-dimensional properties of subordinators.
\newblock {\em {J. Theoret. Probab.}}, 32(3):1202--1219, 2019.

\bibitem{bb11}
{Benjamini, I.} and {Berestycki, N.}
\newblock An integral test for the transience of a {B}rownian path with limited
  local time.
\newblock {\em {Ann. Inst. H. Poincar\'e Probab. Statist.}}, 47(2):539--558,
  2011.

\bibitem{bhp18}
{Berger, Q.}, {den Hollander, F.}, and {Poisat, J.}
\newblock Annealed scaling for a charged polymer in dimensions two and higher.
\newblock {\em Journal of Physics A: Mathematical and Theoretical},
  51(5):054002, 2018.

\bibitem{b92}
{Bertoin, J.}
\newblock An extension of {Pitman's} theorem for spectrally positive {L\'evy}
  processes.
\newblock {\em The Annals of Probability}, pages 1464--1483, 1992.

\bibitem{b93}
{Bertoin, J.}
\newblock Splitting at the infimum and excursions in half-lines for random
  walks and {L\'evy} processes.
\newblock {\em Stochastic processes and their applications}, 47(1):17--35,
  1993.

\bibitem{b98}
{Bertoin, J.}
\newblock {\em {L\'evy} processes}.
\newblock Cambridge university press, 1998.

\bibitem{bgt89}
{Bingham, N.}, {Goldie, C.}, and {Teugels, J.}
\newblock {\em Regular variation}, volume~27.
\newblock Cambridge university press, 1989.

\bibitem{bl16}
{Biskup, M.} and {Louidor, O.}
\newblock Full extremal process, cluster law and freezing for two-dimensional
  discrete {G}aussian free field.
\newblock {\em {Adv. Math.}}, 330:589--687, 2018.

\bibitem{bbc03}
{Bogdan, K.}, {Burdzy, K.}, and {Chen, Z.}
\newblock Censored stable processes.
\newblock {\em {Probab. Theory Related Fields}}, 127(1):89--152, 2003.

\bibitem{b02}
{Bolthausen, E.}
\newblock Part i: Large deviations and interacting random walks.
\newblock {\em Lectures on Probability Theory and Statistics}, pages 1--124,
  2002.

\bibitem{chp12}
{Caravenna, F.}, {den Hollander, F.}, and {P{\'e}tr{\'e}lis, N.}
\newblock Lectures on random polymers.
\newblock In {\em {Clay Math. Proc.}}, volume~15, pages 319--393, 2012.

\bibitem{c18}
{\c{C}etin, U.}
\newblock Path transformations for local times of one-dimensional diffusions.
\newblock {\em {Stoch. Process. Their Appl.}}, 128(10):3439--3465, 2018.

\bibitem{c96}
{Chaumont, L.}
\newblock Conditionings and path decompositions for {L\'e}vy processes.
\newblock {\em {Stochastic Process. Appl.}}, 64(1):39--54, 1996.

\bibitem{dw10}
{Denisov, D.} and {Wachtel, V.}
\newblock Conditional limit theorems for ordered random walks.
\newblock {\em {Electron. J. Probab.}}, 15:292--322, 2010.

\bibitem{dw14}
{Denisov, D.} and {Wachtel, V.}
\newblock Exact asymptotics for the instant of crossing a curve boundary by an
  asymptotically stable random walk.
\newblock {\em arXiv preprint arXiv:1403.5918}, 2014.

\bibitem{dw15}
{Denisov, D.} and {Wachtel, V.}
\newblock Random walks in cones.
\newblock {\em {Ann. Probab.}}, 43(3):992--1044, 2015.

\bibitem{d98}
{Djur{\v{c}}i{\'c}, D.}
\newblock O-regularly varying functions and strong asymptotic equivalence.
\newblock {\em {J. Math. Anal. Appl.}}, 1998.

\bibitem{egv79}
{Embrechts, P.}, {Goldie, C.}, and {Veraverbeke, N.}
\newblock Subexponentiality and infinite divisibility.
\newblock {\em {Probab. Theory Related Fields}}, 49(3):335--347, 1979.

\bibitem{g09}
{Garbit, R.}
\newblock Brownian motion conditioned to stay in a cone.
\newblock {\em {J. Math. Kyoto Univ.}}, 49(3):573--592, 2009.

\bibitem{gk97}
{Golding, I.} and {Kantor, Y.}
\newblock Two-dimensional polymers with random short-range interactions.
\newblock {\em {Phys. Rev. E}}, 56(2):R1318, 1997.

\bibitem{hkw11}
{Hu, Y.}, {Khoshnevisan, D.}, and {Wouts, M.}
\newblock Charged polymers in the attractive regime: a first-order transition
  from {Brownian} scaling to four-point localization.
\newblock {\em {J. Stat. Phys.}}, 144(5):948, 2011.

\bibitem{k06}
{Kallenberg, O.}
\newblock {\em Foundations of modern probability}.
\newblock Springer Science \& Business Media, 2006.

\bibitem{kl13}
{Kerkhoff, U} and {Lerche, H.}
\newblock Boundary crossing distributions of random walks related to the law of
  the iterated logarithm.
\newblock {\em {Statist. Sinica}}, pages 1697--1715, 2013.

\bibitem{ks13}
{Kolb, M.} and {Savov, M.}
\newblock Transience and recurrence of a {B}rownian path with limited local
  time.
\newblock {\em Ann. Probab.}, 44(6):4083--4132, 11 2016.

\bibitem{ks17}
{Kolb, M.} and {Savov, M.}
\newblock Conditional survival distributions of {B}rownian trajectories in a
  one dimensional {P}oissonian environment in the critical case.
\newblock {\em {Electron. J. Probab.}}, 22, 2017.

\bibitem{ks10}
{K{\"o}nig, W.} and {Schmid, P.}
\newblock Random walks conditioned to stay in {W}eyl chambers of type c and d.
\newblock {\em {Electron. Comm. Probab.}}, 15:286--296, 2010.

\bibitem{l13}
{Lerche, H.}
\newblock {\em Boundary crossing of {B}rownian motion: Its relation to the law
  of the iterated logarithm and to sequential analysis}, volume~40.
\newblock Springer Science \& Business Media, 2013.

\bibitem{m16}
{Mallein, B.}
\newblock Asymptotic of the maximal displacement in a branching random walk.
\newblock {\em {Grad. Journ. Math.}}, 1(2):92--104, 2016.

\bibitem{p08}
{Millan, J.}
\newblock On the rate of growth of {L\'evy} processes with no positive jumps
  conditioned to stay positive.
\newblock {\em {Electron. Commun. Probab.}}, 13:494--506, 2008.

\bibitem{n18}
{Nolan, J.}
\newblock {\em Stable Distributions - Models for Heavy Tailed Data}.
\newblock Birkhauser, Boston, 2018.
\newblock In progress, Chapter 1 online at
  http://fs2.american.edu/jpnolan/www/stable/stable.html.

\bibitem{p13}
{Pant{\'\i}, H.}
\newblock On {L\'evy} processes conditioned to avoid zero.
\newblock {\em {Lat. Am. J. Probab. Math. Stat.}}, 14:657--690, 2017.

\bibitem{p75}
{Pitman, J.}
\newblock One-dimensional {B}rownian motion and the three-dimensional {B}essel
  process.
\newblock {\em {Adv. Appl. Probab.}}, 7(3):511--526, 1975.

\bibitem{pw01}
{P{\"o}tzelberger, K.} and {Wang, L.}
\newblock Boundary crossing probability for {B}rownian motion.
\newblock {\em {J. Appl. Probab.}}, 38(1):152--164, 2001.

\bibitem{rvy06}
{Roynette, B.}, {Vallois, P.}, and {Yor, M.}
\newblock Some penalisations of the {W}iener measure.
\newblock {\em {Jpn. J. Math.}}, 1(1):263--290, 2006.

\bibitem{rvy07}
{Roynette, B.}, {Vallois, P.}, and {Yor, M.}
\newblock Penalizing a {BES(d)} process {(0< d< 2)} with a function of its
  local time, {V}.
\newblock {\em {Studia Sci. Math. Hungar.}}, 45(1):67--124, 2007.

\bibitem{sv09}
{Salminen, P.} and {Vallois, P.}
\newblock On subexponentiality of the {L\'evy} measure of the diffusion inverse
  local time; with applications to penalizations.
\newblock {\em {Electron. J. Probab}}, 14:1963--1991, 2009.

\bibitem{hk01a}
{van der Hofstad, R.} and {Klenke, A.}
\newblock Self-attractive random polymers.
\newblock {\em {Ann. Appl. Probab.}}, pages 1079--1115, 2001.

\bibitem{hk01}
{van der Hofstad, R.} and {K{\"o}nig, W.}
\newblock A survey of one-dimensional random polymers.
\newblock {\em {J. Stat. Phys.}}, 103(5):915--944, 2001.

\bibitem{w80}
{Westwater, M.}
\newblock On {E}dwards' model for long polymer chains.
\newblock {\em {Comm. Math. Phys.}}, 72(2):131--174, 1980.

\bibitem{y94}
{Yakymiv, A.}
\newblock On the asymptotics of the density of an infinitely divisible
  distribution at infinity.
\newblock {\em {Theory Probab. Appl.}}, 47(1):114--122, 2003.

\bibitem{yyy09}
{Yano, K.}, {Yano, Y.}, and {Yor, M.}
\newblock Penalising symmetric stable {L\'evy} paths.
\newblock {\em {J. Math. Soc. Japan}}, 61(3):757--798, 2009.

\end{thebibliography}

\end{document}